\newcommand{\pbaddress}{biran@math.tau.ac.il}
\newcommand{\lpaddress}{polterov@math.tau.ac.il}
\newcommand{\dsaddress}{salamon@math.ethz.ch}

\documentclass[12pt]{amsart}
\usepackage{epsfig}
\usepackage{amscd}
\usepackage[mathscr]{eucal}
\usepackage{amssymb}
\usepackage{amsxtra}
\usepackage{amsmath}
\usepackage[all]{xy}

\theoremstyle{plain}
\newtheorem{thm}{Theorem}[subsection]
\newtheorem{cor}[thm]{Corollary}
\newtheorem{lem}[thm]{Lemma}
\newtheorem{prop}[thm]{Proposition}

\theoremstyle{definition}
\newtheorem{dfn}[thm]{Definition}
\newtheorem*{dfnnonum}{Definition}

\theoremstyle{remark}
\newtheorem{rem}[thm]{Remark}
\newtheorem*{remnonum}{Remark}

\newtheorem{ex}[thm]{Example}
\newtheorem*{exnonum}{Example}
\newtheorem*{exsnonum}{Examples}

\newcommand{\R}{{\mathbb{R}}}
\newcommand{\C}{{\mathbb{C}}}     

\newcommand{\T}{{\mathbb{T}}}
\newcommand{\Z}{{\mathbb{Z}}}
\newcommand{\N}{{\mathbb{N}}}
\newcommand{\D}{{\mathbb{D}}}
\newcommand{\Ha}{{\mathcal{H}}}
\newcommand{\Da}{{\mathcal{D}}}
\newcommand{\Ma}{{\mathcal{M}}}
\newcommand{\id}{{\rm id}}        
\newcommand{\tih}{{\widetilde h}}

\newcommand{\one}{{{\mathchoice {\rm 1\mskip-4mu l} {\rm 1\mskip-4mu l}
{\rm 1\mskip-4.5mu l} {\rm 1\mskip-5mu l}}}}
\newcommand{\A}{\boldsymbol{\mathscr{A}}}
\newcommand{\Aa}{{\mathcal A}}
\newcommand{\Dd}{{\mathcal D}}
\newcommand{\Hh}{{\mathcal H}}
\newcommand{\Jj}{{\mathcal J}}
\newcommand{\Kk}{{\mathcal K}}
\newcommand{\Ll}{{\mathcal L}}
\newcommand{\Mm}{{\mathcal M}}
\newcommand{\Pp}{{\mathcal P}}
\newcommand{\Ss}{{\mathcal S}}
\newcommand{\Uu}{{\mathcal U}}
\newcommand{\Xx}{{\mathcal X}}
\newcommand{\G}{{\rm G}}
\newcommand{\Cinf}{C^\infty}
\newcommand{\Vect}{{\rm Vect}}
\newcommand{\Diff}{{\rm Diff}}

\newcommand{\Hom}{{\rm Hom}}
\newcommand{\om}{{\omega}}
\newcommand{\Om}{{\Omega}}
\newcommand{\eps}{{\varepsilon}}
\renewcommand{\phi}{{\varphi}}

\newcommand{\can}{{\rm can}}
\newcommand{\loc}{{\rm loc}}
\newcommand{\IMP}{\Longrightarrow}

\newcommand{\p}{{\partial}}
\newcommand{\inner}[2]{\langle #1, #2\rangle}   
\newcommand{\INNER}[2]{\left\langle #1, #2\right\rangle} 
\def\Nabla#1{\nabla\kern-.5ex{}_{#1}}
\newcommand{\CF}{{\mathrm{CF}}}
\newcommand{\HF}{{\mathrm{HF}}}
\newcommand{\SHi}{\underleftarrow{\mathstrut{\rm SH}}}
\newcommand{\SHd}{\underrightarrow{\mathstrut{\rm SH}}}

\newcommand{\Jreg}{{\mathcal J}_{\rm reg}}
\newcommand{\TJreg}{\widetilde{\mathcal J}_{\rm reg}}
\newcommand{\IFF}{\Longleftrightarrow}

\newcommand{\Qed}{\hfill \qedsymbol \medskip}

\begin{document}

\title{\,\,Propagation in Hamiltonian dynamics\,\,\linebreak and
  relative symplectic homology}


\date{\today}

\author{Paul Biran, Leonid Polterovich, and Dietmar Salamon}

\address{Paul Biran, School of Mathematical Sciences, Tel-Aviv
  University, Ramat-Aviv, Tel-Aviv 69978, Israel} \email{\pbaddress}
\address{Leonid Polterovich, School of Mathematical Sciences, Tel-Aviv
  University, Ramat-Aviv, Tel-Aviv 69978, Israel} \email{\lpaddress}
\address{Dietmar Salamon, Department of Mathematics, ETHZ, 8092
  Z\"{u}rich, Switzerland} \email{\dsaddress}

%
%
% ...

\maketitle

%
%


\section{Introduction} 


\subsection{Stable propagation}\label{sec:stable} 

Let $M := U^*\T^n$ be the open unit cotangent bundle of the
$n$-dimensional Euclidean torus $\T^n = \R^n/\Z^n$. We express the
elements of $M$ in terms of the canonical coordinates
$(q_1,\dots,q_n,p_1,\dots,p_n)$, where $q_i\equiv q_i+1$.  Thus we
identify $M$ with $\T^n\times\D^n$, where $\D^n = \{|p|<1\}$ is the
open unit ball in $\R^n$ and $|v|$ stands for the Euclidean norm of a
vector $v\in\R^n$.  Consider the space $\Hh$ of all smooth compactly
supported functions on $[0,1] \times M$.  Every function $H \in \Hh$
gives rise to the Hamiltonian system on $M$
\begin{equation}\label{eq:ham}
   \dot q = \frac{\partial H }{\partial p }(t,q,p),\qquad
   \dot p = -\frac{\partial H }{ \partial q}(t,q,p).
\end{equation}
The flow $h_t$ which sends any initial condition $x(0)=(q(0),p(0))$ to
the solution $x(t)=(q(t),p(t))$ at time t is called the {\bf
  Hamiltonian flow} (or {\bf isotopy}) generated by $H$, and the
time-one-map $h_1$ is called the {\bf Hamiltonian diffeomorphism}
generated by $H$. The set of all Hamiltonian diffeomorphisms of $M$
form a group denoted by $\Da$.

Let $f_* = \{f_k\}_{k = 1,2,\dots}$ be an infinite sequence of
Hamiltonian diffeomorphisms. We consider $f_*$ as a dynamical system
as follows. Put $f^{(k)} = f_k\cdots f_1$.  This sequence is called
the {\bf evolution} of $f_*$.  The orbit $\{x_k\}_{k\in\N}$ of a point
$x\in M$ is defined by $x_k = f^{(k)}x$.  Classical dynamical systems
(iterations of a single map $f$) correspond to constant sequences $f_k
\equiv f$.  Sequential systems arise naturally as perturbations of the
classical ones.  In a number of interesting situations stability has
been observed in the sense that these perturbations inherit dynamical
properties of the original system (cf.~\cite{PRud}).  In this paper we
present a new stability phenomenon of this type -- stable propagating
behaviour -- which we are going to describe next.

Every Hamiltonian diffeomorphism $h$ has a {\em canonical} lift
$\widetilde h$ to the universal cover $\widetilde M = \R^n \times
\D^n$ of $M$.  Throughout we denote
$$
Q := \left\{(q,p) \in \widetilde{M} \;{\big |} \; |p| < 1,\; \max_i
   |q_i| \le 1/2\right\};
$$
this is a fundamental domain of the covering.

\medskip
\noindent{\bf Definition.} 
{\it A sequential system $f_*$ {\bf propagates to infinity with speed}
  (at least) $c$ if for every vector $v \in \R^n$ with $|v| \leq c$
  there exists a sequence of points $\widetilde x_k \in Q$ such that
  the sequence $(p_k,q_k):={\widetilde f}^{(k)}(\widetilde x_k)$
  satisfies
  $$
  \lim_{k \to \infty} \frac{q_k}{k} = v.
  $$
  In other words, the projection to $\R^n$ of ${\widetilde
    f}^{(k)}(Q)$ covers the Euclidean ball of radius $kc$ up to an
  error which is small with respect to $k$.}

\medskip
\noindent
Let us illustrate this notion in the following simple example.
Consider a Hamiltonian function $H \in \Ha$ which depends on momenta
variables only: $H=H(p)$.  The lift of the corresponding Hamiltonian
flow to the universal cover is given by ${\widetilde h}^t (q,p) =
(q+t\nabla H(p),p)$.  Therefore the projection of ${\widetilde
  h}^k(Q)$ to $\R^n$ coincides up to a bounded error with the set
$kI$, where $I$ denotes the image of the gradient map $p\mapsto\nabla
H(p)$. Assume now that $H(0) > c$ for some $c > 0$. We claim that $I$
contains the Euclidean ball of radius $c$ centered at zero.  Indeed,
for every $v\in\R^n$ with $|v| \leq c$ the function $F(p) := H(p) -
pv$ satisfies $F(0) > c$ and $F\leq c$ near
$\partial\D^n=\mathbb{S}^{n-1}$.  Therefore $F$ attains its maximum at
a point $p_0\in\D^n$.  Hence $\nabla H(p_0) = v$ and the claim
follows.  This shows that in our example we have propagation with
speed $c$.

The group $\Da$ carries a remarkable {\it biinvariant} metric $\rho$
called {\it Hofer's metric} (see~\cite{HoferMetric}).  The
corresponding geometry provides us with a suitable language for the
study of various stable phenomena in Hamiltonian dynamics.  Given a
diffeomorphism $f \in \Da$, write $\rho(\id,f)$ for $\inf (\max F -
\min F)$ where the infimum is taken over all Hamiltonians $F \in \Ha$
generating $f$.  Define Hofer's distance $\rho(f,g)$ between two
elements $f,g \in \Da$ as $\rho(\id,fg^{-1})$.

\medskip
\noindent{\bf Theorem~A (Stable propagation).} 
{\it Let $0<a<c$ and suppose that $h$ is a Hamiltonian diffeomorphism
  generated by a Hamiltonian $H \in \Ha$ such that $H(t,q,0) \ge c$
  for all $t$ and~$q$.  Let $f_*$ be any sequential system such that
  $\rho(f_i,h) < a$ for all $i\in\N$.  Then $f_*$ propagates to
  infinity with speed $c-a$.  }

\medskip
\noindent
Theorem~A will be proved in Section~\ref{sb: sdr} below.
Interestingly enough, such a propagating behaviour may be completely
destroyed by an appropriate arbitrarily $C^{\infty}$-small dissipative
perturbation even in the framework of classical dynamics.  In
Section~\ref{sb: dis} below we elaborate this in the case $n=1$. We
show that every Hamiltonian diffeomorphism $h$ generated by the
Hamiltonian $H=H(p)$ admits an arbitrarily small smooth perturbation
$f$ such that the images $\widetilde{f}^k (Q)$, $k \in \N$, of the set
$Q$ under the iterates of $\widetilde f$ remain in a compact part
of~$\widetilde M$.

Note that Theorem~A 
does not provide any information about propagation of {\it individual}
trajectories on the universal cover (and we doubt that such
information is available at all in this generality).  The situation
improves when one considers sequences $f_* = \{f_i\}$ which roughly
speaking are uniformly distributed with respect to some ``nice"
measure on $\Da$ whose support is close to $h$ in the sense of Hofer's
metric. It turns out that such sequential systems have trajectories
which propagate to infinity with constant velocity.  We refer the
reader to Section~\ref{sb: sdr} for the details.


\subsection{Noncontractible closed orbits}\label{sec:main} 

The main tool for studying stable propagation as described in the
previous section is an existence result for noncontractible periodic
solutions of compactly supported Hamiltonian systems under quite
robust assumptions on the Hamiltonian functions. This result is based
on Floer homology filtered by the symplectic action and its proof is
quite involved.  A solution $x(t) =(q(t),p(t))$ of a Hamiltonian
system generated by a function $H \in \Ha$ is called {\bf periodic} if
$p(1) = p(0)$ and $q(1) = q(0) + e$ for some integer vector $e \in
\Z^n$.  The lattice $\Z^n$ is identified in a natural way with the
fundamental group of $M$, hence we refer to $e$ as the {\bf homotopy
  class} of the solution.  An important quantity associated to a
periodic solution is its {\bf action}
$$
\Aa_H(x) = \int_0^1 \left( H(t,q(t),p(t)) - \sum_{i=1}^n
   p_i(t){\dot q_i}(t) \right)\,dt.
$$
Denote by $Z\subset M$ the zero section $\{p=0\}$.

\medskip
\noindent{\bf Theorem~B.}  
{\it For every compactly supported smooth Hamiltonian function
  $H\in\Ha$ and every $e\in\Z^n$ such that
\begin{equation}\label{eq:main}
     |e| \le c:=\inf_{[0,1]\times Z}H,
\end{equation}
the Hamiltonian system~(\ref{eq:ham}) has a periodic solution $x(t)$
in the homotopy class $e$ with action $\Aa_H(x)\ge c$.}

\medskip
\noindent
The result above is {\it sharp} in the following sense.
Firstly, the inequalities in Theorem~B 
which guarantee the existence of periodic solutions in the class $e$
cannot be improved.  For instance, for every $\eps > 0$ it is easy to
produce a Hamiltonian of the form $H = H(|p|)$ whose restriction to
$[0,1]\times Z$ equals $1-\eps$ and such that all periodic solutions
are contractible (such a function $H(|p|)$ can be obtained from
$1-|p|$ by an appropriate smoothing).  Secondly, the zero section $Z$
in the inequality~(\ref{eq:main}) cannot be replaced by an arbitrary
smooth section.  Consider, for instance an arbitrary $\Cinf$-small
perturbation $S = \{p=u(q)\} $ of $Z$ such that the 1-form $u(q)dq$ on
$\T^n$ is {\bf not} closed (in symplectic terms this means that the
section $S$ is non-Lagrangian; this is only possible when $n\geq 2$).

\medskip
\noindent{\bf Theorem~C.} 
{\it Given any $c>0$, there exists a Hamiltonian function $H \in \Ha$
  such that $H(t,x)\ge c$ for every $t\in[0,1]$ and every $x\in S$ and
  such that every $1$-periodic solution of~(\ref{eq:ham}) is
  contractible.}

\medskip
\noindent
If one takes $S$ as the graph of {\it an exact} 1-form
on $\T^n$, the assertion of Theorem~B  
remains valid with $Z$ replaced by $S$.  Thus the existence mechanism
for noncontractible periodic
solutions described in Theorem~B 
is quite sensitive to the choice of a subset where the Hamiltonian is
large enough. This phenomenon will be studied below in terms of a {\it
  relative symplectic capacity} (see Section~\ref{sec:rel}).

On the other hand, the statement of Theorem~B 
is {\it robust} from the following viewpoint: if the restriction of a
Hamiltonian $H$ to $[0,1]\times Z$ is bigger than a certain positive
number, (most of) the periodic solutions guaranteed by the theorem
persist under $C^0$-perturbations of the Hamiltonian.  This robustness
plays a crucial role in the study of stable propagation.

In Section~\ref{sec:rel} below we prove Theorem~B 
and its generalization where the torus is replaced by a hyperbolic
manifold.  The proof uses Floer homology for the action functional on
the space of noncontractible loops. The main difficulty we have to go
round is as follows. Since we deal with compactly supported
Hamiltonians, noncontractible solutions are ``non-essential" from the
viewpoint of Floer homology -- they may not persist under deformations
of the Hamiltonian.  In brief, the idea of the proof can be described
as follows. One can squeeze the Hamiltonian function $H$ between two
more or less standard functions, $H_- \le H \le H_+$ where $H_-$ and
$H_+$ depend only on $|p|$. The filtered Floer homologies of $H_-$ and
$H_+$ as well as the natural morphism between them can be computed
explicitly.  This morphism turns out to be nontrivial, and since it
factors through the Floer homology of $H$, we get nontrivial
information about the noncontractible periodic solutions
corresponding to $H$.  The calculations are quite involved, hence we
introduce a convenient algebraic tool -- relative symplectic homology
-- which helps us perform them in a more organized way.

Noncontractible orbits of Hamiltonian systems were studied earlier in
a number of interesting situations.  In a beautiful paper~\cite{GL}
Gatien and Lalonde considered the following setting.  Let $L_0$ and
$L_1$ be two disjoint closed Lagrangian submanifolds of a symplectic
manifold. Assume that $H$ is an autonomous Hamiltonian function which
is ``small" on $L_0$ and ``large" on $L_1$. It turns out that under
certain additional assumptions of a topological nature one can prove
the existence of noncontractible periodic orbits for the Hamiltonian
flow generated by $H$.
This result was the starting point for Theorem~B. 
Another important idea of a symplectic capacity which feels the
fundamental group is contained in Schwarz's work~\cite{Sc}. We will
develop it further in Section~\ref{sec:rel} below.  Noncontractible
orbits of autonomous Hamiltonians on cotangent bundles whose levels
are starshaped were considered by a number of authors (see
e.g.~\cite{Cieliebak-1}). In contrast to our case, these closed orbits
are homologically essential in the sense of Floer homology.  Let us
mention finally that the interest in noncontractible periodic orbits
on cotangent bundles comes from classical mechanics where one
considers Hamiltonians $H(t,q,p)$ which are convex with respect to the
momenta variable $p$. In this case the existence of closed orbits in
given homotopy classes can be derived with methods from the classical
calculus of variations.

\medskip
\noindent
{\bf Acknowledgment.} We thank D. Burago and E. Gluskin for useful
discussions.


\section{Stable propagation in sequential Hamiltonian  dynamics} 

In this section we study stable propagation along the lines mentioned
in the Introduction.  The main tool is Theorem~B.


\subsection{Preliminaries on Hamiltonian diffeomorphisms}
\label{sec:ham}

Let $M$ be an open manifold (i.e.  $M$ is connected, noncompact and
has no boundary) and denote by
$$
\pi:\widetilde M\to M
$$
its universal cover.  Let $\Diff_0(M)$ denote the group of
compactly supported diffeomorphisms of $M$ that are isotopic to the
identity by isotopies with compact support in $[0,1]\times M$.  Then
every $h\in\Diff_0(M)$ has a canonical lift $\tih:\widetilde
M\to\widetilde M$ to the universal cover.  To see this, choose a
compactly supported isotopy $[0,1]\to\Diff_0(M):[0,1]\mapsto h_t$ from
$h_0=\id$ to $h_1=h$.  Given $\widetilde x_0\in\widetilde M$, lift the
path $[0,1]\to M:t\mapsto h_t(\pi(\widetilde x_0))$ to a path
$[0,1]\to\widetilde M:t\mapsto\widetilde x(t)$ such that $\widetilde
x(0)=\widetilde x_0$ and define $\tih(\widetilde x_0):=\widetilde
x(1)$.  The following remark shows that this definition is independent
of the choice of the isotopy $t\mapsto h_t$.

\begin{rem} \label{r: cl}
   Let $t\mapsto\phi_t$ and $t\mapsto\psi_t$ be two compactly
   supported isotopies such that $\phi_1=\psi_1=h\in\Diff_0(M)$.
   Given a point $x\in M$, consider the paths $t\mapsto\phi_t(x)$ and
   $t\mapsto\psi_t(x)$ connecting $x$ to $h(x)$. We claim that they
   are homotopic with fixed endpoints. To see this choose a path
   $[0,1]\to M:s\mapsto\gamma(s)$ such that $\gamma(0)=x$ and
   $\gamma(1)$ lies outside the support of the isotopies. Looking at
   $\varphi_t(\gamma(s))$ and $\psi_t(\gamma(s))$ it is easy to see
   that both paths $t\mapsto\phi_t(x)$ and $t\mapsto\psi_t(x)$ are
   homotopic with fixed endpoints to the path $\Gamma$ obtained by
   going first from $x$ to $\gamma(1)$ along $\gamma$ and then from
   $\gamma(1)$ to $h(x)$ along $h(\gamma^{-1})$. Here $\gamma^{-1}$
   stands for the reverse of $\gamma$.
\end{rem}

Now suppose that $M$ is equipped with a symplectic form $\om$ and
denote by $\Dd\subset\Diff_0(M)$ the group of Hamiltonian
diffeomorphisms that are generated by compactly supported Hamiltonian
functions on $[0,1]\times M$.  Given a compact subset $A\subset M$ and
a real number $c$ let us denote by $\Dd_c=\Dd_c(M,A)\subset\Dd$ the
subset of Hamiltonian diffeomorphisms that are generated by compactly
supported Hamiltonian functions $H\in\Cinf_0([0,1]\times M)$ with
$\inf_{[0,1]\times A}H\ge c$.  As before denote by $\rho$ the Hofer
metric on $\Da$.

\begin{prop} \label{P: hof}
   Let $c > a$ be positive numbers, and let $f,g\in\Da$ be Hamiltonian
   diffeomorphisms with $f\in\Da_c$ and $\rho(f,g)<a$. Then
   $g\in\Da_{c-a}$.
\end{prop}

\begin{proof}
   This is an immediate consequence of the following product formula.
   If the functions $\Phi_t$ and $\Psi_t$ generate Hamiltonian flows
   $\phi_t$ and $\psi_t$, respectively, then the product flow
   $\phi_t\circ\psi_t$ is generated by the Hamiltonian $\Phi_t +
   \Psi_t\circ\phi_t^{-1}$.
\end{proof}

In what follows we shall deal also with time periodic Hamiltonians,
namely functions $H$ satisfying $H(t, \cdot)=H(t+1, \cdot)$. It is
useful to think of these Hamiltonians as smooth functions $H:S^1
\times M \to \mathbb{R}$, where we identify $S^1\cong \R/\Z$.  The
following proposition shows that every $h\in\Dd_c$ can be generated by
a periodic Hamiltonian bounded below by $c$ on $S^1\times A$.

\begin{prop}{\bf (Periodic Hamiltonians)}\label{prop:per}
   Let $M$ be an open symplectic manifold and $A\subset M$ be a
   compact subset. Let $h$ be a Hamiltonian diffeomorphism of $M$
   generated by a Hamiltonian $H\in\Cinf_0([0,1]\times M)$ with
   $\inf_{[0,1]\times A}H\ge c$. Then there exists a Hamiltonian
   $\overline H\in\Cinf_0(S^1\times M)$ with $\inf_{S^1\times
     A}\overline H\ge c$ that generates $h$.
\end{prop}

\begin{proof} 
   Let $h_t$ denote the Hamiltonian isotopy generated by $H$ and $f_t$
   denote the Hamiltonian flow generated by a time independent
   compactly supported function $F:M\to\R$.  Let $\tau:[0,1]\to[0,1]$
   be a smooth nondecreasing function which equals $0$ near $t=0$ and
   equals $1$ near $t=1$.  Then the Hamiltonian isotopy
   $$
   \overline h_t := f_{t-\tau(t)} \circ h_{\tau(t)}
   $$
   is generated by the Hamiltonian functions
   $$
   \overline H_t := F + \tau'(t)(H_{\tau(t)}-F)\circ f_{\tau(t)-t}.
   $$
   The function $\overline H_t$ equals $F$ near $t=0$ and $t=1$ and
   hence defines a smooth Hamiltonian on $S^1\times M$.  Moreover,
   $\overline h_1=h_1=h$.  If $F$ is chosen to be equal to $c$ in a
   neighbourhood of $A$, then $f_t$ is equal to the identity on $A$
   and $H_{\tau(t)}-F$ is nonnegative on $A$, and hence
   $\inf_A\overline H_t\ge c$ for every $t$.
\end{proof}

In the remainder of this section we assume that $M=U^*\T^n$ is the
open unit cotangent bundle and $A=Z=\T^n\subset U^*\T^n$ is the zero
section.  Given $e\in\Z^n$, we denote by $T_e$ the deck transformation
$(q,p)\to(q+e,p)$ of the covering. Note that 1-periodic solutions of
the Hamiltonian system associated to $H$ are in one-to-one
correspondence with the fixed points of $h$.  The homotopy class of
the periodic solution $x(t)$ corresponding to the fixed point $x \in
M$ (or, in brief, {\it the homotopy class of the fixed point} $x$) can
be determined as follows. Pick a lift $y$ of $x$.  Then $\tih(y) =
T_e(y)$ for some $e\in\Z^n$ and this $e$ is the homotopy class in
question.
In this terminology Theorem~B 
asserts that every diffeomorphism $h \in\Da_c=\Dd_c(U^*\T^n,\T^n)$ has
a fixed point in the homotopy classes $e$ for every $e$ such that
$|e|\le c$. It turns out that the same result allows us to get
information on the fixed points of the iterates $h^k$ for $k\in\N$.

\begin{prop} \label{P: iter}
   If $h\in\Da_c$ then $h^k\in\Da_{kc}$ for every $c>0$ and every
   $k\in\N$. In particular, $h^k$ has a fixed point in every homotopy
   class $e\in\Z^n$ such that $|e|\le kc$.
\end{prop}

\begin{proof} 
   We have seen in Proposition~\ref{prop:per} that, for every
   $h\in\Da_c$, one can find a function $H:\R\times M\to\R$ which is
   1-periodic in time, satisfies $\inf_{\R\times A}H\geq c$, and
   generates $h$ as the time-$1$-map.  It follows that the $k$th
   iterate $h^k$ is generated by $kH(kt,x)$, hence $h^k\in\Da_{kc}$,
   and hence the result follows from Theorem~B.
\end{proof}


\subsection{Preliminaries from ergodic theory}\label{sec:ergodic}

Let $X$ be a compact metrizable topological space and $f:X\to X$ be a
homeomorphism. Let $\Mm$ denote the set of Borel probability measures
on $X$ and $\Mm(f)\subset\Mm$ denote the subset of $f$-invariant Borel
probability measures. Both sets are convex and compact with respect to
the weak-$*$ topology ($\mu_i\to\mu$ iff $\int u\,d\mu_i\to\int
u\,d\mu$ for every continuous function $u:X\to\R$). An $f$-invariant
measure $\mu\in\Mm(f)$ is called {\bf ergodic} if every Borel set
$\Lambda\subset X$ such that $f(\Lambda)=\Lambda$ has measure
$\mu(\Lambda)\in\{0,1\}$.

\begin{thm}{\bf (Birkhoff's Ergodic Theorem~\cite{CFS})} 
   \label{thm:birkhoff}
   If $\mu\in\Mm(f)$ is an ergodic $f$-in\-va\-ri\-ant Borel
   probability measure and $u\in C^0(X)$, then there exists a Borel
   set $Y\subset X$ such that $\mu(Y)=1$ and
   $$
   \int u\,d\mu = \lim_{N\to\infty}\frac{1}{N}
   \sum_{k=0}^{N-1}u(f^k(x))
   $$
   for every $x\in Y$.
\end{thm}

\begin{rem}\label{rmk:extremal}
   Let $\Xx$ be a Banach space and $\Kk\subset \Xx^*$ be a weak-$*$
   compact convex subset of its dual space.  A point $x^*\in\Kk$ is
   called {\bf extremal} if, for all $y^*,z^*\in\Kk$ and all real
   numbers $t$,
   $$
   x^*=(1-t)y^*+tz^*,\quad 0<t<1 \qquad\IMP\qquad y^*=z^*.
   $$
   The Krein--Milman theorem~\cite[p. 242]{ROYDEN} asserts that
   every weak-$*$ compact convex subset of $\Xx^*$ is equal to the
   weak-$*$ closure of the convex hull of its extremal points.  In
   particular, every compact convex subset $\Kk\subset\R^n$ with
   nonempty interior has at least $n+1$ extremal points.  Namely, if
   there were less than or equal to $n$ extremal points then $\Kk$
   would be contained in an affine subspace of dimension less than $n$
   and so have empty interior.
\end{rem}

\begin{rem}\label{rmk:ergodic}
   Every extremal point of $\Mm(f)$ is ergodic.  To see this suppose
   that $\mu\in\Mm(f)$ is not ergodic.  Then there exists an invariant
   Borel set $\Lambda=f(\Lambda)$ such that $0<\mu(\Lambda)<1$.  Let
   $t:=\mu(\Lambda)$.  Then
   $\mu=t\mu_\Lambda+(1-t)\mu_{X\setminus\Lambda}$ where
   $\mu_\Lambda\in\Mm(f)$ is defined by
   $\mu_\Lambda(A):=\mu(\Lambda\cap A)/\mu(\Lambda)$ and similarly for
   $\mu_{X\setminus\Lambda}$.  Hence $\mu$ is not an extremal point
   of~$\Mm(f)$.
\end{rem}

\begin{rem}\label{rmk:ue}
   A homeomorphism $f:X\to X$ is called {\bf uniquely ergodic} if
   $\Mm(f)$ consists of a single point, i.e. there exists a unique
   $f$-invariant Borel probability measure $\mu$ on $X$.  In this
   case, by Remark~\ref{rmk:ergodic}, the invariant measure $\mu$ is
   necessarily ergodic.
\end{rem}

\begin{ex}\label{ex:irrational}
   Let $\alpha=(\alpha_1,\dots,\alpha_d)\in\R^d$ be a vector with
   rationally independent components. Then the homeomorphism
   $\phi:\T^d\to\T^d$ induced by the map $\R^d\to\R^d:y\mapsto
   y+\alpha$ is uniquely ergodic and the unique $\phi$-invariant
   measure is the Lebesgue measure (cf.~\cite{CFS}).
\end{ex}

\begin{lem}\label{le:ergodic}
   Let $f:X\to X$ be a homeomorphism of a compact metric space and
   $u:X\to\R^n$ be a continuous function.  Consider the linear
   projection
   $$
   \Mm(f)\to\R^n:\mu\mapsto R(\mu) := \int u\,d\mu.
   $$
   Let $r$ be an extremal point of the compact convex set
   $K:=R(\Mm(f))\subset\R^n$.  Then there exists an ergodic
   $f$-invariant Borel probability measure $\mu\in\Mm(f)$ such that $
   r = \int u\,d\mu.$
\end{lem}

\begin{proof}
   The set $\Mm_r(f):=\{\mu\in\Mm(f)\,|\,R(\mu)=r\}$ is a nonempty
   weak-$*$ compact convex subset of the dual space of $C^0(X)$.
   Hence, by the Krein--Milman theorem, it has an extremal point $\mu$
   (Remark~\ref{rmk:extremal}).  We prove that $\mu$ is an extremal
   point of $\Mm(f)$.  To see this, suppose that
   $\mu=(1-t)\mu_0+t\mu_1$ such that $\mu_0,\mu_1\in\Mm(f)$ and
   $0<t<1$.  We must prove that $\mu_0=\mu_1$.  To see this, let
   $r_0:=R(\mu_0)$ and $r_1:=R(\mu_1)$.  Then $r=(1-t)r_0+tr_1$ and
   $r_0,r_1\in K$.  Since $r$ is an extremal point of $K$ it follows
   that $r_0=r_1=r$ and hence $\mu_0,\mu_1\in\Mm_r(f)$.  Since $\mu$
   is an extremal point of $\Mm_r(f)$ this implies that $\mu_0=\mu_1$.
   Thus we have proved that $\mu$ is an extremal point of $\Mm(f)$ and
   hence, by Remark~\ref{rmk:ergodic}, it is an ergodic $f$-invariant
   Borel probability measure.
\end{proof}


\subsection{Stable propagation revisited} \label{sb: sdr} 

Let us now return to the case where $M=U^*\T^n$ is the open unit
cotangent bundle of the $n$-torus and $\Dd$ is the group of compactly
supported Hamiltonian diffeomorphisms of $M$.  For $a > 0$ and
$h\in\Da$ denote by
$$
B(h,a) := \left\{f\in\Dd\,|\,\rho(f,h) < a \right\}
$$
the open ball of radius $a$ in Hofer's metric.
Recall that Theorem~A 
states the following.  Let $c >a >0$ be real numbers, $h \in \Da_c$ be
a Hamiltonian diffeomorphism, and $f_*$ be a sequential system such
that $f_i\in B(h,a)$ for every $i\in\N$.  Then $f_*$ propagates to
infinity with speed $c-a$.  We now prove Theorem~A, assuming
Theorem~B.

\begin{proof}[Proof of Theorem~A] 
   Write $f_i = h\psi_i$ with $\rho(\id,\psi_i)<a$.  Then the
   evolution of $f_*$ can be written as follows:
   $$
   f^{(k)} = (h\psi_kh^{-1})(h^2\psi_{k-1}h^{-2})\cdots
   (h^k\psi_1h^{-k})h^k.
   $$
   Combining the fact that $\rho$ is biinvariant with the triangle
   inequality we get
   $$
   \rho(f^{(k)},h^k) < ka.
   $$
   Now let $w\in\R^n$ be any vector such that $|w|\le c-a$ and
   choose a sequence $v_k\in\Z^n$ such that
   $$
   w = \lim_{k \to \infty}\frac{v_k}{k},\qquad |v_k| < k(c-a).
   $$
   Since $h^k\in\Da_{kc}$ (see Proposition~\ref{P: iter}) we obtain
   from Proposition~\ref{P: hof} that
   $f^{(k)}\in\Da_{k(c-a)}$. Hence, by Theorem~B, 
   there exists a point $\widetilde x_k = (q_k,p_k) \in Q$ so that
   ${\widetilde f}^{(k)}(\widetilde x_k) = (q_k + v_k,p_k)$.  Since
   the sequence $q_k$ is bounded, we have $\lim_{k\to\infty}(q_k +
   v_k)/k = w$, and hence $f_*$ propagates to infinity with speed
   $c-a$.
\end{proof}

Now we turn to the study of the long time behaviour of {\it individual
  trajectories}.  We start with the following general definition.
Consider the trajectory $x_k=f^{(k)}(x)$ of a point $x \in M$ under a
sequential system $f_*$.  Let $\widetilde x\in\widetilde M$ be a lift
of $x$ and consider the lifted trajectory $ \widetilde x_k={\widetilde
  f}^{(k)}(\widetilde x) $ in the universal cover.

\medskip
\noindent{\bf Definition.}
{\it A trajectory $x_k$ of a sequential system with a lift $\widetilde
  x_k=(q_k,p_k)\in\widetilde M$ is said to {\bf propagate with
    velocity vector} $r\in\R^n$ if
  $$
  \lim_{k \to \infty} \frac{q_k}{k} = r.
  $$
  The Euclidean norm of $r$ is called the {\bf speed} of the
  trajectory $x_k$.  }

\medskip
\noindent
As an illustration, let us mention that every $k$-periodic orbit of a
Hamiltonian diffeomorphism $h$ in a nonzero homotopy class $e\in\Z^n$
propagates with velocity vector $e/k$.  Note also that propagation in
the universal cover $\widetilde M$ corresponds to ``rotation" in $M$.
In classical dynamics the velocity vector of a propagating orbit is
called the rotation vector.
  
Consider now the $d$-dimensional torus $\T^d = \R^d/\Z^d$. Fix a
vector $\alpha = (\alpha_1,\dots,\alpha_d)$ whose coordinates are
independent over the rationals. Let $g: \T^d \to \Da$ be a map which
is continuous with respect to the $C^0$-topology on $\Da$.  For a
point $y \in \T^d$ define a sequential system
$$
f_* (y) := \{g(y + k\alpha)\}_{k\in\N}.
$$
Roughly speaking, the sequence $f_*(y)$ is ``random'' -- it is
uniformly distributed on the subset $g(\T^d)$.
The next result is an improvement of Theorem~A 
for (almost all) such sequences.

\begin{thm}\label{T: rot}
   Let $0 <a < c$ be real numbers and $h\in\Da_c$ be a Hamiltonian
   diffeomorphism.  Let $g:\T^d\to\Dd$ be a continuous map whose image
   is contained in $B(h,a)$. Then, for set of Lebesgue measure one of
   points $y \in \mathbb{T}^d$, the sequential system $f_*(y)$ has at
   least $n+1$ trajectories which propagate with speed at least $c-a$.
   
   In fact there exists a compact convex set $K\subset\R^n$ containing
   the closed ball of radius $c-a$ in $\mathbb{R}^n$ centered at the
   origin so that the following holds: for every extremal point $r\in
   K$ there exists a subset $Y_r\subset\T^d$ of Lebesgue measure one
   such that, for every $y\in Y_r$, the system $f_*(y)$ has a
   trajectory which propagates with velocity vector $r$.
\end{thm}

The basic feature of the irrational shift of the torus which is
crucial for our purposes is unique ergodicity (see
Example~\ref{ex:irrational}).  The assertion of Theorem~\ref{T: rot}
continues to hold, and the proof remains the same, if one replaces the
torus $\T^d$ by an arbitrary compact metric space $Y$, the Lebesgue
measure by any Borel probability measure $\sigma$, and the irrational
shift by any $\sigma$-preserving uniquely ergodic homeomorphism of
$Y$.

\begin{proof}[Proof of Theorem~\ref{T: rot}]
   Let $\overline{M}$ denote the closed unit cotangent bundle and
   $\widetilde M$ its universal cover.  Write $Y := \T^d$ and
   $\phi(y):=y+\alpha$ and denote by $\sigma$ the Lebesgue measure
   on~$Y$.  Consider the skew-product map
   $$
   \overline{M} \times Y \to \overline{M} \times Y: (x,y)\mapsto
   S(x,y) := (g(y)(x),\phi(y))
   $$
   and its canonical lift $\widetilde S$ to $\widetilde M\times Y$.
   There exist functions $u,v:M\times Y\to\R^n$ so that
   $$
   \widetilde S(q,p,y) = (q+u(q,p,y),v(q,p,y),\phi(y)).
   $$
   Here we slightly abuse notation and write $u(q,p,y)$ instead of
   $u(x,y)$, whenever $\widetilde x=(q,p)$.  The function $u$ which
   measures displacement along the $q$-plane on $\widetilde M$ will
   play an important role below.  Note that $u$ has compact support in
   $M\times Y$.
   
   Let $\Ma(S)$ denote the set of all $S$-invariant Borel probability
   measures on $\overline{M} \times Y$. This space is convex and
   weak-$*$ compact (see Section~\ref{sec:ergodic}).  For a measure
   $\mu \in \Ma$ define its {\bf rotation vector} $R(\mu)\in\R^n$ by
   $$
   R(\mu) := \int u\,d\mu \in \R^n.
   $$
   The map $R:\Ma(S)\to\R^n$ is affine and continuous.  Hence its
   image $K:=R(\Mm(S))$ is a compact convex subset of $\R^n$.  This
   set has all the properties formulated in the theorem.
   
   We prove that $K$ contains a ball of radius $c-a$.  Pick a vector
   $v\in\R^n$ with $|v|\le c-a$ and a point $y \in Y$.  Then, by
   Theorem~A, there exists a sequence $x_k\in M$ such that
   $$
   \lim_{k\to\infty}\frac{1}{k}\sum_{i=0}^{k-1} u(S^i(x_k,y)) = v.
   $$
   Denote by $\delta_{i,k}$ the Dirac measure on
   $\overline{M}\times Y$ concentrated at $S^i(x_k,y)$ and consider
   the sequence of measures $\mu_k$ on $\overline{M}\times Y$ defined
   by
   $$
   \mu_k := \frac{1}{k} \sum_{i=0}^{k-1} \delta_{i,k}.
   $$
   By Alaoglu's theorem, this sequence has a limit point
   $\mu=\lim_{\nu\to\infty}\mu_{k_\nu}$ with respect to the weak-$*$
   topology. Since
   $$
   \lim_{k\to\infty}\left( \int F\circ S\,d\mu_k - \int
      F\,d\mu_k\right) =
   \lim_{k\to\infty}\frac{F(S^k(x_k,y))-F(x_k,y)}{k} = 0
   $$
   for every continuous function $F:\overline{M}\times Y\to\R$ the
   limit point $\mu$ is $S$-invariant.  Moreover, it satisfies
   $$
   R(\mu) = \int u\,d\mu = \lim_{\nu\to\infty} \int u\,d\mu_{k_\nu}
   = \lim_{\nu\to\infty} \frac{1}{k_\nu}\sum_{i=0}^{k_\nu-1}
   u(S^i(x_{k_\nu},y)) = v.
   $$
   Hence $v\in K$.
   
   We prove that the extremal points of $K$ satisfy the requirements
   of the theorem.  Let $r$ be an extremal point of $K$. Then, by
   Lemma~\ref{le:ergodic}, there exists an ergodic measure
   $\mu_r\in\Mm(S)$ such that $R(\mu_r) = r$.  By Birkhoff's Ergodic
   Theorem~\ref{thm:birkhoff}, there exists a Borel set $Z_r\subset
   \overline{M}\times Y$ such that $\mu_r(Z_r)=1$ and
   $$
   \lim_{N\to\infty}\frac{1}{N}\sum_{i=0}^{N-1} u(S^i(z)) = \int
   u\,d\mu_r = r
   $$
   for every $z\in Z_r$.  Note that $Z_r\subset M\times Y$.  Let
   $Y_r\subset Y$ be the image of $Z_r$ under the obvious projection
   $\overline{M}\times Y\to Y$ and let $\sigma_r$ be the pushforward
   of the measure $\mu_r$ (i.e.
   $\sigma_r(B):=\mu_r(\overline{M}\times B)$ for every Borel set
   $B\subset Y$).  Then $\sigma_r$ is a $\phi$-invariant Borel
   probability measure on $Y$ and so, by unique ergodicity,
   equals~$\sigma$.  Hence
   $$
   \sigma(Y_r) = \sigma_r(Y_r) = \mu_r(\overline{M}\times Y_r) \ge
   \mu_r(Z_r) = 1.
   $$
   Moreover, for every $y\in Y_r$ there exists a point $x\in M$
   such that $(x,y)\in Z_r$.  By definition of the set $Z_r$, the
   trajectory of the point $x$ under the system $f_*(y)$ propagates
   with velocity vector~$r$.
\end{proof}

It is instructive to compare the proof above with Aubry--Mather theory
for Hamiltonian diffeomorphisms of the annulus generated by a
Hamiltonian which is convex with respect to the momentum $p$.  Aubry
and Mather constructed remarkable orbits which lie on special
Cantor-type sets called cantori. These orbits propagate in our sense.
Roughly speaking, these orbits are obtained as the limits of periodic
ones. A crucial point for such a limiting procedure is the fact that,
due to convexity, one can control the relative positions of closed
orbits in phase space. In our setting the Hamiltonians are not assumed
to be convex, and we have no information on the positions of the
closed orbits.  We went round this difficulty by using an idea which
goes back to Mather~\cite{MinMeas} -- to look at limits of invariant
measures
which sit on periodic orbits provided by Theorem~B. 


\subsection{A dissipative counter-example} \label{sb: dis} 

Assume $n=1$ so $M= S^1 \times (-1,1)$ is the annulus.  Let
$H:(-1,1)\to\R$ be any compactly supported function (which depends on
the momenta variable only and has an arbitrarily large value at
$p=0$).  Then the corresponding Hamiltonian diffeomorphism $h$ has the
form
$$
h(q,p) = (q+H'(p),p),
$$
where $H'$ denotes the derivative of $H$. Recall that $Q$ denotes
the fundamental domain $Q:=[-1/2,1/2]\times(-1,1)$ in the universal
cover $\widetilde M=\R\times(-1,1)$.  We present an example of an
arbitrarily small smooth perturbation $f$ of $h$ such that the images
${\widetilde f}^k(Q)$ under the iterates of $\widetilde f$ remain in a
compact part of $\widetilde M$, thus the propagating behaviour
disappears.  Assume without loss of generality that the support of $H
= H(p)$ is contained in the open interval $(-1/2,1/2)$ and let
$\gamma:=\max_p|H'(p)|$.
   
Let $u:[-1,1]\to[-1,1]$ be an orientation preserving diffeomorphism
such that $u(s) > s$ for $-3/4 < s < 3/4$ and $u(s) = s$ for $|s|\geq
3/4$.  Note that $u$ can be chosen arbitrarily close to the identity.
Choose $N\in\N$ such that $u^N(-2/3) > 2/3$.

Define a map $\phi:M \to M$ by $\phi(q,p) := (q,u(p))$ and let $f :=
\phi h$. For a point $x=(q_0,p_0)\in Q$ consider its orbit
$(q_i,p_i):={\widetilde f}^i(x)$. We claim that $|q_i-q_0| \le
(N+1)\gamma$ for all $x \in Q$ and $i\in\Z$.  This universal estimate
means the absence of propagation.

To prove this, note that $h$ preserves the $p$-coordinate, and $\phi$
preserves the $q$-coordinate of each point. Moreover, the sequence
$p_i$ is nondecreasing. If $|p_0|\ge 3/4$ then $(p_i,q_i)=(p_0,q_0)$
for all $i\in\Z$.  If $|p_0|<3/4$ then
$\lim_{i\to\pm\infty}p_i=\pm3/4$.  In this case let $j_0\in\Z$ be the
largest integer such that $p_{j_0} < -2/3$ and $j_1>j_0$ be the
smallest integer such that $p_{j_1} > 2/3$.  Then $|j_1-j_0|\le N+2$
and $q_i=q_{j_0}$ for $i\le j_0$ and $q_i=q_{j_1}$ for $i\ge j_1$.
For every $i\in[j_0,j_1-1]$ we have $|q_{i+1}-q_i| =
|H'(p_i)|\le\gamma$.  Hence, for all $j,j'\in\Z$ such that $j'>j$,
$$
|q_{j'}-q_j| \le \sum_{i=j}^{j'-1}|q_{i+1}-q_i| \le
\sum_{i=j_0}^{j_1-1}|q_{i+1}-q_i| \le (N+1)\gamma.
$$


\section{Relative symplectic topology}\label{sec:rel} 

Below we discuss existence and nonexistence results for
noncontractible closed orbits in a more general topological
framework.  The main player in this section is a relative symplectic
capacity -- a symplectic invariant which provides a convenient
language for thinking about these results. Using this language, we
prove Theorems~B and~C 
as stated in the Introduction.


\subsection{Symplectic action}\label{sec:action}

This is a preparatory section in which we set notation.  Let $(M,\om)$
be an open symplectic manifold.  We assume throughout that the
symplectic form $\om$ is exact and fix a $1$-form $\lambda\in\Om^1(M)$
such that
$$
d\lambda = \om.
$$
Let $S^1:=\R/\Z$ and denote the free loop space of $M$ by
$LM:=\Cinf(S^1,M)$. We denote the set of free homotopy classes of
loops in $M$ by $\widetilde{\pi}_1(M)$ and for $x\in LM$ we write
$[x]\in\widetilde{\pi}_1(M)$ for its free homotopy class. Given a
subset $\alpha\subset\widetilde{\pi}_1(M)$, we write
$$
L_\alpha M:= \left\{x\in LM\,|\,[x]\in\alpha\right\}.
$$
We shall mostly consider single elements
$\alpha\in\widetilde{\pi}_1(M)$, however, for some of our applications
it is useful to consider more general subsets of
$\widetilde{\pi}_1(M)$.  We denote the space of compactly supported
time dependent $1$-periodic Hamiltonian functions on $M$ by
$\Hh(M):=\Cinf_0(S^1\times M)$.  We do not distinguish in notation
between the function $H\in\Hh(M)$ and its lift $H:\R\times M\to\R$
and, for $H\in\Hh(M)$ and $t\in\R$, we define $H_t:M\to\R$ by
$H_t(x):=H(t,x)$. Every Hamiltonian function $H\in\Hh(M)$ determines a
$1$-periodic family of Hamiltonian vector fields
$X_t=X_{t+1}\in\Vect(M,\om)$ via $\iota(X_t)\om=-dH_t$. The space of
$1$-periodic solutions of the corresponding Hamiltonian differential
equation representing a class in the set
$\alpha\subset\widetilde{\pi}_1(M)$ will be denoted by
\begin{equation}\label{eq:P(H)}
   \Pp(H;\alpha) := \left\{x\in L_\alpha M\,|\,
      \dot x(t) = X_t(x(t))
   \right\}.
\end{equation}
The elements of $\Pp(H;\alpha)$ are the critical points of the
symplectic action $\Aa_H:L_\alpha M\to\R$, defined by
\begin{equation}\label{eq:action}
   \Aa_H(x) := \int_0^1\left(H_t(x(t)) - \lambda(\dot x(t))\right)\,dt
\end{equation}
for $x\in L_\alpha M$. The sole purpose of the $1$-form $\lambda$ is
to fix a normalization of the symplectic action (which otherwise is
only well defined up to an additive constant).  Some of the invariants
discussed in this paper do not depend on this normalization and this
will be pointed out at the appropriate places. However, we do not
indicate the dependence of $\Aa_H$ on $\lambda$ in the notation.

\begin{rem}\label{rmk:action}
   The notation in this section differs from the one used in
   Section~\ref{sec:main}, where the Hamiltonian functions are not
   required to be periodic in the $t$-variable.
   Proposition~\ref{prop:per} shows that this does not effect the
   class of Hamiltonian diffeomorphisms to which the theory applies.
   It also does not effect the value of the symplectic action at the
   periodic orbits. To see this, fix a compactly supported, but not
   necessarily periodic, Hamiltonian function $H\in\Cinf_0([0,1]\times
   M)$ and let $[0,1]\times M\to M:(t,x)\mapsto h_t(x)$ be the
   Hamiltonian isotopy generated by $H$. Fix any point $x_0\in M$ and
   consider the path $x:[0,1]\to M$ given by
   $$
   x(t) := h_t(x_0).
   $$
   Let $\gamma:[0,1]\to M$ be a smooth path such that $\gamma(0)$
   lies outside of the support of $H$ and $\gamma(1)=x_0$.
   Differentiating the function $s\mapsto\Aa_H(x_s)$, where
   $x_s(t):=h_t(\gamma(s))$, we find
   $$
   \Aa_H(x) = \int_0^1(h^*\lambda-\lambda)(\dot\gamma(s))\,ds.
   $$
   Thus the action of $x$ depends only on $x_0$, $h$, and
   $\lambda$.
\end{rem}


\subsection{A relative symplectic capacity} 

Let $M$ be an open symplectic manifold with symplectic form
$\om=d\lambda$ and $A\subset M$ be a compact subset.  In this section
we define a relative symplectic capacity
$$
C(M,A):2^{\widetilde{\pi}_1(M)} \times[-\infty,\infty) \to
[0,\infty]
$$
for the pair $(M,A)$.  This capacity has the following feature:
given $\alpha\subset\widetilde{\pi}_1(M)$ and $a\in\R\cup\{-\infty\}$
such that $C(M,A;\alpha,a) < \infty$, every $H\in\Hh(M)$ with
$\inf_{S^1\times A}H > C(M,A;\alpha,a)$ must have a $1$-periodic orbit
representing one of the homotopy classes in $\alpha$ with symplectic
action $\Aa_H(x)\ge a$.  Moreover, $C(M,A;\alpha,a)$ is the optimal
bound for the existence of such an orbit.  More precisely, for $c > 0$
let
$$
\Hh_c(M,A) := \left\{H\in\Hh(M)\,|\, \inf_{S^1\times A}H \ge
   c\right\}.
$$
For a subset $\alpha\subset\widetilde{\pi}_1(M)$ and a number
$a\ge-\infty$ we define the {\bf relative symplectic capacity}
$C(M,A;\alpha,a)\ge 0$ by
$$
C(M,A;\alpha,a) := \inf \left\{c>0 \biggm| \forall H\in\Hh_c(M,A)
   \; \exists x\in\Pp(H;\alpha)\; \textnormal{ such that } \Aa_H(x)\ge
   a \right\}.
$$
We use the convention that $\inf\emptyset=\infty$. The infimum in
the definition of $C(M,A;\alpha,a)$ is always achieved (see
Proposition~\ref{P:inf} below).

In order to relax the notation we identify the single element subset
$\{ \alpha\} \subset \widetilde{\pi}_1(M)$ with $\alpha \in
\widetilde{\pi}_1(M)$. Also, for $a=-\infty$ we abbreviate
$$
C(M,A;\alpha) := C(M,A;\alpha,-\infty) = \inf\left\{c>0 \biggm|
   \Pp(H;\alpha)\ne\emptyset \mbox{ for every } H\in\Hh_c(M,A)
\right\}.
$$
Note that the invariant $C(M,A;\alpha)$ is independent of the
choice of $\lambda$. It turns out that $C(M,A;\alpha)$ is {\em finite}
in some interesting cases. To begin with lt us consider the case
$A=\textnormal{pt}$. Then the invariant $C(M,A; \alpha)$, with
$\alpha=\widetilde{\pi}_1(M)$, is analogous to the Hofer-Zehnder
capacity~\cite{H-Z}. The difference is that Hofer and Zehnder
considered nonnegative and time-independent Hamiltonian functions.
Lalonde~\cite{La} suggested to consider more general subsets $A$.

Here are some examples in which our capacity can be computed. Let $X$
be a closed (i.e. compact without boundary) connected Riemannian
manifold. Denote by
$$
U^*X := \left\{v\in T^*X \bigm| |v| <1 \right\} \subset T^*X
$$
the open unit cotangent bundle. This manifold is equipped with the
canonical symplectic form $\omega_\can=d\lambda_\can$ of $T^*X$.  We
identify $X$ with the zero section of $U^*X$ and
$\widetilde{\pi}_1(X)$ with $\widetilde{\pi}_1(U^*X)$.  Moreover, in
the case of the standard torus $X=\T^n=\R^n/\Z^n$ there is a natural
isomorphism
$$
\widetilde{\pi}_1(U^*\T^n)\cong\widetilde{\pi}_1(\T^n)\cong\Z^n.
$$
We denote the Euclidean norm of an integer vector
$k=(k_1,\dots,k_n)\in\Z^n$ by
$$
|k| := \sqrt{\sum_{j=1}^n k_j^2}.
$$

\begin{thm} \label{T:cotangent} 
   {\bf (i)} Let $\T^n=\R^n/\Z^n$ be the flat equilateral torus with
   the metric induced by the standard metric on $\R^n$. Then, for
   every $k\in\Z^n\cong\widetilde{\pi}_1(U^*\T^n)$ and every $a\in\R$,
   $$
   C(U^*\T^n,\T^n;k,a) =\max\{|k|,a\}.
   $$
   {\bf (ii)} Let $X$ be a closed Riemannian manifold with negative
   sectional curvature. Then, for every
   $\alpha\in\widetilde{\pi}_1(X)\cong\widetilde{\pi}_1(U^*X)$ and
   every $a\in\R$,
   $$
   C(U^*X, X;\alpha,a) =
   \max\{\textnormal{length}(\gamma_\alpha),a\},
   $$
   where $\gamma_\alpha:S^1\to X$ is the (unique up to time shift)
   closed geodesic that represents the homotopy class $\alpha$.
\end{thm}


\subsection{Properties}\label{sec:axioms} 

In this section we establish some basic properties of the relative
capacity and prove Theorems~B and~C, assuming
Theorem~\ref{T:cotangent}~(i).

\begin{prop}{\bf (Monotonicity)} \label{P: mon} 
   Let $A_2 \subset A_1 \subset M_1 \subset M_2$, where $M_1$ is an
   open subset of $M_2$. Let $\iota_*:\widetilde{\pi}_1(
   M_1)\to\widetilde{\pi}_1( M_2)$ denote the map induced by the
   inclusion $\iota:M_1 \subset M_2$.  Then, for
   $\alpha_i\subset\widetilde{\pi}_1(M_i)$ and $a_i\ge-\infty$, we
   have
   $$
   \iota_*^{-1}(\alpha_2)\subset\alpha_1,\;\; a_1\leq a_2
   \qquad\IMP\qquad C(M_1,A_1;\alpha_1,a_1)\le
   C(M_2,A_2;\alpha_2,a_2).
   $$
   In particular, if $\iota_*$ is injective then for every $\alpha
   \subset \widetilde{\pi}_1(M_1)$ and $a_1\leq a_2$ we have $$C(M_1,
   A_1; \alpha, a_1) \leq C(M_2, A_2; \iota_*(\alpha), a_2).$$
\end{prop}

\begin{proof}
   Let $c>C(M_2,A_2;\alpha_2,a_2)$. Then every Hamiltonian function
   $H\in\Hh_c(M_2,A_2)$ has a $1$-periodic orbit representing a class
   in $\alpha_2$ with action at least $a_2$. In particular, for every
   $H\in \Hh_c(M_1,A_1)$ there exists a $1$-periodic orbit $x:S^1\to
   M_1$ such that $\iota_*[x]\in\alpha_2$ (and hence $[x]\in\alpha_1$)
   and $\Aa_H(x)\ge a_2\ge a_1$.  This shows that $c\ge
   C(M_1,A_1;\alpha_1,a_1)$.
\end{proof}

Our next result provides a useful condition which guarantees that the
relative capacity is trivial.

\begin{prop}{\bf (Displacement)}\label{P: disp} 
   Suppose there exists a compactly supported Hamiltonian isotopy
   $f_t:M\to M$, $0\le t\le1$, such that $f_1(A)\cap A = \emptyset$
   and all 1-periodic orbits of $f_t$ are contractible.  Then, for
   every constant $c>0$, there exists a Hamiltonian $H\in\Hh_c(M,A)$
   such that $\Pp(H;\alpha)=\emptyset$ for every nontrivial class
   $\alpha\in\widetilde{\pi}_1(M)$. In particular, $C(M,A;\alpha) =
   \infty$ for every nontrivial class $\alpha$.
\end{prop}

\begin{proof}
   Let $F\in\Hh$ be the Hamiltonian function generating $f_t$ and
   suppose that $U$ is a neighbourhood of $A$ such that $f_1(U)\cap U
   =\emptyset$.  Choose a smooth function $G:M\to\R$ such that
   $$
   {\rm supp}\,G\subset U,\qquad \inf_A G + \inf_{[0,1]\times M} F
   \ge c.
   $$
   Let $g_t$ denote the Hamiltonian flow of $G$ and consider the
   Hamiltonian flow $h_t:=f_tg_t$.  Since $f_1(U)\cap U=\emptyset$,
   the Hamiltonian diffeomorphisms $h_1$ and $f_1$ have the same fixed
   points.  Hence all $1$-periodic orbits of the flow $h_t$ are
   contractible: they have the form $x(t)=x(t+1)=h_t(x_0)=f_t(x_0)$
   for some $x_0\in M\setminus U$.  The same holds for the conjugate
   isotopy $\phi_t:=f_1^{-1}h_tf_1$. Moreover, $\phi_1 = g_1f_1 =
   \psi_1$, where $\psi_t = g_tf_t$.  Hence, by Remark~\ref{r: cl},
   all $1$-periodic orbits of the isotopy $\psi_t$ are contractible as
   well.  Observe now that the flow $\psi_t$ is generated by the
   Hamiltonian functions $\Psi_t := G + F_t\circ g_t^{-1}$ and they
   satisfy $\inf_A\Psi_t \ge c$ for $x \in A$.  By
   Proposition~\ref{prop:per}, there is a periodic Hamiltonian in
   $\Hh_c(M,A)$ with the same time-$1$-map as $\Psi_t$ and, by
   Remark~\ref{r: cl}, the $1$-periodic orbits of this periodic
   Hamiltonian are also contractible.
\end{proof}

Let us now look at the relative capacity from the geometric viewpoint.
The group $\Da$ of all compactly supported Hamiltonian diffeomorphisms
of $M$ carries Hofer's metric (see Section~\ref{sec:stable} above).
Denote by $\Da_c\subset \Da$ the subset consisting of all Hamiltonian
diffeomorphisms generated by functions from $\Hh_c(M,A)$.  Recall that
every $h\in\Dd$ has a canonical lift to the universal cover and hence
the homotopy class (in $\widetilde{\pi}_1(M)$) of a fixed point is
independent of the choice of the compactly supported Hamiltonian
generating $h$ (see Remark~\ref{r: cl}).

\begin{prop}{\bf (Stability)} \label{P: hof1}
   {\bf (i)} Suppose that $C(M,A;\alpha) \le c < \infty$ for some
   nontrivial class $\alpha\in\widetilde{\pi}_1(M)$.  Let
   $f\in\Da_{c+a}$ be a Hamiltonian diffeomorphism, where $a>0$.  Then
   every compactly supported Hamiltonian diffeomorphism $h\in\Dd$ with
   $\rho(f,h) < a$ has a fixed point in the class $\alpha$. In
   particular, $\rho(f,\id)\ge a$.
   
   \smallskip
\noindent{\bf (ii)}
Suppose that
$$
\lim_{\substack{a>0 \\ a\to 0}}C(M,A;0,a)=0
$$
and let $c>0$.  Then $\rho(f,\id)\ge c$ for every $f\in\Dd_c$.
\end{prop}

\begin{proof}
   We prove~(i). Let $f\in\Da_{c+a}$ and $h\in\Dd$ such that
   $\rho(f,h)<a$. Then, by Proposition~\ref{P: hof}, $h\in\Dd_c$ and
   hence there exists a Hamiltonian isotopy $h_t$ with $h_1=h$ whose
   Hamiltonian function belongs to $\Hh_c(M,A)$.  By definition of the
   relative capacity, this flow has a 1-periodic orbit in the class
   $\alpha$. In particular $h\ne\id$.  This proves~(i).
   
   We prove~(ii).  Suppose, by contradiction, that $\rho(f,\id)<c$ for
   some $f\in\Dd_c$.  Then, by Propositions~\ref{P: hof}
   and~\ref{prop:per}, there exists a Hamiltonian function
   $H\in\Hh(M)$ that generates the identity and satisfies
   $\inf_{S^1\times A}H > 0$.  By assumption, there exists a constant
   $a>0$ such that $C(M,A;0,a)<\inf_{S^1\times A}H$.  Hence there
   exists a contractible periodic orbit $x\in\Pp(H;0)$ with action
   $\Aa_H(x)\ge a$.  On the other hand, since $H$ generates the
   identity, every orbit is a contractible periodic orbit with action
   zero.  This contradiction proves~(ii).
\end{proof}

\begin{prop} \label{P:inf}
   Let $\alpha\subset\widetilde{\pi}_1(M)$ and $a\ge-\infty$. Then
   every Hamiltonian $H\in\Hh(M)$ with $\inf_{S^1\times A}H \ge
   C(M,A;\alpha,a)$ must have a $1$-periodic orbit representing one of
   the homotopy classes in $\alpha$ with symplectic action
   $\Aa_H(x)\ge a$.  In other words, the set $ \left\{c>0\,|\,\forall
      H\in\Hh_c(M,A) \; \exists x\in\Pp(H;\alpha)\; \textnormal{ such
        that } \Aa_H(x)\ge a \right\} $ is either empty or has a
   minimum.
\end{prop}

\begin{proof}
   Without loss of generality we may assume that either $0 \not \in
   \alpha$ or $a>0$, since otherwise there is nothing to prove (all
   points $x$ outside the support of $H$ are constant periodic orbits
   in the class $0$ and with action $0$). Choose a compact subset
   $K\subset M$ such that $A \subset \textnormal{Int\,} K$ and
   $S^1\times K \supset \textnormal{supp} H$. Next, let $\sigma:M \to
   \mathbb{R}$ be a Hamiltonian function supported in $K$ and such
   that $\sigma|_A > 0$. Consider the sequence of Hamiltonians
   $H_n:=H+\frac{1}{n} \sigma$. Clearly $\inf_{S^1\times A} H_n >
   C(M,A;\alpha, a)$, hence for every $n$, $H_n$ has a $1$-periodic
   orbit $x_n$ with $[x_n] \in \alpha$ and $\Aa_{H_n}(x_n)\ge a$. Note
   that $x_n \subset K$ for every $n$ due to the assumption that
   either $0\not \in \alpha$ or $a>0$.
   
   Now $H_n$ converges to $H$ in the $\Cinf$ topology and $K$ is
   compact. Hence, replacing $H_n$ by a suitable subsequence if
   necessary, we obtain a sequence of periodic orbits $x_n$, that
   converges to a periodic orbit $x$ of $H$.  If follows that, for $n$
   sufficiently large, the loops $x_n$ all represent the same homotopy
   class.  Hence $[x]=\lim_{n\to\infty}[x_n] \in \alpha$ and $\Aa_H(x)
   = \lim_{n\to \infty}\Aa_{H_n}(x_n)\ge a$.
\end{proof}

\begin{proof}[Proof of Theorem~B]  
   Due to Proposition~\ref{prop:per} and the discussion on homotopy
   classes following its proof we may assume our Hamiltonians to be
   defined on $S^1\times U^*\mathbb{T}^n$ rather than on $[0,1]\times
   U^*\mathbb{T}^n$. Let $H:S^1\times U^*\T^n\to\R$ be a compactly
   supported Hamiltonian function and $k\in\Z^n$ an integer vector
   such that
   $$
   |k|\le c := \inf_{S^1\times\T^n}H.
   $$
   By Theorem~\ref{T:cotangent}~(i), $ C(U^*\T^n,\T^n;k,c) = c.$ It
   follows from the definition of our capacity and
   Proposition~\ref{P:inf} that there exists a periodic solution $x
   \in\Pp(H;k)$ of~(\ref{eq:ham}) such that $x$ represents the
   homotopy class $k$ and $ \Aa_{H}(x)\ge c.$ This proves the
   assertion of Theorem~B for periodic Hamiltonian functions.  The
   assertion in the non-periodic case follows from the periodic case,
   Proposition~\ref{prop:per}, and Remark~\ref{rmk:action}.
\end{proof}

\begin{proof}[Proof of Theorem~C] 
   Let $S\subset U^*\T^n$ be a non-Lagrangian section and $c>0$ be any
   real number.  It is shown in~\cite{PFlo, LSik} that there exists a
   Hamiltonian function $H:U^*\T^n\to\R$ such that the vector field
   $X_H$ is nowhere tangent to $S$:
   $$
   x\in S\qquad\IMP\qquad X_H(x)\notin T_xS.
   $$
   We may assume without loss of generality that $H$ has compact
   support.  Now let $\phi_t:U^*\T^n\to U^*\T^n$ denote the flow
   generated by $X_H$.  Then there exists an $\eps>0$ such that
   $$
   0<t<\eps\qquad\IMP\qquad \phi_t(S)\cap S=\emptyset.
   $$
   If $\delta$ is sufficiently small then the only $1$-periodic
   solutions of the Hamiltonian flow $t\mapsto\phi_{\delta t}$ are
   constant and $\phi_\delta(S)\cap S=\emptyset$.  Hence, by
   Proposition~\ref{P: disp}, there exists, for every $c>0$, a
   Hamiltonian function $F\in\Hh_c(U^*\T^n,S)$ such that
   $\Pp(F;k)=\emptyset$ for every nonzero integer vector $k\in\Z^n$.
\end{proof}


\subsection{Existence of closed orbits on hypersurfaces} \label{sb:hyp}
As a by-product of our study of the relative capacity we obtain
existence of closed orbits on hypersurfaces in various situations.

\begin{thm} \label{T:Hyp1}
   Let $X$ be either $\mathbb{T}^n$ or a closed negatively curved
   manifold and $H:T^*X \to \mathbb{R}$ a proper and bounded below
   Hamiltonian. Suppose that the sublevel set $\{ H < c\}$ contains
   the zero section. Then for every nontrivial homotopy class $\alpha
   \in \widetilde{\pi}_1( T^*X)$ there exists a dense subset
   $S_{\alpha} \subset (c,\infty)$ with the property that for every
   $s\in S_{\alpha}$ the level set $\{H=s\}$ contains a closed orbit
   $x_s$ in the class $\alpha$ and with $\int_{x_s}
   \lambda_{\textnormal{can}} > 0$.
\end{thm}
\begin{proof}
   Fix a nontrivial homotopy class $\alpha \in
   \widetilde{\pi}_1(T^*X)$. To prove the statement of the theorem we
   shall show that for every $b>a>c$ there exists $s \in (a,b)$ such
   that the level set $\{H=s\}$ carries a closed orbit $x$ in the
   class $\alpha$ and with $\int_x \lambda_{\textnormal{can}} > 0$.
   
   Pick a Riemannian metric $g$ on $X$ which in case $X=\mathbb{T}^n$
   is the flat equilateral metric, or in the other case has negative
   curvature. Denote by $U^*X \subset T^*X$ the unit cotangent bundle
   of $(X,g)$ endowed with the canonical symplectic structure.
   
   Let $b>a>c$. Due to a rescaling argument we may assume without loss
   of generality that the sublevel set $\{ H\leq b \}$ is contained in
   $U^*X$. Next, put $C_{\alpha} = C(U^*X, X; \alpha^{-1})$. Note that
   by Theorem~\ref{T:cotangent} the number $C_{\alpha}$ is finite (and
   in fact equals the length of a geodesic in the class
   $\alpha^{-1}$). Now choose a smooth function $\sigma: \mathbb{R}
   \to \mathbb{R}$ with the following properties:
   \begin{enumerate}
     \item[$\bullet$] $\sigma(r) = C_{\alpha}$ for $r\leq a$.
     \item[$\bullet$] $\sigma(r) = 0$ for $r \geq b$.
     \item[$\bullet$] $\sigma'(r) < 0$ for $a<r<b$.
   \end{enumerate}
   Consider the compactly supported Hamiltonian $F:U^*X \to
   \mathbb{R}$ defined by $F=\sigma \circ H$. As $F|_X = C_{\alpha}$
   it follows from the definition of $C_{\alpha}$ that $F$ has a
   $1$-periodic orbit, say $y$, in the class $\alpha^{-1}$.  Note that
   $y$ is nonconstant (because $\alpha$ is nontrivial). Therefore
   $y$ must lie on one of the level sets $\{F=\rho\}$ where $0< \rho <
   C_{\alpha}$. Since $\sigma$ takes the interval $(a,b)$ injectively
   onto $(0,C_{\alpha})$ there exists $s \in (a,b)$ such that
   $\{F=\rho\} = \{ H=s\}$, hence $y$ lies on the level set $\{H=s\}$.
   Applying to $y$ a suitable {\em orientation reversing}
   reparametrization will give us a closed orbit $x_s$ of $H$ in the
   class $\alpha$. This follows easily from the fact that $X_F =
   (\sigma' \circ H) X_H$ and $\sigma'(s) < 0$.
   
   Finally note that by Theorem~\ref{T:cotangent}, the orbit $y$ has
   action $\Aa_F(y) \geq C_{\alpha}$, hence
   $$C_{\alpha} \leq \int_0^1\left(F(y(t)) -
      \lambda_{\textnormal{can}}(\dot y(t))\right)\,dt = \rho +
   \int_{x_s} \lambda_{\textnormal{can}} < C_{\alpha} + \int_{x_s}
   \lambda_{\textnormal{can}}.$$
   Consequently $\int_{x_s}
   \lambda_{\textnormal{can}} > 0$.
\end{proof}

\subsubsection{Symplectically convex boundaries}
Let $(M, \omega)$ be a symplectic manifold and $U \subset M$ a
relatively compact domain with smooth {\bf convex} boundary
$Q=\partial \overline{U}$.  Recall that this means, by definition,
that there exists a Liouville vector field $Y$ (namely $\Ll_Y \omega =
\omega$), defined on a neighbourhood of $Q$ in $M$, such that $Y$
points outside of $U$ along $Q$.  Note that in this case $Q$ is a
hypersurface of contact type since the vector field $Y$ gives rise to
a contact form $\lambda_Q = (\iota(Y)\omega)|_{TQ}$ on $Q$ with
$d\lambda_Q = \omega|_{TQ}$.

Denote by $\mathcal{L}_Q = \ker (\omega|_{TQ}) \subset TQ$ the
characteristic line field of $Q$. The Reeb vector field $R$ of
$\lambda_Q$ is a nonvanishing section of $\Ll_Q$ and so defines an
orientation on $\Ll_Q$.  We call this the {\bf canonical orientation}
of $\Ll_Q$.  It induces an orientation on each leaf of the
characteristic foliation of $Q$ (namely the foliation corresponding to
$\mathcal{L}_Q$).

\begin{cor} \label{C:convex1}
   Let $X$ be either $\mathbb{T}^n$ or a closed negatively curved
   manifold. Let $U \subset T^*X$ be a relatively compact domain
   containing the zero section, and with smooth convex boundary
   $Q=\partial \overline{U}$. Let $\mathcal{L}_Q$ be equipped with its
   canonical orientation. Then for every nontrivial homotopy class
   $\alpha \in \widetilde{\pi}_1(X)$ the characteristic foliation of
   $Q$ has a closed leaf $x \subset Q$ with $j_*[x] = \alpha$, where
   $j_*: \widetilde{\pi}_1(Q) \to \widetilde{\pi}_1(X)$ is the map
   induced by the composition $j: Q \subset T^*X
   \stackrel{\pi}{\longrightarrow} X$.
\end{cor}
\begin{remnonum}
   We were informed by C.~Viterbo that the above corollary should
   follow also from the variational techniques of~\cite{H-V}.
\end{remnonum}

\begin{proof}[Proof of Corollary~\ref{C:convex1}]
   Denote by $\omega$ the (canonical) symplectic form on $T^*X$ and
   let $\lambda_Q$ be the contact form on $Q$ induced by a local
   Liouville vector field $Y$ as described above. Since $U$ has convex
   boundary there exist a neighbourhood $W$ of
   $Q=\partial\overline{U}$ in $M$ and a diffeomorphism
   $\Phi:(-\eps,\eps)\times Q \to W$ such that (see~\cite{MS}):
   \begin{enumerate}
     \item[$\bullet$] $\Phi(0,q) = q$ for every $q \in Q$.
     \item[$\bullet$] $\Phi\big((-\eps,0)\times Q\big) \subset U$, and
      $\Phi\big((0,\eps)\times Q\big) \subset W \setminus
      \overline{U}$.
     \item[$\bullet$] $\Phi^*\omega = d(e^{\tau}\lambda_Q)$, where
      $\tau$ is the coordinate on $(-\eps,\eps)$.
   \end{enumerate}
   Choose a smooth function $h:(-\eps,\eps) \to \mathbb{R}$ with the
   following properties:
   \begin{enumerate}
     \item[$\bullet$] $h(\tau) = 0 $ for $\tau \leq -\eps/2$, and
      $h(\tau) = 1$ for $\tau \ge \eps/2$.
     \item[$\bullet$] $h'(\tau) > 0$ for $-\eps/2<\tau<\eps/2$.
   \end{enumerate}
   Define now $H:W \to \mathbb{R}$ by $H(\Phi(\tau, q)) = h(\tau)$.
   Next, extend $H$ to a smooth function $H:T^*X \to \mathbb{R}$ in
   such a way that:
   \begin{enumerate}
     \item[$\bullet$] $H \equiv 0$ on $U \setminus W$.
     \item[$\bullet$] $H \geq 1$ on $T^*X \setminus (U \cup W)$.
     \item[$\bullet$] $H$ is proper.
   \end{enumerate}
   Pick $a,b$ with $0<a<b<1$. By Theorem~\ref{T:Hyp1}, there exists an
   $s_0 \in(a,b)$ such that the level set $\{H=s_0\}$ carries a closed
   orbit of $H$, say $y$, with $\pi_*[y]=\alpha$.
   
   Since $0<s_0<1$, there exists $\tau_0 \in (-\eps/2,\eps/2)$ such
   that $\Phi^{-1}(\{H=s_0\}) = \tau_0\times Q$.  Thus $\Phi^{-1}(y) =
   (\tau_0, x)$ where $x\subset Q$ is a closed curve. On the other
   hand $\Phi^{-1}(y)$ is a closed orbit for the Hamiltonian
   $h(\tau,q) = h(\tau)$ on $\Big( (-\eps,\eps) \times Q, \,
   d(e^{\tau}\lambda_Q) \Big)$.  A simple computation shows that the
   Hamiltonian vector field $X_h$ of $h$ is given along $\tau_0\times
   Q$ by $X_h|_{\tau_0 \times Q} = h'(\tau_0)e^{-\tau_0} R$, where $R$
   is the Reeb vector field of $\lambda_Q$.  This proves that
   $x\subset Q$ is a closed leaf of the characteristic foliation of
   $Q$. The orientation of $x$ agrees with the one of the
   characteristic foliation because $h'(\tau_0)>0$.
   
   Finally, the map $\Phi(\tau_0, \cdot): Q \to T^*X$ is homotopic to
   the inclusion $Q\subset X$, hence $j_*[x] = \pi_*[y] = \alpha$.
\end{proof}

The main point in the proofs of both Theorem~\ref{T:Hyp1} and
Corollary~\ref{C:convex1} is that when $X$ is either $\mathbb{T}^n$ or
a negatively curved manifold we have $C(U^*X, X; \alpha) < \infty$
{\em for every} $\alpha \in \widetilde{\pi}_1(X)$. This suggests the
following definition.
\begin{dfn}
   Let $X$ be a smooth closed manifold. We say that a nontrivial free
   homotopy class $\alpha \in \widetilde{\pi}_1(X)$ is {\bf
     symplectically essential} if there exists a domain $U \subset
   T^*X$ containing the zero-section such that $C(U,X;\iota_*(\alpha))
   < \infty$. Here $\iota_*:\widetilde{\pi}_1(X) \to
   \widetilde{\pi}_1(U)$ is the map induced by the inclusion $\iota:X
   \to U$ of the zero section into $U$.
\end{dfn}
\begin{exnonum}
   Let $X$ be either $\mathbb{T}^n$ or a closed negatively curved
   manifold. It follows from Theorem~\ref{T:cotangent} that every
   nontrivial homotopy class $\alpha \in \widetilde{\pi}_1(X)$ is
   symplectically essential.
\end{exnonum}

In this language, Theorem~\ref{T:Hyp1} has the following obvious
generalization.
\begin{thm} \label{T:Hyp2}
   Let $X$ be a smooth closed manifold and $H:T^*X \to \mathbb{R}$ be
   a proper and bounded below Hamiltonian. Suppose that the sublevel
   set $\{ H < c\}$ contains the zero section. Then for every
   nontrivial symplectically essential class $\alpha \in
   \widetilde{\pi}_1( X)$ there exists a dense subset $S_{\alpha}
   \subset (c,\infty)$ with the property that for every $s\in
   S_{\alpha}$ the level set $\{H=s\}$ contains a closed orbit whose
   projection to the zero section is in the class $\alpha^{-1}$.
\end{thm}

The proof is an obvious modification of the one of
Theorem~\ref{T:Hyp1}. Similarly, the proof of the next result is
analogous to that of Corollary~\ref{C:convex1}.

\begin{cor} \label{C:convex2}
   Let $X$ be a closed manifold and $U \subset T^*X$ a relatively
   compact domain containing the zero section, and with smooth convex
   boundary $Q=\partial \overline{U}$. Let $\mathcal{L}_Q$ be equipped
   with its canonical orientation. Then for every nontrivial
   symplectically essential homotopy class $\alpha \in
   \widetilde{\pi}_1( X)$ the characteristic foliation of $Q$ has a
   closed leaf $x \subset Q$ with $j_*[x] = \alpha$, where $j_*:
   \widetilde{\pi}_1(Q) \to \widetilde{\pi}_1(X)$ is the map induced
   by the composition $j: Q \subset T^*X
   \stackrel{\pi}{\longrightarrow} X$.
\end{cor}


\subsection{Topological applications} 

\begin{dfnnonum}[Manifolds of type $\mathcal{F}$]
   Let $X$ be a smooth closed manifold. We say that $X$ is of {\bf
     type $\mathcal{F}$} if there exist two {\em nontrivial} free
   homotopy classes $\alpha_1, \alpha_2 \in \widetilde{\pi}_1(X)$ such
   that:
   \begin{enumerate}
     \item[(i)] $\alpha_1$ and $\alpha_2$ are not {\em positively
        proportional}, namely there exist no $k_1, k_2 \in \mathbb{N}$
      with $\alpha_1^{k_1} = \alpha_2^{k_2}$.
     \item[(ii)] $\alpha_1$ and $\alpha_2$ are both symplectically
      essential. It is easy to see that this is equivalent to
      existence of one domain $U \subset T^*X$ containing the zero
      section and such that both $C(U,X;\iota_*(\alpha_1))$ and $C(U,
      X; \iota_*(\alpha_2))$ are finite.  Here $\iota_*$ is the map
      induced by the inclusion of zero section $X$ into $U$.
   \end{enumerate}
\end{dfnnonum}

\begin{exsnonum}
   {\bf (i)} Let $X$ be a closed negatively curved manifold.  Then $X$
   is of type $\mathcal{F}$.  To see this note first that $X$ is not
   simply connected because the universal cover is diffeomorphic to
   Euclidean space.  Now let $\alpha \in\widetilde{\pi}_1(X)$ be a
   nontrivial free homotopy class and denote by $\gamma$ the unique
   (up to parametrization) closed geodesic representing the class
   $\alpha$. We claim that $\alpha$ and $\alpha^{-1}$ are not
   positively proportional. Suppose otherwise that there exist
   positive integers $k$ and $\ell$ such that
   $\alpha^k=\alpha^{-\ell}$. Then the $k$'th and $\ell$'th iterates
   $\gamma^k$ and $\gamma^{-\ell}$ of $\gamma$ and $\gamma^{-1}$,
   respectively, are two different geodesics representing the same
   free homotopy class. This is impossible for negatively curved
   manifolds because every geodesic has Morse index zero.
   
   Now let $\alpha_1, \alpha_2 \in \widetilde{\pi}_1(X)$ be two
   nontrivial positively non-proportional classes. By
   Theorem~\ref{T:cotangent}~(ii), $C(U^*X,X;\alpha_1)$ and
   $C(U^*X,X;\alpha_2)$ are finite, hence $X$ if of type
   $\mathcal{F}$.
   
   \smallskip\noindent{\bf (ii)} It follows from
   Theorem~\ref{T:cotangent}(i) that $\mathbb{T}^n$ is also of type
   $\mathcal{F}$.
\end{exsnonum}

\subsubsection{Applications to Hamiltonian circle actions}

\begin{thm} \label{T:moment}
   Let $(M,\omega)$ be a compact symplectic manifold (possibly with
   boundary) with $\dim_{\mathbb{R}}M \geq 4$ and $L \subset
   \textnormal{Int\,}(M)$ a compact connected Lagrangian submanifold.
   Suppose that $M \setminus L$ admits a Hamiltonian circle action
   with a surjective moment map $\mu: M \setminus L \to [r,R)$ that is
   proper onto its image ($R$ may possibly be $\infty$).  Then $L$
   cannot be of type $\mathcal{F}$.  In particular, $L$ cannot be
   diffeomorphic to $\mathbb{T}^n$, $n\geq 2$, or to any negatively
   curved manifold.
\end{thm}

Note that the Theorem, as stated, is not true when
$\dim_{\mathbb{R}}M=2$. For example take $M=U^*S^1$ and let $L=S^1$ be
the zero section.  Then $U^*S^1\setminus S^1$ is diffeomorphic to
$S^1\times \big((-1,0) \cup (0,1)\big)$ and has a circle action with
moment map $\mu(q,p)=-|p|$ which satisfies all the assumptions of the
theorem.

The crucial difference between dimension two and higher ones is that
in dimension two, $M \setminus L$ might be disconnected whereas in
higher dimension this never happens. Indeed the following
$2$-dimensional version of Theorem~\ref{T:moment} holds:

\begin{thm} \label{T:moment2dim}
   Let $(M,\omega)$ be a compact symplectic $2$-manifold (possibly
   with boundary) and $S^1 \approx L \subset \textnormal{Int\,}(M)$ an
   embedded circle.  Suppose that $M \setminus L$ admits a Hamiltonian
   circle action with a surjective moment map $\mu:M\setminus L \to
   [r,R)$ that is proper onto its image ($R$ may possibly be
   $\infty$).  Then $M\setminus L$ is disconnected.
\end{thm}

Before we turn to the proofs of Theorems~\ref{T:moment}
and~\ref{T:moment2dim} let us remark that in both Theorems it is
impossible to drop the assumption that $\mu$ is proper onto its image
and that its image is a half open interval.  Indeed take
$M=\mathbb{T}^2 \approx S^1 \times S^1$ and $L=\textnormal{pt} \times
S^1$.  Then $M \setminus L$ is diffeomorphic to $S^1\times (0,1)$ and
so is connected.  It has an obvious circle action whose moment map
$\mu: S^1 \times (0,1) \to (0,1)$ is the projection onto the second
factor.  However, there is no Hamiltonian circle action on
$S^1\times(0,1)$ such that the image of the moment map is a half open
interval and the moment map is proper.  One can easily produce higher
dimension examples as well, e.g. multiplying $M$ by another
$\mathbb{T}^2$ factor and $L$ by $S^1$.

\begin{proof}[Proof of Theorem~\ref{T:moment}]
   By Darboux' theorem there exist an open relatively compact
   neighbourhood $U_0$ of $L$ in $T^*L$, an open neighbourhood $W_0$
   of $L$ in $M$ and a symplectomorphism $f:(U_0,
   \omega_{\textnormal{can}}) \to (W_0, \omega)$ taking $L \subset
   U_0$ identically to $L\subset W_0$.
   
   Note that $M\setminus U_0$ is compact and let
   $$
   R_0 := \max_{M\setminus W_0}\mu,\qquad U := f^{-1} (M\setminus
   \{\mu \leq R_0\})\subset U_0.
   $$
   Then $U\setminus L$ carries a circle action as well, with moment
   map $\mu \circ f : U\setminus L \to (R_0, R)$ which is proper onto
   its image. By reducing $U_0$ if necessary we may assume that $U$ is
   connected. Note that $U \setminus L$ is also connected because the
   codimension of $L$ in $U$ is at least two. It follows that any two
   $S^1$-orbits in $U\setminus L$ represent positively proportional
   homotopy classes in $\widetilde{\pi}_1(U)$.
   
   Now let $R_0<R_1<R_2<R$ and choose a smooth function
   $\sigma:\R\to\R$ with the following properties:
   \begin{itemize}
     \item $\sigma(r) = 0$ for every $r \leq R_1$.
     \item $\sigma(r) = c$ for every $r \geq R_2$, where $c>0$ is a
      constant that will be determined later.
     \item $\sigma'(r)>0$ for every $R_1 < r < R_2$.
   \end{itemize}
   Define now a compactly supported Hamiltonian $H:U \to \mathbb{R}$
   by putting $H := \sigma \circ \mu \circ f$ on $U\setminus L$, and
   extending $H$ to be $c$ on $L$. Note that the vector field $X_H$ is
   everywhere tangent to the orbits of the circle action on $U
   \setminus L$ and moreover along $\mu^{-1}(R_1, R_2)$, $X_H$ points
   in the direction induced by the circle action (because $\sigma'>0$
   there). Since any two orbits of the circle action represent
   positively proportional homotopy classes in $\widetilde{\pi}_1(U)$
   it follows that any two nonconstant closed orbits of $X_H$ must
   also be positively proportional in $\widetilde{\pi}_1(U)$.
   
   Suppose now that $L$ is of type $\mathcal{F}$. We shall get a
   contradiction by showing that once the constant $c$ from the
   definition of $\sigma$ is chosen to be large enough, $X_H$ must
   carry two closed orbits whose classes are positively
   non-proportional. Indeed, if $L$ is of type $\mathcal{F}$ then
   there exist two positively non-proportional classes $\alpha_1,
   \alpha_2 \in \widetilde{\pi}_1(L)$ and a domain $U'$ containing the
   zero section such that both capacities $C(U',L;\iota'_*(\alpha_1))$
   and $C(U',L;\iota'_*(\alpha_2))$ are finite, where $\iota'_*$ is
   the map induced by the inclusion of the zero section $L$ into $U'$.
   By reducing the size of the Darboux neighbourhood $W_0$ if
   necessary we may assume that $U \subset U'$.  Denote by $\iota:L
   \to U$ and $j:U\to U'$ the obvious inclusions and by $\pi:U\to L$
   the obvious projection.  Then $\pi_*\iota_*=\id$ and hence
   ${j_*}^{-1}(\iota'_*(\alpha_i)) \subset\pi_*^{-1}(\alpha_i)$ for
   $i=1,2$.  Hence, by monotonicity (Proposition~\ref{P: mon}), we
   have
   \begin{align*}
      & C(U,L;{\pi_*}^{-1}(\alpha_1))\le C(U_1,L;\iota'_*(\alpha_1))
      <\infty \\
      & C(U,L;{\pi_*}^{-1}(\alpha_2))\le
      C(U_2,L;\iota'_*(\alpha_2))<\infty.
   \end{align*} 
   Now choose the constant $c$ in the definition of $\sigma$ to be
   larger than both of $C(U, L;{\pi_*}^{-1}(\alpha_1))$ and
   $C(U,L;{\pi_*}^{-1}(\alpha_2))$. Since $H|_L=c$, the Hamiltonian
   $H$ must have two periodic orbits $x_1,x_2\subset U$ with
   $\pi_*[x_1]=\alpha_1$ and $\pi_*[x_2]=\alpha_2$. As $\alpha_1$ and
   $\alpha_2$ are not positively proportional, neither are the classes
   $[x_1]$ and $[x_2]$.  This contradicts the fact, established above,
   that any two periodic orbits of $X_H$ represent proportional
   homotopy classes in $\widetilde{\pi}_1(U)$.
\end{proof}

\begin{proof}[Proof of Theorem~\ref{T:moment2dim}]
   Let $W_0$ be a Darboux neighbourhood of $L$ in $M$ as in the proof
   of Theorem~\ref{T:moment}. Choose an orientation on $L \approx S^1$
   and denote by $\gamma_0 \in \widetilde{\pi}_1(W_0)$ the homotopy
   class represent by the oriented circle $L$. By
   Theorem~\ref{T:cotangent} both $C(W_0, L; \gamma_0)$ and $C(W_0, L;
   \gamma_0^{-1})$ are finite.
   
   Let $\sigma:\mathbb{R} \to \mathbb{R}$ be the function defined in
   the proof of Theorem~\ref{T:moment}, where the constant $c$ used in
   its definition is now taken to be larger than both of $C(W_0, L;
   \gamma_0)$ and $C(W_0, L; \gamma_0^{-1})$. Define now a compactly
   supported Hamiltonian $H:W_0 \to \mathbb{R}$ by putting $H:= \sigma
   \circ \mu$ on $W_0 \setminus L$ and extending $H$ to be $c$ on $L$.
   By the choice of the constant $c$, $H$ must have two periodic
   orbits $x_1, x_2 \subset W_0 \setminus L$ with $[x_1]=\gamma_0$ and
   $[x_2]=\gamma_0^{-1}$.
   
   As in the proof of Theorem~\ref{T:moment}, every nonconstant
   periodic orbit of $H$ represent a class that is positively
   proportional in $\widetilde{\pi}_1(W_0)$ to a class represented by
   an orbit of the $S^1$ action. Thus both $\gamma_0$ and
   $\gamma_0^{-1}$ are positively proportional to classes represented
   by orbits of the $S^1$ action.
   
   Suppose, by contradiction, that $M \setminus L$ is connected. Then
   all the classes represented by orbits of the $S^1$ action are
   positively proportional in $\widetilde{\pi}_1(M\setminus L)$, in
   particular also in $\widetilde{\pi}_1(M)$. It follows that
   $\gamma_0=[x_1]$ and $\gamma_0^{-1}=[x_2]$ become positively
   proportional when viewed as classes in $\widetilde{\pi}_1(M)$,
   namely the classes $\iota_*(\gamma_0), (\iota_*(\gamma_0))^{-1} \in
   \widetilde{\pi}_1(M)$ are positively proportional, where $\iota:W_0
   \to M$ is the inclusion. Hence there exists $k>0$ such that
   $(\iota_*(\gamma_0))^k$ is the trivial class in
   $\widetilde{\pi}_1(M)$. Passing to homology we get that $k[L]=0 \in
   H_1(M; \mathbb{Z})$, where $[L] \in H_1(M; \mathbb{Z})$ is the
   homology class represented by $L$. But $H_1(M; \mathbb{Z})$ has no
   torsion, hence $[L]=0$. Finally, note that if an {\em embedded}
   circle $L$ in an orientable surface $M$ is zero in homology then
   $M\setminus L$ must be disconnected.  Contradiction.
\end{proof}

\subsubsection{Applications to Stein manifolds}
Let $(W,J)$ be a Stein manifold. Recall that a smooth function
$\varphi: W \to \mathbb{R}$ is called {\bf plurisubharmonic} if the
2-form
$$
\omega_{\varphi} := - d(d\varphi\circ J)
$$
is a $J$-positive symplectic form, i.e.~$\om_\phi(v,Jv)>0$ for
every nonzero tangent vector $v\in TW$. We denote by
$$
g_{\varphi}(\cdot,\cdot):=\omega_{\varphi}(\cdot,J\cdot)
$$
the associated K\"ahler metric.  Let $\varphi: W \to \mathbb{R}$ be
an {\bf exhausting plurisubharmonic function}, namely in addition to
being plurisubharmonic $\varphi$ is also proper, bounded from below
and has no critical points outside some compact subset of $W$.  Let
$$
X_\varphi := \textnormal{grad}_{g_{\varphi}}\varphi
$$
be the gradient vector field of $\varphi$ with respect to the
metric $g_{\varphi}$. Then
$$
\Ll_{X_\phi}\om_\phi = d\iota(X_\phi)\om_\phi =
-d(g_\phi(X_\phi,J\cdot)) = -d(d\phi\circ J) = \om_\phi,
$$
and hence the flow $X_\phi^t:W\to W$ of $X_\phi$ satisfies
$(X_{\varphi}^t)^* \omega_{\varphi} = e^t\omega_{\varphi}$.  Denote by
$\Delta_{\varphi}$ the union of all the stable submanifold of the flow
of $X_{\varphi}$:
\begin{equation} \label{eq:skel}
   \Delta_\varphi := \bigcup_{p\in \textnormal{Crit}(\varphi)}
   W^s_p(X_{\varphi}).
\end{equation}
Note that $\Delta_{\varphi}\subset W$ is the maximal compact invariant
subset for the flow of $X_{\varphi}$. We call $\Delta_{\varphi}$ the
{\bf associated skeleton} of $\varphi$.  When $\varphi$ is Morse, each
stable submanifold in the union~\eqref{eq:skel} is isotropic with
respect to $\omega_{\varphi}$ (see~\cite{E-G}), in particular $\dim
\Delta_{\varphi} \leq \frac{1}{2} \dim_{\mathbb{R}}W$.  We remark
that, even if $\varphi$ is Morse, $\Delta_{\varphi}$ may have a quite
``wild'' structure. However, for a generic exhausting plurisubharmonic
function $\varphi$, $\Delta_{\varphi}$ has the structure of an
isotropic CW-complex (see~\cite{Bi}). In rare situations it may even
happen that some $\varphi$'s (Morse or not) have smooth skeletons.

We now turn to a special class of Stein manifolds, namely those
obtained from removing an ample divisor from a smooth algebraic
variety.  More precisely, let $(M,J)$ be a closed algebraic manifold
and $\Sigma \subset M$ a smooth ample divisor. It is well known that
$(M\setminus \Sigma, J)$ is an affine variety, in particular Stein.
The following theorem deals with topological restrictions on the
possible smooth manifolds that may arise as skeletons of Stein
manifolds of the type just mentioned.

\begin{thm} \label{T:affine}
   Let $(M,J)$ be a closed algebraic manifold and $\Sigma \subset M$ a
   smooth and reduced ample divisor. If the Stein manifold
   $(M\setminus \Sigma, J)$ admits an exhausting plurisubharmonic
   function $\varphi:M\setminus\Sigma \to \mathbb{R}$ whose skeleton
   $\Delta_{\varphi}$ is a smooth connected Lagrangian submanifold,
   then $\Delta_{\varphi}$ cannot be of type $\mathcal{F}$. In
   particular, $\Delta_{\varphi}$ cannot be diffeomorphic to
   $\mathbb{T}^n$, $n\geq 2$, or to any closed negatively curved
   manifold.
\end{thm}

\begin{proof}
   Put $W:=M \setminus \Sigma$ and endow $W$ with the complex
   structure $J$.  The proof has two steps.
   
   \smallskip
   \noindent{\bf Step~1.} 
   {\sl For every exhausting plurisubharmonic function $\varphi:W \to
     \mathbb{R}$ there exists an open relatively compact domain $W_0
     \subset W$ with the following properties:
   \begin{enumerate}
     \item[(i)] $W_0 \supset \Delta_{\varphi}$.
     \item[(ii)] The boundary $P:=\p\overline{W}_0$ is smooth,
      connected, and convex with respect to~$\omega_{\varphi}$.
     \item[(iii)] The leaves of the characteristic foliation (with
      respect to $\omega_{\varphi}$) on $P$ are all orbits of a free
      circle action.
   \end{enumerate}
   }
 
 \smallskip
 \noindent
 The idea of the proof is the following. Since $W$ is the complement
 of smooth ample divisor it is possible to endow $W$ with an
 exhausting plurisubharmonic function $\varphi_0$ for which the above
 statement holds. In order to pass from $\varphi_0$ to any given
 plurisubharmonic function $\varphi$ we modify $\varphi_0$ and
 $\varphi$ at infinity to obtain new plurisubharmonic functions
 $\overline{\varphi}_0$ and $\overline{\varphi}$ for which the vector
 fields $X_{\overline{\varphi}_0}$ and $X_{\overline{\varphi}}$ are
 complete. Then by the theory of Eliashberg-Gromov~\cite{E-G} the
 symplectic forms $\omega_{\overline{\varphi}_0}$ and
 $\omega_{\overline{\varphi}}$ are diffeomorphic. Hence the statement
 being true for $\varphi_0$ is true also for $\varphi$.
 
 Here are the precise details. We first adjust $\varphi$ at infinity
 so that the vector field $X_{\varphi}$ becomes complete. This is a
 standard procedure (\cite{E-G}, see also~\cite{Bi-Ci} Lemma~3.1).
 More precisely, let $\overline{\varphi}: W \to \mathbb{R}$ be an
 exhausting plurisubharmonic function with the following properties:
   \begin{enumerate}
     \item[$\bullet$] $\overline{\varphi} = \varphi$ on a relatively
      compact domain $W' \subset W$ that contains $\Delta_{\varphi}$.
     \item[$\bullet$] $\Delta_{\overline{\varphi}} =
      \Delta_{\varphi}$.
     \item[$\bullet$] The vector field $X_{\overline{\varphi}}$ is
      complete.
   \end{enumerate}
  
   Next we endow $W$ with another plurisubharmonic function
   $\varphi_0$ that arises from $W$ being the complement of an ample
   divisor. For this purpose, denote by
   $\mathcal{L}=\mathcal{O}_M(\Sigma) \to M$ the holomorphic line
   bundle defined by $\Sigma$, and let $s:M \to \mathcal{L}$ be a
   holomorphic section with $\Sigma=\{s=0\}$. Since $\Sigma$ is ample
   there exists a hermitian metric $\| \cdot \|$ on $\mathcal{L}$ for
   which the associated metric connection $\nabla$ has positive
   curvature $R^{\nabla}$. Then the (real) $2$-form $\omega=\frac{1}{2
     \pi i} R^{\nabla}$ is a $J$-compatible symplectic form.  Let
   $\varphi_0:W \to \mathbb{R}$ be the function defined by $\varphi_0
   = - \frac{1}{4\pi} \log \|s\|^2$. Then on $W$ we have
   $\omega=-d(d\varphi_0\circ J) = \omega_{\varphi_0}$, hence
   $\varphi_0$ is plurisubharmonic. Note that $\varphi_0$ is
   exhausting. Indeed $\varphi_0$ is proper, bounded below, and, since
   $\Sigma$ is smooth and reduced, $\varphi_0$ has no critical points
   near $\Sigma$.
   
   Let $E(\rho) := \{ v\in \mathcal{L}|_{\Sigma} \bigm| \|v\|<\rho\}$
   and $P(\rho) := \{ v\in \mathcal{L}|_{\Sigma} \bigm| \|v\|=\rho\}$
   be the radius-$\rho$ disc and circle subbundles of
   $\mathcal{L}|_{\Sigma}$ and denote by $\pi:E(\rho) \to \Sigma$ the
   projection. Pick a connection $1$-form $\gamma$ on $P(1)$ such that
   $d\gamma = - \pi^*(\omega|_{\Sigma})$ and consider the symplectic
   form $\omega_0 = \pi^*(\omega|_{\Sigma}) + d(r^2 \gamma)$ where $r$
   is the radial coordinate on the fibres induced by the hermitian
   metric. A simple computation shows that the vector field
   $Z:=\frac{r^2-1}{2r} \frac{\partial}{\partial r}$ defined on $E(1)
   \setminus \Sigma$ satisfies $\Ll_Z \omega_0 = \omega_0$. Moreover
   for every $0<\rho<1$, $Z$ is transverse to $P(\rho)$ and points
   towards the inside of $E(\rho)$. Finally note that the leaves of
   characteristic foliation on $P(\rho)$ are all orbits of a free
   circle action, in fact they all coincide with the fibres of the
   principal circle bundle $P(\rho) \to \Sigma$.  Since $\Sigma$ is
   connected so is $P(\rho)$.
   
   By the symplectic tubular neighbourhood theorem there exist a
   neighbourhood $B$ of $\Sigma$ in $M$, an $\eps>0$, and a
   symplectomorphism $f:(E(\eps),\omega_0) \to (B, \omega)$ that sends
   $\Sigma \subset M$ identically onto $\Sigma \subset E(\eps)$. Pick
   any $0< \eps_0 < \eps$ and put $B_0 := f(E(\eps_0))$ and
   $\mathcal{U}_0 := M \setminus \overline{B}_0$.  Then $\mathcal{U}_0
   \subset W$ is a relatively compact domain with convex boundary
   $\partial \overline{\mathcal{U}}_0$ on which the characteristic
   foliation coincides with the orbits of a free circle action. Taking
   $\eps_0$ to be smaller if necessary we may assume that
   $\mathcal{U}_0$ contains $\Delta_{\varphi_0}$.
   
   Similarly to $\varphi$ we adjust $\varphi_0$ at infinity so that
   the vector field $X_{\varphi_0}$ becomes complete. More precisely,
   let $\overline{\varphi}_0: W \to \mathbb{R}$ be an exhausting
   plurisubharmonic function with the following properties:
   \begin{enumerate}
     \item[$\bullet$] $\overline{\varphi}_0 = \varphi_0$ on an open
      subset containing $\overline{\mathcal{U}}_0$.
     \item[$\bullet$] $\Delta_{\overline{\varphi}_0} =
      \Delta_{\varphi_0}$.
     \item[$\bullet$] The vector field $X_{\overline{\varphi}_0}$ is
      complete.
   \end{enumerate}  
   By~\cite{E-G} there exists a symplectomorphism $F:(W,
   \omega_{\overline{\varphi}_0}) \to (W,
   \omega_{\overline{\varphi}})$. Denote by
   $X_{\overline{\varphi}_0}^t$ the flow of
   $X_{\overline{\varphi}_0}$. Then we have $\bigcup_{t\geq 0}
   X_{\overline{\varphi}_0}^t(\mathcal{U}_0)=W$, hence for $t_0>0$
   large enough $F(X_{\overline{\varphi}_0}^{t_0} (\mathcal{U}_0))
   \supset \Delta_{\varphi}$. Clearly, the domain
   $F(X_{\overline{\varphi}_0}^{t_0} (\mathcal{U}_0)) \subset (W,
   \omega_{\overline{\varphi}})$ has a smooth convex connected
   boundary whose characteristic foliation has leaves which are orbits
   of a free circle action. Since $
   (X_{\overline{\varphi}}^{-t})^*\omega_{\overline{\varphi}} = e^{-t}
   \omega_{\overline{\varphi}} $ the same holds also for each of the
   domains $X_{\overline{\varphi}}^{-t} \circ F \circ
   X_{\overline{\varphi}_0}^{t_0} (\mathcal{U}_0)$. Pick now $t_1>0$
   large enough so that $W_0:=X_{\overline{\varphi}}^{-t_1} \circ F
   \circ X_{\overline{\varphi}_0}^{t_0} (\mathcal{U}_0) \subset W'$.
   Recall that on $W'$ we have $\omega_{\overline{\varphi}} =
   \omega_{\varphi}$. Hence the domain $W_0$ satisfies all three
   conditions claimed in Step~1.
   
   \smallskip
   \noindent{\bf Step~2.} 
   {\sl We prove the theorem.}
   
   \smallskip
   \noindent
   Let $\varphi:W \to \mathbb{R}$ be an exhausting plurisubharmonic
   function with skeleton $\Delta_{\varphi}$ which is a connected
   smooth manifold of type $\mathcal{F}$.
   
   By Darboux' theorem there exist a neighbourhood
   $V(\Delta_{\varphi})$ of $\Delta_{\varphi}$ and a symplectic
   embedding $g:(V(\Delta_{\varphi}), \omega_{\varphi}) \to
   (T^*\Delta_{\varphi}, \omega_{\textnormal{can}})$ which takes
   $\Delta_{\varphi} \subset V(\Delta_{\varphi})$ identically onto the
   zero section $\Delta_{\varphi} \subset T^*\Delta_{\varphi}$.
     
   Let $X_{\varphi} := \textnormal{grad}_{g_{\varphi}}\varphi$, denote
   by $X_{\varphi}^t$ the flow of $X_{\varphi}$, and recall that
   $(X_{\varphi}^t)^* \omega_{\varphi} = e^t\omega_{\varphi}$.  Let
   $W_0\subset W$ be the domain defined by Step~1 and let $W_1\subset
   W$ be a relatively compact domain containing $\overline{W}_0$. From
   the definition of $\Delta_{\varphi}$ it follows that for $T>0$
   large enough we have $\Delta_{\varphi} \subset
   X_{\varphi}^{-T}(W_1) \subset V(\Delta_{\varphi})$.
   
   Denote by $Y=p \frac{\partial}{\partial p}$ the standard Liouville
   vector field on $T^* \Delta_{\varphi}$. Then its flow $Y^t$
   satisfies $(Y^t)^* \omega_{\textnormal{can}} = e^t
   \omega_{\textnormal{can}}$.  Hence $Y^T \circ g \circ
   X_{\varphi}^{-T} : (W_1, \omega_{\varphi}) \to
   (T^*\Delta_{\varphi}, \omega_{\textnormal{can}})$ is a symplectic
   embedding. It follows that $ U :=Y^T \circ g \circ
   X_{\varphi}^{-T}(W_0) \subset T^* \Delta_{\varphi} $ is an open
   relatively compact domain containing the zero section and with
   convex boundary $ Q := Y^T \circ g \circ X_{\varphi}^{-T}(P).  $
   Moreover, by condition~(iii) in Step~1, the leaves of the
   characteristic foliation on $Q$ coincide with the orbits of a free
   circle action on $Q$.  Since $Q$ is connected, it follows that any
   two leaves of the characteristic foliation on $Q$ represent
   positively proportional homotopy classes.  On the other hand,
   $\Delta_{\varphi}$ is of type $\mathcal{F}$ and so, by
   Corollary~\ref{C:convex2}, the characteristic foliation of $Q$
   contains two closed leaves whose homotopy classes are positively
   non-proportional. Contradiction.
\end{proof}


\section{A homological capacity}\label{s:Homological} 

In order to study and compute the relative capacity $C(M,A;\alpha)$ we
shall define another quantity $ \widehat{C}(M,A;\alpha) \ge
C(M,A;\alpha) $ which captures the existence of {\em homologically
  essential} periodic orbits in given homotopy classes.  Here the term
``homologically essential'' refers to Floer homology.  The homological
capacity $\widehat{C}(M,A;\alpha)$ will be defined in purely
Floer-homological terms.  It is easier to compute and enjoys some nice
functorial properties.  We begin with a brief discussion of convex
boundaries.


\subsection{Convex boundaries}\label{sb:convex} 

Let $(\overline M,\omega)$ be a compact connected symplectic manifold
with convex boundary and denote $M:=\overline M\setminus\p\overline
M$.  Recall (\cite{E-G}, see also Section~\ref{sb:hyp} above) that the
boundary is called {\bf convex} if there exist a vector field
$X\in\Vect(\overline M)$ and a neighbourhood $U$ of $\p\overline M$
such that $X$ points out on the boundary and is dilating on $U$,
namely $\Ll_X\om=\om$ on $U$. Let $\phi_t$ denote the flow of $X$,
suppose that $ U = \left\{\phi_t(x)\,|\,x\in\p
   M,\,-\eps<t\le0\right\}, $ and denote by $ \xi :=
\ker(\iota(X)\omega|_{_{T\partial\overline{M}}}) $ the contact
structure on the boundary determined by $X$ and $\om$.  Under these
hypotheses (the existence of $X$ and $U$) there is an $\om$-compatible
almost complex structure $J$ on $\overline M$ such that
\begin{equation} \label{eq:Jconvex}
   J\xi=\xi, \qquad \om(X(x),J(x)X(x)) = 1,\qquad 
   D\phi_t(x)J(x) = J(\phi_t(x))D\phi_t(x)
   \notag
\end{equation}
for all $x\in\p\overline M$ and $t\in(-\eps,0]$.  Such an almost
complex structure is called {\bf convex near the boundary}.

Consider the function $f:U\to\R$ given by
$$
f(\phi_t(x)):=e^t
$$
for $x\in\p\overline M$ and $-\eps<t\le 0$. A simple computation
shows that its gradient with respect to the metric $\om(\cdot,J\cdot)$
is $X$ and hence $-d(df \circ J)=\omega$.  This means that $f$ is
plurisubharmonic on $U$, or in other words the $2$-form $-d(df\circ
J)$ is positive on every $J$-complex tangent line in $TU$.  Let
$u:\Om\to U$ be a nonconstant $J$-holomorphic curve, defined on a
connected open subset $\Om\subset\C$. Then the function $f\circ
u:\Omega \to \mathbb{R}$ is subharmonic and, by the mean value
inequality, cannot have a strict interior maximum.  Hence a
nonconstant $J$-holomorphic curve in $\overline M$ cannot intersect
$\p\overline M$ at an interior point of its domain $\Omega$.

\begin{rem}\rm
   Since $J$ is invariant under the flow of $X$ (near $\p\overline M$)
   we have $ 0 = (\Ll_XJ)v = (\Nabla{X}J)v + J\Nabla{v}X -
   \Nabla{Jv}X, $ where $\nabla$ denotes the Levi-Civita connection of
   the metric $\om(\cdot,J\cdot)$.  Now consider the $1$-form $
   \alpha:=-df\circ J =\iota(X)\om=\inner{JX}{\cdot}.  $ A simple
   calculation, using $\Ll_XJ=0$, shows that
   $$
   d\alpha(v,w) = \inner{\Nabla{v}(JX)}{w} -
   \inner{\Nabla{w}(JX)}{v} = 2\inner{J\Nabla{v}X}{w}.
   $$
   Hence the identity $d\alpha=\om$ is equivalent to
   $\Nabla{v}X=v/2$.  Since $X$ is the gradient of $f$ it follows that
   the Laplacian of $f\circ u$ is given by $\Delta(f\circ u)=|du|^2/2$
   for every $J$-holomorphic curve $u:\Omega\to \overline M$.
\end{rem}

\begin{rem}\label{rmk:convex}\rm
   The space of almost complex structures on $\overline{M}$ that are
   convex near the boundary is connected.  To see this fix first the
   dilating vector field $X$ on a neighbourhood of $\partial
   \overline{M}$. Then the space of $\omega$-compatible almost complex
   structures on the symplectic bundle $\xi =
   \ker(\iota(X)\omega|_{_{T\partial \overline{M}}})$ is connected
   (see~\cite{MS}), hence the space of $\omega$-compatible almost
   complex structures $J$ satisfying~\eqref{eq:Jconvex} is also
   connected.  Finally, we may allow $X$ to vary since the space of
   vector fields that are dilating near $\partial \overline{M}$ and
   point out on $\partial \overline{M}$ is convex hence connected.
\end{rem}


\subsection{The setting}\label{sb:setting} 

{From} now on our standing hypotheses are that $(\overline M,\om)$ is
a compact connected symplectic manifold with convex boundary
$\p\overline M$ and $A\subset M := \overline M \setminus \p\overline
M$ is a compact subset. We assume that the symplectic form is exact
and fix a $1$-form $\lambda\in\Om^1(\overline M)$ such that $ d\lambda
= \om.  $ We call $\lambda$ an {\bf $\om$-primitive}.

As in section~\ref{sec:ham}, we denote by
$\Hh:=\Hh(M):=\Cinf_0(S^1\times M)$ the space of smooth compactly
supported smooth Hamiltonian functions on $S^1\times M$ and by
$\Dd\subset\Diff_0(M)$ the group of Hamiltonian diffeomorphisms of $M$
generated by functions from $\Hh$.  For $c>0$ we denote by
$\Hh_c=\Hh_c(M,A)$ the subspace of all Hamiltonian functions $H\in\Hh$
that satisfy $\inf_{S^1\times A} h \ge c$ and by $\Dd_c = \Dd_c(M,A)$
the set of all Hamiltonian diffeomorphisms that are generated by
functions from $\Hh_c$.


\subsection{Floer homology}\label{sb:FH} 

Floer homology is an essential ingredient in the definition of our
invariants. The purpose of this section is to summarize the main
building blocks of this theory needed for our applications. The reader
is referred to~\cite{F, F-H-2, C-F-H-1, Vi-1} for a detailed
foundation of the subject (see also~\cite{Sa} for a general
exposition).

Fix a {\em nontrivial} free homotopy class
$\alpha\in\widetilde{\pi}_1(M)$ and recall from Section~\ref{sec:ham}
that $\Pp(H;\alpha)\subset L_\alpha M$ denotes the set of periodic
solutions of the Hamiltonian system associated to $H\in\Hh$ and that
these periodic solutions are the critical points of the symplectic
action functional $\Aa_H:L_\alpha M\to\R$ defined
by~(\ref{eq:action}).  The set of critical values of $\Aa_H$ is called
the {\bf action spectrum} and will be denoted by
$$
\Ss(H;\alpha) := \Aa_H(\Pp(H;\alpha)) = \left\{\Aa_H(x) \,|\, x\in
   L_\alpha M,\, \dot x(t) = X_{H_t}(x(t))\right\}.
$$
Here $X_{H_t}\in\Vect(M)$ is given by $\iota(X_{H_t})\om=-dH_t$.
Now let $ -\infty \le a < b \le\infty $ and denote by
$\Pp^{[a,b)}(H;\alpha)$ the set of $1$-periodic solutions of the
Hamiltonian system associated to $H$ that represent the class $\alpha$
and whose action lies in the interval $[a,b)$:
$$
\Pp^{[a,b)}(H;\alpha) :=
\Pp^b(H;\alpha)\setminus\Pp^a(H;\alpha),\qquad \Pp^a(H;\alpha) :=
\left\{x\in\Pp(H;\alpha)\,|\,\Aa_H(x) < a\right\}.
$$
Suppose that $H\in\Hh$ is a Hamiltonian function that satisfies the
following hypothesis:
\begin{enumerate}
  \item[(H)] {\it $a,b\notin\Ss(H;\alpha)$ and every $1$-periodic
     orbit $x\in\Pp(H;\alpha)$ is nondegenerate.}
\end{enumerate}
Then the Floer homology group $\HF^{[a,b)}(H;\alpha)$ is defined as
the homology of a chain complex over $\Z_2$ generated by the
$1$-periodic orbits in $\Pp^{[a,b)}(H;\alpha).$\footnote{We use the
  convention that the complex generated by the empty set is $0$.}  It
is useful to think of this chain complex as the quotient
$$
\CF^{[a, b)}(H;\alpha) := \CF^b(H;\alpha)/\CF^a(H;\alpha),\qquad
\CF^a(H;\alpha) := \bigoplus_{x\in\Pp^a(H;\alpha)}\Z_2 x.
$$
The Floer boundary operator is defined as follows.  Let
$J_t=J_{t+1}\in\Jj(M,\om)$ be a $t$-dependent smooth family of
$\om$-compatible almost complex structures on $\overline M$ such that
$J_t$ is convex and independent of $t$ near the boundary $\p\overline
M$ (see Section~\ref{sb:convex}).  Consider the Floer differential
equation
\begin{equation}\label{eq:floer}
   \p_su + J_t(u)(\p_tu - X_{H_t}(u)) = 0.
\end{equation}
For a smooth solution $u:\R\times S^1\to M$ of~(\ref{eq:floer}) define
its energy to be
$$
E(u):=\int_0^1\int_{-\infty}^\infty|\p_su|^2\,ds dt.
$$
Then if $u:\R\times S^1\to M$ is a smooth solution
of~(\ref{eq:floer}) with finite energy then the limits
\begin{equation}\label{eq:limits}
   \lim_{s\to\pm\infty}u(s,t) = x^\pm(t),\qquad
   \lim_{s\to\pm\infty}\p_su(s,t) = 0
\end{equation}
exist and are uniform in the $t$-variable. Moreover, $x^\pm\in\Pp(H)$
and we have
$$E(u) = \Aa_H(x^-) - \Aa_H(x^+).$$
The following observations allow
us to define Floer homology groups in the present situation.
\begin{enumerate}
  \item[(i)] Since every periodic solution $x\in\Pp(M;\alpha)$ is
   nonconstant (the class $\alpha$ is nontrivial) and $J$ is convex
   near the boundary, there exists an open set $U\subset M$ such that
   $M\setminus U$ is compact and $u(\R\times S^1)\cap U=\emptyset$ for
   every finite energy solution of~(\ref{eq:floer}).
  \item[(ii)] By~(i) and the energy identity, the space of finite
   energy solutions of~(\ref{eq:floer}) is compact with respect to
   $\Cinf$-convergence on compact sets, i.e. only the splitting into a
   finite sequence of adjacent Floer connecting orbits can occur in
   the limit.
  \item[(iii)] For a generic family of almost complex structures
   $J=\{J_t\}$ (that are convex and independent of $t$ on~$U$) the
   linearized operator for equation~(\ref{eq:floer}) is surjective for
   every finite energy solution of~(\ref{eq:floer}) in the homotopy
   class $\alpha$ (see~\cite{F-H-S}).  Such a family of almost complex
   structures is called {\bf regular} and the space of regular
   families of almost complex structures will be denoted by
   $\Jreg(H;\alpha)$.
\end{enumerate}
For every $J\in\Jreg(H;\alpha)$ and every pair $x^\pm\in\Pp(H;\alpha)$
the space $\Mm(x^-,x^+;H,J)$ of solutions of~(\ref{eq:floer})
and~(\ref{eq:limits}) is a smooth manifold whose dimension near a
solution $u$ of~(\ref{eq:floer}) and~(\ref{eq:limits}) is given by the
difference of the Conley--Zehnder indices (see~\cite{S-Z}) of $x^-$
and $x^+$ (relative to $u$).  The subspace of solutions of index one
will be denoted by $\Mm^1(x^-,x^+;H,J)$.  For $J\in\Jreg(H;\alpha)$ it
follows from~(i) and~(ii) that the quotient $\Mm^1(x^-,x^+;H,J)/\R$
(modulo time shift) is a finite set for every pair
$x^\pm\in\Pp(H;\alpha)$.  The Floer boundary operator $\p^{H,J}$ on
$\CF^b(H;\alpha)$ is defined by
$$
\p^{H,J}x := \sum_{y\in\Pp^b(H;\alpha)} \#(\Mm^1(x,y;H,J)/\R)\,y
$$
for every $x\in\Pp^b(H;\alpha)$.  That this is indeed a boundary
operator is proved as in Floer's original work~\cite{F}.  The energy
identity shows that $\CF^a(H;\alpha)$ is a subcomplex, namely it is
invariant under the Floer boundary operator. We thus get an induced
boundary operator $[\p^{H,J}]$ on the quotient $\CF^{[a,
  b)}(H;\alpha)$. We denote the homology of the quotient complex by
$$
\HF^{[a,b)}(H,J;\alpha) := \frac{\ker ( [\p^{H,J}]:\CF^{[a,
    b)}(H;\alpha) \to \CF^{[a, b)}(H;\alpha))} {{\rm
    im}([\p^{H,J}]:\CF^{[a, b)}(H;\alpha)\to\CF^{[a, b)}(H;\alpha))}.
$$
These Floer homology groups are independent of the choice of the
almost complex structure $J=\{J_t\}_{t\in S^1}$ in the sense that for
any two almost complex structures $J_0,J_1\in\Jreg(H;\alpha)$ there is
a natural isomorphism
$$
\tau_{J_1J_0}:\HF^{[a,b)}(H,J_0;\alpha)\to
\HF^{[a,b)}(H,J_1;\alpha).
$$
If the two almost complex structure agree near the boundary then
this follows from the standard arguments as in Floer's original
paper~\cite{F} (choose a homotopy of almost complex structures
$\{J_{s,t}\}$ from $J_0$ to $J_1$, independent of $s$ and $t$ near the
boundary, and use the solutions of equation~(\ref{eq:homotopy}) below
with $H_{s,t}=H_t$ to construct the isomorphism between the two Floer
homology groups; see also~\cite{Sa,S-Z}).  To show that the Floer
homology groups are also independent of the choice of the {\it convex}
almost complex structure near the boundary one can use the fact that
the space of convex almost complex structures near the boundary is
connected (Remark~\ref{rmk:convex}) and that the Floer chain complex
associated to a regular almost complex structure remains unchanged
under sufficiently small perturbations of~$J$.  The upshot is that the
Floer homology groups are independent of $J$ up to natural
isomorphisms.  For this reason we shall sometimes drop the argument
$J$ and refer to $\HF^{[a,b)}(H;\alpha):=\HF^{[a,b)}(H,J;\alpha)$ as
the {\em Floer homology associated to~$H$}.


\subsection{Homotopy invariance}\label{sb:homotopy} 

Following the work of Floer--Hofer~\cite{F-H-2},
Cieliebak-Floer-Hofer~\cite{C-F-H-1}, and Viterbo~\cite{Vi-1} we
describe the local isomorphisms of Floer homology in a given interval
of the action spectrum.  Consider the space
$$
\Hh^{a,b}(M;\alpha) :=
\left\{H\in\Hh(M)\,|\,a,b\notin\Ss(H;\alpha)\right\}
$$
of all Hamiltonians $H\in\Hh$ that do not contain $a$ and $b$ in
their action spectrum.  We consider the space $\Hh$ with the strong
Whitney $\Cinf$-topology. Note that the action spectrum
$\Ss(H;\alpha)$ is compact for every $H$ and is a lower semicontinuous
function of $H$ (i.e. for every open neighbourhood $V\subset\R$ of
$\Ss(H;\alpha)$ there exists a neighbourhood $\Uu\subset\Hh$ of $H$
such that $\Ss(H';\alpha)\subset V$ for every $H'\in\Uu$).  Hence the
set $\Hh^{a,b}(M;\alpha)$ is open in $\Hh$.  We now explain why the
Floer homology groups $\HF^{[a,b)}(H;\alpha)$ are independent of $H$
in every component of $\Hh^{a,b}(M;\alpha)$.

Fix a Hamiltonian function $H\in\Hh^{a,b}(M;\alpha)$ and choose a
(convex) neighbourhood $\Uu$ of $H$ such that
$\Uu\subset\Hh^{a,b}(M;\alpha)$.  Now suppose that $H^+,H^-\in\Uu$
satisfy~$(H)$, i.e.  all periodic solutions $x\in\Pp(H^\pm;\alpha)$
are nondegenerate.  Connect $H^-$ and $H^+$ by a smooth homotopy
$\R\mapsto\Uu:s\mapsto H_s=\{H_{s,t}\}$ such that $H_{s,t}=H^-_t$ for
$s\le -T$ and $H_{s,t}=H^+_t$ for $s\ge T$. Consider the equation
\begin{equation}\label{eq:homotopy}
   \p_su + J_{s,t}(u)(\p_tu - X_{H_{s,t}}(u)) = 0,
\end{equation}
where $s\mapsto\{J_{s,t}\}$ is a {\bf regular homotopy} of families of
almost complex structures.  This means that $J_{s,t}$ satisfies the
following conditions.
\begin{enumerate}
  \item[$\bullet$] $J_{s,t}$ is convex and independent of $s$ and $t$
   near the boundary of $\overline M$.
  \item[$\bullet$] $J_{s,t}=J^-_t$ is regular for $H^-_t$ for $s\le
   -T$.
  \item[$\bullet$] $J_{s,t}=J^+_t$ is regular for $H^+_t$ for $s\ge
   T$.
  \item[$\bullet$] The finite energy solutions of~(\ref{eq:homotopy})
   are transverse (i.e. the associated Fredholm operators are
   surjective) and hence form finite dimensional moduli spaces.
\end{enumerate}
The key observation is the energy identity
\begin{equation}\label{eq:energy}
   E(u) = \Aa_{H^-}(x^-) - \Aa_{H^+}(x^+) 
   + \int_0^1\int_{-\infty}^\infty (\p_sH)(s,t,u(s,t))\,ds dt
\end{equation}
for every solution of~(\ref{eq:homotopy}) and~(\ref{eq:limits}).  It
follows from~(\ref{eq:energy}) that
$$
\Aa_{H^+}(x^+) \le \Aa_{H^-}(x^-) +
\int_{-\infty}^\infty\max_{S^1\times M}\p_sH_s\,ds.
$$
In particular, if the homotopy has the form
$H_{s,t}:=H_{0,t}+\beta(s)(H^+_t-H^-_t)$ for a nondecreasing function
$\beta:\R\to[0,1]$ we obtain $\p_sH_s=\dot\beta(s)(H^+-H^-)$ and hence
\begin{equation}\label{eq:key}
   \Aa_{H^+}(x^+) 
   \le \Aa_{H^-}(x^-)
   + \max_{S^1\times M}(H^+-H^-).
\end{equation}
Now choose $\eps>0$ such that
$$
\Ss(H^\pm;\alpha)\cap[a-4\eps,a+4\eps]=\emptyset,\qquad
\Ss(H^\pm;\alpha)\cap[b-4\eps,b+4\eps]=\emptyset,
$$
and suppose that $\sup_{S^1\times M}|H^\pm-H|\le\eps$.  Then
$\sup_{S^1\times M}|H^+-H^-|\le2\eps$ and hence, by~(\ref{eq:key}),
the Floer chain map (see~\cite{F,F-H-2,C-F-H-1,Vi-1,Sa,S-Z}) from
$\CF(H^-;\alpha)$ to $\CF(H^+;\alpha)$ defined by the solutions
of~(\ref{eq:homotopy}) preserves the subcomplexes $\CF^a$ and $\CF^b$.
The same applies to the Floer chain map from $\CF(H^+;\alpha)$ to
$\CF(H^-;\alpha)$ and to the chain homotopy equivalence associated to
a suitable homotopy of homotopies.  Hence the solutions
of~(\ref{eq:homotopy}) define a homomorphism
$\CF^{[a,b)}(H^-;\alpha)\to \CF^{[a,b)}(H^+;\alpha)$ which induces an
isomorphism of Floer homology, whenever $H^\pm$ are sufficiently close
to a given Hamiltonian function $H\in\Hh^{a,b}(M;\alpha)$.
  
\begin{rem}{\bf (Local isomorphisms)}\label{rmk:degenerate}
   The above discussion shows that every Hamiltonian function
   $H\in\Hh^{a,b}(M;\alpha)$ has a neighbourhood $\Uu$ such that the
   Floer homology groups $\HF^{[a,b)}(H',J';\alpha)$, for every
   $H'\in\Uu$ that satisfies~$(H)$ and every regular almost complex
   structure~$J'\in\Jreg(H';\alpha)$, are naturally isomorphic.  We
   can use these local isomorphisms to define the Floer homology
   groups $\HF^{[a,b)}(H;\alpha)$ for every Hamiltonian
   $H\in\Hh^{a,b}(M;\alpha)$, whether or not it satisfies~$(H)$.
\end{rem}

\begin{rem}{\bf (Contractible loops)}\label{rmk:contractible}
   When $\alpha\in\widetilde{\pi}_1(M)$ is the homotopy class of the
   constant loops we are not allowed to work with intervals $[a,b)$
   that contain zero, since the Hamiltonians we work with always have
   degenerate periodic orbits with action zero as they vanish at
   infinity.  In this case we are forced to work with either $0< a < b
   \leq \infty$ or $-\infty \leq a < b < 0$.
\end{rem}

\begin{rem}{\bf (Composition)}\label{rmk:local}
   We emphasize that the canonical isomorphism
   $$
   \HF^{[a,b)}(H^-,J^-;\alpha)\to\HF^{[a,b)}(H^+,J^+;\alpha)
   $$
   only exists locally, when $H^\pm$ are sufficiently close to a
   given Hamiltonian function $H\in\Hh^{a,b}(M;\alpha)$.  It is easy
   to construct Hamiltonian functions $H_0,H_1\in\Hh^{a,b}(M;\alpha)$
   such that $\HF^{[a,b)}(H_0;\alpha)$ is not isomorphic to
   $\HF^{[a,b)}(H_1;\alpha)$.  If $H_0$ and $H_1$ belong to the same
   component of $\Hh^{a,b}(M;\alpha)$ then there is a smooth path $
   [0,1]\to\Hh^{a,b}(M;\alpha):s\mapsto H_s $ connecting $H_0$ to
   $H_1$.  Hence in this case $\HF^{[a,b)}(H_0;\alpha)$ is isomorphic
   to $\HF^{[a,b)}(H_1;\alpha)$.  However, in general the isomorphism
   cannot be defined directly in terms of the solutions
   of~(\ref{eq:homotopy}).  It can only be constructed as a
   composition of isomorphisms
   $$
   \HF^{[a,b)}(H_{s_i};\alpha)\to\HF^{[a,b)}(H_{s_{i+1}};\alpha)
   $$
   for a regular homotopy, where each of these isomorphisms is
   defined in terms of the solutions of~(\ref{eq:homotopy}).
   Moreover, it is an open question if this composition along a loop
   $s\mapsto H_s$ with $H_0=H_1$ is always the identity.
\end{rem}


\subsection{Monotone homotopies}\label{sb:monotone}

Suppose that $H_0,H_1\in\Hh^{a,b}(M;\alpha)$ satisfy
$$
H_0(t,x) \ge H_1(t,x)
$$
for all $(t,x)\in\R\times M$ as well as~$(H)$.  Then there exists a
homotopy $s\mapsto H_s$ from $H_0$ to $H_1$ such that $ \p_sH_s\le 0.
$ We call such a homotopy of Hamiltonian functions {\bf monotone}.  In
the monotone case it follows from~(\ref{eq:energy}) that the Floer
chain map $\CF(H_0;\alpha)\to\CF(H_1;\alpha)$, defined in terms of the
solutions of~(\ref{eq:homotopy}) preserves the subcomplexes $\CF^a$
and $\CF^b$.  Hence every monotone homotopy $s\mapsto H_s$ induces a
natural homomorphism
$$
\sigma_{H_1H_0}:\HF^{[a,b)}(H_0;\alpha)\to\HF^{[a,b)}(H_1;\alpha).
$$
We call such a homomorphism {\bf monotone}.  The standard arguments
in Floer homology~\cite{F,F-H-2,C-F-H-1,Vi-1,Sa,S-Z} show that this
homomorphism is independent of the choice of the monotone homotopy of
Hamiltonians, used to define it, and that
$$
\sigma_{H_2H_1}\circ\sigma_{H_1H_0} = \sigma_{H_2H_0},
$$
whenever $H_0,H_1,H_2\in\Hh^{a,b}(M;\alpha)$ satisfy $H_0\ge H_1\ge
H_2$, and $\sigma_{HH}=\id$ for every $H \in \Hh^{a,b}(M;\alpha)$.

The homomorphism $\sigma_{H_1H_0}$ is in general neither injective nor
surjective. For example, it may happen that during the homotopy the
action of some periodic orbit of $H_s$ leaves or enters the interval
$[a,b)$. It turns out that this is the only possible reason for
$\sigma_{H_1H_0}$ not to be an isomorphism.  More precisely, we have
the following proposition which is an easy consequence of the theory
developed in~\cite{F-H,C-F-H-1} (as outlined above) and appears in an
explicit form in~\cite{Vi-1}.

\begin{prop} \label{P:homotopy}
   Let $-\infty\le a < b\le\infty$, $\alpha\in\widetilde{\pi}_1(M)$ be
   a nontrivial homotopy class, and $H_0,H_1\in\Hh^{a,b}(M;\alpha)$ be
   such that $H_0\ge H_1$.  Suppose that there exists a monotone
   homotopy $\{H_s\}_{0\le s\le 1}$ from $H_0$ to $H_1$ such that
   $H_s\in\Hh^{a,b}(M;\alpha)$ for every $s\in[0,1]$.  Then
   $\sigma_{H_1H_0}:\HF^{[a,b)}(H_0;\alpha)\to\HF^{[a,b)}(H_1;\alpha)$
   is an isomorphism. This continues to hold for the trivial homotopy
   class $\alpha=0$ provided that $0\notin[a,b]$.
\end{prop}

\begin{proof}
   The monotone homomorphism $ \sigma_{H_{s_1}H_{s_0}} $ agrees with
   the local isomorphism of Section~\ref{sb:homotopy} whenever $s_0$
   and $s_1$ are both sufficiently close to a number $s\in[0,1]$ such
   that $\Hh_s\in\Hh^{a,b}(H_s;\alpha)$.  By assumption, we have
   $H_s\in\Hh^{a,b}(H_s;\alpha)$ for every $s\in[0,1]$.  Hence we can
   write $\sigma_{H_1H_0}$ as a composition of finitely many
   isomorphisms of the form $\sigma_{H_{s_{i+1}}H_{s_i}}$, where
   $0=s_0<s_1<\cdots<s_{N-1}<s_N=1$.
\end{proof}


\subsection{Direct and inverse limits} 
\label{sb:Limits} 

The next step towards defining the relative capacity is to define two
kinds of {\em symplectic homologies} especially suited for our
purposes. The definitions of these invariants require the algebraic
notions of direct and inverse limits.  In this subsection we recall
the basic definitions (for more details see~\cite{G-M}, but note that
below we use somewhat different conventions than theirs).

Let $(I,\preceq)$ be a partially ordered set.  Think of $I$ as a
category with precisely one morphism from $i$ to $j$ whenever
$i\preceq j$. Let $R$ be a commutative ring. A {\bf partially ordered
  system} of $R$-modules over $I$ is a functor from $(I,\preceq)$ into
the category of $R$-modules. We write this functor as a pair
$(\G,\sigma)$ where $\G$ assigns to each $i\in I$ an $R$-module $\G_i$
and $\sigma$ assigns to each pair $i,j\in I$ with $i\preceq j$ a
homomorphism $\sigma_{ji}:\G_i\to\G_j$ such that
$\sigma_{kj}\circ\sigma_{ji}=\sigma_{ki}$ for $i\preceq j\preceq k$
and $\sigma_{ii}=\id$ is the identity map on $\G_i$.

The partially ordered set $(I,\preceq)$ is called {\bf upward
  directed} if for every pair $i,j\in I$ there exists an $\ell\in I$
such that $i\preceq\ell$ and $j\preceq\ell$.  In this case the functor
$(\G,\sigma)$ is called a {\bf directed system of $R$-modules}. The
{\bf direct limit} of such a directed system is defined as the
quotient
$$
\varinjlim\G := \underset{i\in I}{\varinjlim}\G_i :=
\left\{(i,x)\,|\,i\in I,\,x\in\G_i\right\}/\sim
$$
where $(i,x)\sim(j,y)$ iff there exists an $\ell\in I$ such that
$i\preceq\ell$, $j\preceq\ell$ and $\sigma_{\ell i}(x)=\sigma_{\ell
  j}(y)$.  Since $I$ is upward directed, this is an equivalence
relation. The direct limit is an $R$-module with the operations
$[i,x]+[j,y]:=[\ell,\sigma_{\ell i}(x)+\sigma_{\ell j}(y)]$ for
$\ell\in I$ such that $i\preceq\ell$ and $j\preceq\ell$ and
$r[i,x]:=[i,rx]$ for every $r \in R$. For $i\in I$ we denote by
$\iota_i:\G_i\to\varinjlim\G$ the homomorphism given by
$\iota_i(x):=[i,x]$.  Then $ \iota_i=\iota_j\circ\sigma_{ji} $ for
$i\preceq j$. Despite the notation, $\iota_i$ need not be injective.
Up to isomorphism the direct limit is characterized by the following
universal property.  If $H$ is any $R$-module and $\tau_i:\G_i\to H$
is a family of homomorphisms, indexed by $i\in I$, such that
$\tau_i=\tau_j\circ\sigma_{ji}$ whenever $i\preceq j$, then there
exists a unique homomorphism $\tau:\varinjlim\G\to H$ such that
$\tau_i=\tau\circ\iota_i$ for every $i\in I$.  (The homomorphism
$\tau$ is given by $[i,x]\mapsto\tau_i(x)$.)

The partially ordered set $(I,\preceq)$ is called {\bf downward
  directed} if for every pair $i,j\in I$ there exists a $k\in I$ such
that $k\preceq i$ and $k\preceq j$.  In this case the functor
$(\G,\sigma)$ is called an {\bf inverse system of $R$-modules}. The
{\bf inverse limit} of such an inverse system is defined as
$$
\varprojlim\G := \underset{i \in I}{\varprojlim}\G_i :=
\left\{\{x_i\}_{i\in I} \in \prod_{i \in I} G_i \biggm| i\preceq
   j\IMP\sigma_{ji}(x_i) = x_j\right\}.
$$
For $j\in I$ we denote by $\pi_j:\varprojlim\G\to\G_j$ the obvious
projection to the $j$th component.  Then $ \pi_j=\sigma_{ji}\circ\pi_i
$ for $i\preceq j$.  Despite the notation, $\pi_j$ need not be
surjective.  Up to isomorphism the inverse limit is characterized by
the following universal property.  If $H$ is any $R$-module and
$\tau_j:H\to\G_j$ is a family of homomorphisms, indexed by $j\in I$,
such that $\tau_j=\sigma_{ji}\circ\tau_i$ whenever $i\preceq j$, then
there exists a unique homomorphism $\tau:H\to\varprojlim\G$ such that
$\tau_j=\pi_j\circ\tau$ for every $j\in I$.  (The homomorphism $\tau$
is given by $y\mapsto\{\tau_i(y)\}_{i\in I}$.)

\medskip
\begin{remnonum}
   Note that inverse and direct limits are related via the following
   duality. Let $(\G,\sigma)$ be a directed system of $R$-modules and
   $H$ be any $R$-module. Denote by $(I^*,\preceq^*)$ the oppositely
   partially ordered set, namely $I^*:=I$ and $i\preceq^*j$ iff
   $i\succeq j$.  Then there exists a canonical isomorphism
   $$
   \Hom_R \Big(\underset{i\in I}{\varinjlim}\G_i,H\Big) \cong
   \underset{i\in I^*}{\varprojlim}\Hom_R(\G_i,H).
   $$
   In particular, if $R$ is a field and the $G_i$ are vector spaces
   over $R$ then $ (\underset{i\in I}{\varinjlim}\G_i)^* \cong
   \underset{i \in I^*}{\varprojlim} \G_i^*.  $
\end{remnonum}

\medskip
\noindent
In most of our applications the partially ordered set $(I,\preceq)$
will be {\bf bidirected}, i.e. both upward and downward directed.  In
this case we call the functor $(\G,\sigma)$ a {\bf bidirected system
  of $R$-modules}.  The next lemma follows directly from the
definitions.

\begin{lem}\label{le:limits}
   Let $(I,\preceq)$ be a downward directed partially ordered set and
   $I'\subset I$ be an upward directed subset (with respect to the
   restriction of the partial order $\preceq$ to $I'$).  Let
   $(G,\sigma)$ be a partially ordered system of $R$-modules over $I$.
   Then there exists a unique homomorphism
   $$
   T:\underset{i\in I}{\varprojlim} G_i \to \underset{i'\in
     I'}{\varinjlim} G_{i'}
   $$
   such that the following diagram commutes for all $j',k'\in I'$
   with $j'\preceq k'$:
   \[
   \xymatrix{ {\underset{i \in I}{\varprojlim}\G_i} \ar[rr]^-T
     \ar[d]_-{\pi_{_{j'}}} && {\underset{i' \in I'}{\varinjlim}\G_{i'}} \\
     \G_{j'} \ar[rr]^-{\sigma_{_{k'j'}}} && {\G_{k'}}
     \ar[u]_-{\iota_{_{k'}}} }
   \] 
\end{lem}

\begin{proof}
   The map $T$ is given by $\{x_i\}_{i\in I}\mapsto[i',x_{i'}]$ for
   $i'\in I'$.  By definition of direct and inverse limits, this map
   is independent of the choice of~$i'$.
\end{proof}

In Lemma~\ref{le:limits} we {\bf do not} assume that for every $i\in
I$ there exists an $i'\in I'$ such that $i\preceq i'$ (and indeed this
condition is not satisfied in our application).  If this holds then
the map $T$ factors through $\pi_i$ for every $i\in I$ and not just
for $i\in I'$.  However, there are examples where $\G_i=\{0\}$ for
some $i\in I$ and $T\ne0$.


\subsection{Exhausting sequences}\label{sb:exhaust}

To compute direct and inverse limits we introduce the notion of
exhausting sequences.  Let $(\G,\sigma)$ be a partially ordered system
of $R$-modules over $(I,\preceq)$ and denote $
\Z^\pm:=\left\{\nu\in\Z\,|\,\pm\nu>0\right\}.  $ A sequence
$\{i_\nu\}_{\nu\in\Z^+}$ is called {\bf upward exhausting} for
$(\G,\sigma)$ iff the following holds
\begin{enumerate}
  \item[$\bullet$] For every $\nu\in\Z^+$ we have $i_\nu\preceq
   i_{\nu+1}$ and
   $\sigma_{i_{\nu+1}i_\nu}:\G_{i_\nu}\to\G_{i_{\nu+1}}$ is an
   isomorphism.
  \item[$\bullet$] For every $i\in I$ there exists a $\nu\in\Z^+$ such
   that $i\preceq i_\nu$.
\end{enumerate}
A sequence $\{i_\nu\}_{\nu\in\Z^-}$ is called {\bf downward
  exhausting} for $(\G,\sigma)$ iff the following holds
\begin{enumerate}
  \item[$\bullet$] For every $\nu\in\Z^-$ we have $i_{\nu-1}\preceq
   i_\nu$ and $\sigma_{i_\nu i_{\nu-1}}:\G_{i_{\nu-1}}\to\G_{i_\nu}$
   is an isomorphism.
  \item[$\bullet$] For every $i\in I$ there exists a $\nu\in\Z^-$ such
   that $i_\nu\preceq i$.
\end{enumerate}
The importance of such sequences is that they can be used to simplify
computations of direct and inverse limits.

\begin{lem} \label{L:exhaust}
   Let $(\G,\sigma)$ be a partially ordered system of $R$-modules over
   $(I,\preceq)$.
   
   \smallskip
   \noindent{\bf (i)} 
   If $\{i_\nu\}_{\nu\in\Z^+}$ is an upward exhausting sequence for
   $(\G,\sigma)$ then $(I,\preceq)$ is upward directed and the
   homomorphism $ \iota_{i_\nu}:\G_{i_\nu}\to\varinjlim\G $ is an
   isomorphism for every $\nu\in\Z^+$.

   \smallskip
   \noindent{\bf (ii)} 
   If $\{i_\nu\}_{\nu\in\Z^-}$ is a downward exhausting sequence for
   $(\G,\sigma)$ then $(I,\preceq)$ is downward directed and the
   homomorphism $ \pi_{i_\nu}:\varprojlim\G\to\G_{i_\nu} $ is an
   isomorphism for every $\nu\in\Z^-$.
\end{lem}

\begin{proof}
   To prove~(i) we fix an integer $\nu\in\Z^+$.  Let $x\in\G_{i_\nu}$
   and suppose that $\iota_{i_\nu}(x)=0$.  Then there exists an $i\in
   I$ such that $i_\nu\preceq i$ and $\sigma_{ii_\nu}(x)=0$.  Choose
   an integer $\nu'\ge\nu$ such that $i\preceq i_{\nu'}$.  Then $
   \sigma_{i_{\nu'}i_\nu}(x) =
   \sigma_{i_{\nu'}i}\circ\sigma_{ii_\nu}(x) = 0 $ and hence $x=0$.
   Hence $\iota_{i_\nu}$ is injective.  Now let $y\in\G_j$ and choose
   an integer $\nu'\ge\nu$ such that $j\preceq i_{\nu'}$. Since
   $\sigma_{i_{\nu'}i_\nu}$ is surjective there exists an
   $x\in\G_{i_\nu}$ such that
   $\sigma_{i_{\nu'}i_\nu}(x)=\sigma_{i_{\nu'}j}(y)$.  Hence
   $(j,y)\sim(i_\nu,x)$.  This shows that $\iota_{i_\nu}$ is
   surjective
   
   To prove~(ii) we fix an integer $\nu\in\Z^-$.  Let $\{x_i\}_{i\in
     I}\in\varprojlim\G$ such that $x_{i_\nu}=0$.  Given $i\in I$
   choose an integer $\nu'\le\nu$ such that $\iota_{\nu'}\preceq i$.
   Then $\sigma_{i_\nu i_{\nu'}}(x_{i_{\nu'}})=x_{i_\nu}=0$, hence
   $x_{i_{\nu'}}=0$, and hence
   $x_i=\sigma_{ii_{\nu'}}(x_{i_{\nu'}})=0$.  This shows that
   $\pi_{i_\nu}$ is injective.  Now let $x\in\G_{i_\nu}$. Given $i\in
   I$, choose an integer $\nu'\le\nu$ such that $i_{\nu'}\le i$ and
   define $x_i\in\G_i$ by
   $$
   x_i:=\sigma_{ii_{\nu'}}(x_{i_\nu'}),\qquad \sigma_{i_\nu
     i_{\nu'}}(x_{i_{\nu'}}) := x.
   $$
   Since $\sigma_{i_\nu i_{\nu'}}$ is surjective the element
   $x_i\in\G_i$ exists.  Since $\sigma_{i_\nu i_{\nu'}}$ is injective,
   the element $x_i$ is unique, and it is independent of the choice of
   $\nu'$.  We prove that $\sigma_{ji}(x_i)=x_j$ whenever $i\preceq
   j$.  To see this choose $\nu'\le\nu$ such that $i_{\nu'}\preceq
   i\preceq j$.  Then
   $$
   x_j = \sigma_{ji_{\nu'}}(x_{i_{\nu'}}) = \sigma_{ji}\circ
   \sigma_{ii_{\nu'}}(x_{i_{\nu'}}) = \sigma_{ji}(x_i).
   $$
   Hence $\{x_i\}_{i\in I}\in\varprojlim\G$ and
   $\pi_{i_\nu}(\{x_i\}_{i\in I})=x_{i_\nu}=x$.  This shows that
   $\pi_{i_\nu}$ is surjective.
\end{proof}


\subsection{Symplectic homology}\label{sb:SH}

The set $\Hh(M)$ of Hamiltonian functions on $S^1\times M$ with
compact support is partially ordered by
$$
H_0 \preceq H_1\qquad\IFF\qquad H_0(t,x) \ge
H_1(t,x)\;\;\;\forall\;(t,x)\in S^1\times M.
$$
This defines a bidirected partial order on $\Hh(M)$.  Let
$\alpha\in\widetilde{\pi}_1(M)$ be a nontrivial homotopy class and
$a,b\in\R\cup\{\pm\infty\}$ such that $a<b$.  As in
Section~\ref{sb:homotopy} we denote by $\Hh^{a,b}(M;\alpha)$ the
subset of all Hamiltonian functions $H\in\Hh(M)$ such that
$a,b\notin\Ss(H;\alpha)$.  In Subsection~\ref{sb:monotone} we have
seen that there is a natural homomorphism
$$
\sigma_{H_1H_0}:\HF^{[a,b)}(H_0;\alpha) \to \HF^{[a,b)}(H_1;\alpha)
$$
whenever $H_0,H_1\in \Hh^{a,b}(M;\alpha)$ satisfy $H_0\preceq H_1$.
These homomorphisms define an inverse (in fact bidirected) system of
Floer homology groups over $(\Hh^{a,b}(M;\alpha),\preceq)$.  The
inverse limit of this system is called the {\bf symplectic homology}
of $M$ in the homotopy class $\alpha$ for the action interval $[a,b)$.
A version of this homology group was introduced
in~\cite{F-H-2,C-F-H-1} for the homotopy class of contractible loops
and later on for general homotopy classes in~\cite{Cieliebak-2}. We
denote it by
$$
\SHi^{[a,b)}(M;\alpha) := \underset{H \in \Hh^{a,b}(M;\alpha)}
\varprojlim\HF^{[a,b)}(H;\alpha).
$$
Now fix a compact subset $A\subset M$ and a constant $c\in\R$.
Consider the set $\Hh^{a,b}_c(M,A;\alpha)$ of all Hamiltonian
functions $H\in\Hh^{a,b}(M;\alpha)$ that satisfy $H>c$ on $S^1\times
A$, namely
$$
\Hh^{a,b}_c(M,A;\alpha) := \left\{H\in\Hh^{a,b}(M;\alpha)\,\Big|\,
   \inf_{S^1\times A}H>c\right\}.
$$
This gives rise to a directed (in fact bidirected) system of Floer
homology groups over $ (\Hh^{a,b}_c(M,A;\alpha),\preceq).  $ The
direct limit of this system is called the {\bf relative symplectic
  homology} of the pair $(M,A)$ at the {\bf level $c$} in the homotopy
class $\alpha$ for the action interval $[a,b)$.  We denote it by
$$
\SHd^{[a,b);c}(M,A;\alpha) := \underset{H \in
  \Hh^{a,b}_c(M,A;\alpha)} \varinjlim\HF^{[a,b)}(H;\alpha).
$$
\begin{remnonum}
   Since we have chosen to work with $\mathbb{Z}_2$-coefficients all
   the Floer homology groups $\HF^{[a,b)}(H;\alpha)$ are in fact
   $\mathbb{Z}_2$-vector spaces. Consequently also the symplectic
   homologies $\SHi^{[a,b)}(M;\alpha)$ and
   $\SHd^{[a,b);c}(M,A;\alpha)$ have the structure of
   $\mathbb{Z}_2$-vector spaces.
\end{remnonum}

\medskip
\noindent
An important feature of absolute and relative symplectic homologies is
the existence of a homomorphism between them which factors through
Floer homology.
\begin{prop} \label{P:factor}
   Let $\alpha\in\widetilde{\pi}_1(M)$ be a nontrivial homotopy class
   and suppose that $-\infty\le a<b\le\infty$.  Then, for every
   $c\in\R$, there exists a unique homomorphism
   $$
   T_\alpha^{[a,b);c}:\SHi^{[a,b)}(M;\alpha)\to\SHd^{[a,b);c}(M,A;\alpha)
   $$
   such that for any two Hamiltonian functions
   $H_0,H_1\in\Hh^{a,b}_c(M,A;\alpha)$ with $H_0\ge H_1$ the following
   diagram commutes:
   \[
   \xymatrix{ {\SHi^{[a,b)}(M;\alpha)} \ar[rr]^-{T_{\alpha}^{[a,b);c}}
     \ar[d]_-{\pi_{_{H_0}}}
     && {\SHd^{[a,b);c}(M,A;\alpha)} \\
     \HF^{[a,b)}(H_0;\alpha) \ar[rr]^-{\sigma_{_{H_1,H_0}}} &&
     {\HF^{[a,b)}(H_1;\alpha)} \ar[u]_-{\iota_{_{H_1}}} }
   \]
   Here
   $\pi_{_{H_0}}:\SHi^{[a,b)}(M;\alpha)\to\HF^{[a,b)}(H_0;\alpha)$ and
   $\iota_{_{H_1}}:\HF^{[a,b)}(H_1;\alpha)\to\SHd^{[a,b);c}(M,A;\alpha)$
   are the canonical homomorphisms introduced in
   Section~\ref{sb:Limits}.  In particular, since $\sigma_{HH}=\id$
   for every $H\in\Hh^{a,b}_c(M,A;\alpha)$, we have
   \[
   \xymatrix{ {\SHi^{[a,b)}(M;\alpha)} \ar[rr]^-{T_{\alpha}^{[a,b);c}}
     \ar[dr]_-{\pi_{_H}}
     && {\SHd^{[a,b);c}(M,A;\alpha)} \\
     & \HF^{[a,b)}(H;\alpha) \ar[ur]_-{\iota_{_H}} }
   \]   
   The statements above continue to hold also for the trivial class
   $\alpha=0$, provided that $0 \notin [a,b]$.
\end{prop}

\begin{proof} The proof follows at once from Lemma~\ref{le:limits}.
\end{proof}


\subsection{The homological relative capacity}\label{sb:Rel-cap}

For every nontrivial homotopy class $\alpha\in\widetilde{\pi}_1(M)$
and every real number $c>0$ we define the set
$$
\A_c(M,A;\alpha) := \left\{a\in\R \,\Big| \textnormal{ The
     homomorphism } T_\alpha^{[a,\infty);c} \textnormal{ does not
     vanish} \right\},
$$
where $T_\alpha^{[a,\infty);c}: \SHi^{[a,\infty)}(M;\alpha) \to
\SHd^{[a,\infty);c}(M,A;\alpha)$ is the homomorphism from
Proposition~\ref{P:factor}.

For the trivial homotopy class $\alpha=0\in\widetilde{\pi}_1(M)$ we
define $\A_c(M,A;0)$ by the same formula except that we only consider
real numbers $a>0$ (for which $T_0^{[a,\infty);c}\ne0$).  The {\bf
  homological relative capacity} of the pair $(M,A)$ is the function
$$
\widehat{C}(M,A):\widetilde{\pi}_1(M)\times[-\infty,\infty)\to[0,\infty]
$$
which assigns to the class $\alpha\in\widetilde{\pi}_1(M)$ and the
number $a\ge-\infty$ the number
$$
\widehat{C}(M,A;\alpha,a) := \inf\left\{c > 0\,\bigm|\, \sup
   \A_c(M,A;\alpha) > a \right\}.
$$
Here we use the convention that $\inf\emptyset = \infty$ and $\sup
\emptyset = -\infty$. For $a=-\infty$ we abbreviate
$$
\widehat{C}(M,A;\alpha) := \widehat{C}(M,A;\alpha,-\infty) =
\inf\left\{c > 0\,\bigm|\, \A_c(M,A;\alpha)\ne\emptyset\right\}.
$$
The latter quantity is independent of the $\om$-primitive $\lambda$
while $\widehat{C}(M,A;\alpha,a)$ does depend on this choice: the set
$\A_c(M,A;\alpha)=:\A_c^\lambda(M,A;\alpha)$ depends on $\lambda$, but
for two $\om$-primitives $\lambda,\lambda'$ we have $
\A_c^{\lambda'}(M,A;\alpha) = \A_c^\lambda(M,A;\alpha) -
\int_\alpha(\lambda'-\lambda).  $

\begin{prop} \label{P:exist}
   Let $\alpha\in\widetilde{\pi}_1(M)$ and $a\in\R$.  If
   $\widehat{C}(M,A;\alpha,a) < \infty$ then every compactly supported
   Hamiltonian $H$ on $S^1 \times M$ with $H|_{S^1\times
     A}\ge\widehat{C}(M,A;\alpha,a)$ has a $1$-periodic orbit in the
   homotopy class $\alpha$ with action $\Aa_H(x)\ge a$.  In
   particular,
   $$
   \widehat{C}(M,A;\alpha,a) \ge C(M,A;\alpha,a).
   $$
\end{prop}

\begin{proof}
   Assume first that $ \inf_{S^1\times A}H>\widehat{C}(M,A;\alpha,a).
   $ Then, by definition of $\widehat{C}(M,A;\alpha,a)$, there exist
   two real numbers $b$ and $c$ such that
   $$
   0 < c < \inf_{S^1\times A}H,\qquad a < b,\qquad
   b\in\A_c(M,A;\alpha).
   $$
   Hence, by definition of the set $\A_c(M,A;\alpha)$, the
   homomorphism
   $$
   T_\alpha^{[b,\infty);c}: \SHi^{[b,\infty)}(M;\alpha)
   \to\SHd^{[b,\infty);c}(M,A;\alpha)
   $$
   is nonzero. Now choose a sequence of Hamiltonian functions
   $H_i\in\Hh(M)$ such that $H_i$ converges to $H$ in the
   $\Cinf$-topology, $ b\notin\Ss(H_i;\alpha), $ and $ \inf_{S^1\times
     A}H_i > c $ for every $i$.  Then
   $H_i\in\Hh^{b,\infty}_c(M,A;\alpha)$ and so, by
   Proposition~\ref{P:factor}, the nonzero homomorphism
   $T_\alpha^{[b,\infty);c}$ factors through the Floer homology group
   $\HF^{[b,\infty)}(H_i;\alpha)$ for every $i$. Hence there exists a
   sequence of periodic orbits $x_i\in\Pp(H_i;\alpha)$ such that
   $\Aa_{H_i}(x_i)>b$. Passing to a converging subsequence we get a
   periodic orbit $x\in\Pp(H;\alpha)$ with $\Aa_H(x)\ge b>a$.  This
   proves the assertion in the case $\inf_{S^1\times
     A}H>\widehat{C}(M,A;\alpha,a)$.  If $\inf_{S^1\times
     A}H=\widehat{C}(M,A;\alpha,a)$ the result follows by another
   approximation argument.
\end{proof}

\begin{rem} \label{rmk:scale}
   Both relative capacities have the following rescaling property.
   Let $c>0$ and replace the symplectic form $\omega$ by $c\omega$ and
   the $\om$-primitive $\lambda$ by $c\lambda$.  Then
   $$
   \widehat{C}(M,A,c\lambda;\alpha,ca) =
   c\widehat{C}(M,A,\lambda;\alpha,a), \qquad
   C(M,A,c\lambda;\alpha,ca) = c C(M,A,\lambda;\alpha,a).
   $$
   To see this, note that the Hamiltonian function $\widetilde
   H:=cH$ has the same Hamiltonian vector field with respect to
   $\widetilde\om:=c\om$ as the Hamiltonian function $H$ with respect
   to $\om$, and that the symplectic action of a periodic orbit
   $x\in\Pp(\widetilde H,\widetilde\om;\alpha)=\Pp(H,\om;\alpha)$ with
   respect to $\widetilde H$ and $\widetilde\lambda:=c\lambda$ is
   equal to $c$ times the action with respect to $H$ and $\lambda$.
\end{rem}


\section{Computation of the capacities} 
\label{s:Prep} 

We are now in a position to compute in certain cases the relative
symplectic homology of the unit cotangent bundle $U^*X$ of a compact
connected Riemannian manifold $X$ without boundary.  We shall always
work with the Liouville form $\lambda_\can$ as a primitive of the
canonical symplectic form $\om_\can$.  We consider the following two
cases.
\begin{enumerate}
  \item[(T)] $X=\T^n=\R^n/\Z^n$ is the flat torus.
  \item[(N)] $X$ has negative sectional curvature.
\end{enumerate}
In either case we identify $\widetilde{\pi}_1(U^*X)$ with
$\widetilde{\pi}_1(X)$ and in the case of the torus we identify
$\widetilde{\pi}_1(\T^n)$ with $\Z^n$.  More precisely, we identify
$k\in\Z^n$ with the homotopy class of the loop $[0,1]\to\T^n:t\mapsto
tk+\Z^n$.


\subsection{The main results}\label{sb:UTn}

In this section we state the main results about the (relative)
symplectic homology of open subsets of cotangent bundles and show how
they can be used to establish Theorem~\ref{T:cotangent}~(i) and~(ii).
The subsequent sections are devoted to their proofs.  The first result
concerns the symplectic homology in the trivial homotopy class.

\begin{thm}\label{T:zero}
   Assume~$(T)$ or~$(N)$ and consider the trivial class
   $\alpha=0\in\widetilde{\pi}_1(X)$. Then, for $a,c>0$, we have
   $$
   \SHi^{[a,\infty)}(U^*X;0) \cong H_*(X;\Z_2),
   $$
   and
   $$
   \SHd^{[a,\infty);c}(U^*X, X;0) \cong
   \begin{cases}
      H_*(X;\Z_2), & \textnormal{if } a\le c, \\
      0, & \textnormal{if } c < a.
   \end{cases}
   $$
   Moreover, the homomorphism $ T^{[a,\infty);c}_0:
   \SHi^{[a,\infty)}(U^*X;0) \to\SHd^{[a,\infty);c}(U^*X,X;0) $ is
   an isomorphism for $0<a\le c$.  In particular, for every $a\in
   \mathbb{R}$ $$\widehat{C}(U^*X,X;0,a)=\max\{0,a\}.$$
\end{thm}

\begin{thm}\label{T:TN}
   Assume~$(T)$ or~$(N)$ and consider a nontrivial homotopy class
   $0\ne\alpha\in\widetilde{\pi}_1(X)$. Denote by $\ell>0$ the
   (unique) length of the geodesics in the class $\alpha$. Let
   $P:=\T^n$ in the case~$(T)$ and $P:=S^1$ in the case~$(N)$.  Then,
   for every $a \in \mathbb{R}$ and $c>0$ we have
   $$
   \SHi^{[a,\infty)}(U^*X;\alpha) \cong
   \begin{cases}
      0, & \textnormal{if } a < \ell,  \\
      H_*(P;\Z_2), & \textnormal{if } a\ge\ell.
   \end{cases}
   $$
   and
   $$
   \SHd^{[a,\infty);c}(U^*X,X;\alpha) \cong
   \begin{cases}
      H_*(P;\Z_2), & \textnormal{if } 0<a\le c, \\
      0, & \textnormal{if } a>c.
   \end{cases}
   $$
   Moreover, the homomorphism $ T^{[a,\infty);c}_\alpha:
   \SHi^{[a,\infty)}(U^*X;\alpha)
   \to\SHd^{[a,\infty);c}(U^*X,X;\alpha) $ is an isomorphism for
   $\ell\le a\le c$.  In particular, for every $a\in \mathbb{R}$
   $$\widehat{C}(U^*X,X;\alpha,a)=\max\{\ell,a\}.$$
\end{thm}

We are now in position to prove Theorem~\ref{T:cotangent}.
\begin{proof}[Proof of Theorem~\ref{T:cotangent}~(i) and~(ii)]
   Assume that $X$ satisfies~$(T)$ or~$(N)$.  If $\alpha=0$ we must
   prove that $C(U^*X,X;0,a)=\max\{0,a\}$ To see this, note that every
   compactly supported Hamiltonian function $H\in\Hh(U^*X)$ has a
   contractible periodic orbit $x$ with action $\Aa_H(x)=0$ and hence
   $C(U^*X,X;0,a)=0$ whenever $a\le 0$.  If $a>0$ then
   Theorem~\ref{T:zero} asserts that $\widehat{C}(U^*X,X;0,a)=a$ and
   hence
   $$
   a = \widehat{C}(U^*X,X;0,a) \ge C(U^*X,X;0,a) \ge a.
   $$
   Here the middle inequality follows from
   Proposition~\ref{P:exist}. To prove the right-hand inequality let
   $0<\delta<a$ and choose any Hamiltonian function $H=H(p)$ that
   depends only on the momenta variables and satisfies $\max
   H=a-\delta$.  Then every contractible periodic orbit $x\in\Pp(H;0)$
   is (up to parametrization) a contractible geodesic and hence, since
   $X$ satisfies~$(T)$ or~$(N)$, is constant and has action
   $\Aa_H(x)=H(x)=a-\delta$.  This shows that $C(U^*X,X;0,a)\ge
   a-\delta$ for every $\delta>0$.  Thus we have proved that
   $C(U^*X,X;0,a)=\max\{0,a\}$ as claimed.
   
   Now assume $\alpha\ne 0$ and abbreviate $\ell:=\ell(\gamma_\alpha)$
   in the case~(N) and $\ell:=|k|$ in the case~(T) with
   $\alpha=k\in\Z^n$.  Then Theorem~\ref{T:TN} asserts that
   $\widehat{C}(U^*X,X;\alpha,a) = \max\{\ell,a\}$ and hence
   $$
   \max\{\ell,a\} = \widehat{C}(U^*X,X;\alpha,a) \ge
   C(U^*X,X;\alpha,a) \ge \max\{\ell,a\}
   $$
   for every real number $a$.  Again the middle inequality follows
   from Proposition~\ref{P:exist} and the rightmost inequality from an
   explicit construction of a Hamiltonian function.  Namely, for any
   $\delta>0$ choose a compactly supported function $f:[0,1)\to\R$
   such that
   $$
   f(r) =
   \begin{cases}
      m-\delta, &\mbox{ for }r\mbox{ near }0, \\
      0, &\mbox{ for }r\mbox{ near }1,
   \end{cases}     
   $$
   where $m:=\max\{\ell,a\}$, and
   $$
   f(r) < (1-r)m,\qquad -m < f'(r) \le 0
   $$
   for every $r$.  Now consider the compactly supported Hamiltonian
   function $H:=f(|p|)$ on $U^*X$.  Its $1$-periodic solutions are
   reparametrized closed geodesics. The sphere bundle $|p|=r$ contains
   a periodic orbit $x$ in the class $\alpha$ if and only if
   $f'(r)=-\ell$ and the action of this periodic orbit is
   $$
   \Aa_H(x) = f(r) - rf'(r) < f(r) + rm < m.
   $$
   (See Lemma~\ref{le:TN} below.)  If $a\le\ell$ then $f'(r)>-\ell$
   for all $r$, hence there is no $1$-periodic solution of length
   $\ell$, and hence none in the class $\alpha$.  If $\ell\le a$ then
   every $1$-periodic solution has action $\Aa_H(x)<a$.  In either
   case there is no $1$-periodic orbit in the class $\alpha$ with
   action at least $a$, and hence $C(U^*X,X;0,a)\ge m-\delta$.  Since
   this holds for every $\delta>0$ we obtain $C(U^*X,X;0,a)\ge m$ as
   claimed.
\end{proof}


\subsection{Morse--Bott theory in Floer homology}\label{sb:MB}

Let us return to the general setting of Section~\ref{sb:setting} where
$(\overline M,\om)$ is a compact connected symplectic manifold with
convex boundary, $\om=d\lambda$ is an exact symplectic form,
$M=\overline M\setminus\p\overline M$, $\Hh=\Hh(M)$ denotes the space
of compactly supported functions on $S^1\times M$, and $\Jj$ denotes
the space of $1$-periodic $\om$-compatible almost complex structures
$J_t=J_{t+1}$ on $M$.

A subset $P\subset\Pp(H)$ is called a {\bf Morse--Bott manifold of
  periodic orbits} if the set $C_0:=\left\{x(0)\,|\,x\in P\right\} $
is a compact submanifold of $M$ and $T_{x_0}C_0 =
\ker\,(D\psi_1(x_0)-\one)$ for every $x_0\in C_0$.

\begin{rem}\label{rmk:MB}\rm
   The Morse--Bott condition can be reformulated as follows.  Firstly,
   a subset $P\subset\Pp(H)$ is a compact submanifold of the loop
   space $LM$ if and only if the set $C_0=\left\{x(0)\,|\,x\in
      P\right\}$ is a compact submanifold of $M$.  Secondly, for every
   $x\in P$ the kernel of the linear map $D\psi_1(x(0))-\one$ on
   $T_{x(0)}M$ is isomorphic to the space of periodic solutions of the
   following {\it linearized Hamiltonian differential equation} for
   vector fields $\xi(t)\in T_{x(t)}M$ along $x$:
\begin{equation}\label{eq:linear-H}
   \Nabla{\dot{x}}{\xi} = \Nabla{\xi}X_H(x),\qquad \xi(t+1)=\xi(t), 
\end{equation}
where $\nabla$ stands for the Levi-Civita connection of the metric
$\omega(\cdot, J \cdot)$. To see this just note that $\xi$
satisfies~(\ref{eq:linear-H}) if and only if
$\xi(t)=D\psi_t(x(0))\xi(0)$ for all $t$.  Note that every tangent
vector of $P$ is a solution of~(\ref{eq:linear-H}).  The Morse--Bott
condition can now be expressed in the form that $P$ is a compact
submanifold of $LM$ and
$$
T_xP = \left\{\xi\in\Cinf(S^1,x^*TM)\,|\,\xi\mbox{ satisfies }
   (\ref{eq:linear-H})\right\}.
$$
We emphasize that the Hessian of the symplectic action functional
$\Aa_H:LM\to\R$ at a critical point $x$ is the linear operator $
\xi\mapsto \Nabla{\xi}(\textnormal{grad}H)-(\Nabla{\xi}J)\dot
x-J\Nabla{\dot{x}}\xi = J(\Nabla{\xi}X_H-\Nabla{\dot{x}}\xi) $ on
$\Cinf(S^1,x^*TM)$ (equipped with the $L^2$-inner product).  Hence the
space of solutions of~(\ref{eq:linear-H}) is the kernel of the Hessian
of $\Aa_H$ at $x$, and the Morse--Bott condition asserts that the
kernel of the Hessian agrees with the tangent space of the critical
manifold $P$.
\end{rem}

\begin{thm}\label{T:global}
   Let $-\infty\le a<b\le\infty$, $\alpha\in\widetilde{\pi}_1(M)$, and
   $H\in\Hh^{a,b}(M;\alpha)$. Suppose that the set
   $P:=\left\{x\in\Pp(H;\alpha)\,|\,a<\Aa_H(x)<b\right\}$ is a
   connected Morse--Bott manifold of periodic orbits.  Then
   $\HF^{[a,b)}(H;\alpha) \cong H_*(P;\Z_2)$.
\end{thm}

This is a version of Pozniak's theorem~\cite{P} which was originally
proved in the context of Floer homology for Lagrangian intersections.
In this section we explain the reduction of this theorem to Pozniak's
original one.

In order to reformulate Floer homology in the Lagrangian setting let
us consider the symplectic manifold
$$
\widetilde M:=M\times M,\qquad \widetilde\om := \om\oplus(-\om) =
d\widetilde\lambda,\qquad \widetilde\lambda :=
\lambda\oplus(-\lambda).
$$
Since $\widetilde\om$ is exact there are no nonconstant
holomorphic spheres in $\widetilde M$, for any
$\widetilde\om$-compatible almost complex structure $\widetilde
J\in\Jj(\widetilde M,\widetilde\om)$.  Since $\widetilde\lambda$
vanishes on the diagonal $\Delta\subset M\times M$, there are also no
nonconstant holomorphic disks with boundary in $\Delta$.  Hence the
standard theory of Floer homology for Lagrangian intersections applies
as in Floer's original work~\cite{F1,F2,F3}.  Given $H\in\Hh(M)$
define $\widetilde H_t:\widetilde M\to\R$ by
$$
\widetilde H_t(x_0,x_1):=H_t(x_0)+H_{1-t}(x_1).
$$
The Hamiltonian isotopy generated by $\widetilde H_t$ with respect
to $\widetilde\om$ is given by $ \widetilde\psi_t(x_0,x_1) =
\left(\psi_t(x_0),\psi_{1-t}\circ\psi_1^{-1}(x_1)\right).  $ Given
$J\in\Jj(M)$ define $\widetilde J_t\in\Jj(\widetilde M,\widetilde\om)$
by
\begin{equation}\label{eq:Jtilde}
   \widetilde J_t
   := \psi_t^*J_t\times(-(\psi_{1-t}\circ\psi_1^{-1})^*J_{1-t})
   = \left(\widetilde\psi_t\right)^*(J_t\times(-J_{1-t}))
\end{equation}
for $0\le t\le1/2$. Given $u:\R\times S^1\to M$ define $\widetilde
u:\R\times[0,1/2]\to\widetilde M$ by
$$
\widetilde u(s,t) :=
\left(\psi_t^{-1}(u(s,t)),\psi_1\circ\psi_{1-t}^{-1}(u(s,1-t))\right)
= \left(\widetilde\psi_t\right)^{-1}(u(s,t),u(s,1-t)).
$$
Then $u$ satisfies~(\ref{eq:floer}) if and only if $\widetilde u$
satisfies the Lagrangian boundary value problem
\begin{equation}\label{eq:floer-L}
   \p_s\widetilde u 
   + \widetilde J_t(\widetilde u)\p_t\widetilde u = 0,\qquad
   \widetilde u(s,0)\in\Delta,\qquad
   \widetilde u(s,1/2)\in{\rm graph}(\psi_1).
\end{equation}
It satisfies~(\ref{eq:limits}) if and only if $\widetilde u$ satisfies
\begin{equation}\label{eq:limits-L}
   \lim_{s\to\pm\infty}\widetilde u(s,t) = \widetilde x^\pm,\qquad
   \lim_{s\to\pm\infty}\p_s\widetilde u(s,t) = 0,
\end{equation}
where $ \widetilde x^\pm :=(x^\pm(0),x^\pm(0))\in\Delta\cap{\rm
  graph}(\psi_1).  $ The solutions of~(\ref{eq:floer-L}) can be
interpreted as the gradient flow lines of the action functional
$$
\widetilde\Aa(\widetilde x) :=
-\int_0^{1/2}\widetilde\lambda(\dot{\widetilde x}(t))\,dt
$$
on the space $\widetilde\Pp$ of paths $\widetilde
x:[0,1/2]\to\widetilde M$ with endpoints $\tilde x(0)\in\Delta$ and
$\tilde x(1/2)\in{\rm graph}(\psi_1)$ (with respect to the $L^2$
metric determined by $\widetilde J$).  Note that the composition of
$\widetilde\Aa$ with the map $LM\to\widetilde\Pp:x\mapsto\widetilde
x$, given by $ \widetilde x(t)
:=(\psi_t^{-1}(x(t)),\psi_1\circ\psi_{1-t}^{-1}(x(1-t))), $ agrees
with $\Aa_H$. Note also that this map induces a bijection
$$
\widetilde\pi_1(M)=\pi_0(LM)\to\pi_0(\widetilde\Pp):
\alpha\mapsto\widetilde\alpha.
$$
Hence the solutions of~(\ref{eq:floer-L}) can be used to define the
Floer homology groups of the pair $(\Delta,{\rm graph}(\psi_1))$ of
Lagrangian submanifolds of $(\widetilde M,\widetilde\om)$.  Moreover,
the Floer homology groups defined by the solutions
of~(\ref{eq:floer-L}) are independent of the choice of the (regular)
almost complex structure $\widetilde J_t$ used to define them.  More
precisely, denote by $\widetilde\Jj$ the space of smooth functions
$[0,1/2]\to\Jj(\widetilde M,\widetilde\om):t\mapsto\widetilde J_t$
such that $\widetilde J_t=J\times (-J)$ near the boundary of
$\widetilde M$, where $J\in\Jj(M,\om)$ is convex.  Given a Hamiltonian
$H\in\Hh$ that satisfies~$(H)$ denote by $\TJreg(H)$ the set of all
almost complex structures $\widetilde J\in\widetilde\Jj$ such that
every finite energy solution $\widetilde u$ of~(\ref{eq:floer-L}) is
regular in the sense that the linearized operator along $\widetilde u$
is surjective.  Then the solutions of~(\ref{eq:floer-L}) give rise to
Lagrangian Floer homology groups $\HF^{[a,b)}(\Delta,{\rm
  graph}(\psi_1);\widetilde J,\widetilde\alpha)$.  Moreover, it
follows as in~\cite{F1,F2,F3} (and as outlined above) that these Floer
homology groups are independent of the almost complex structure
$\widetilde J\in\TJreg(H)$ used to define them.  Note that if
$J\in\Jreg(H)$ and $\widetilde J$ is given by~(\ref{eq:Jtilde}) then
$\widetilde J\in\TJreg(H)$.  Hence there is a natural isomorphism
$$
\HF^{[a,b)}(H;\alpha) \cong\HF^{[a,b)}(\Delta,{\rm
  graph}(\psi_1);\widetilde\alpha)
$$
for every $\alpha\in\widetilde{\pi}_1(M)$ and every Hamiltonian
$H\in\Hh^{a,b}(M;\alpha)$, where $\widetilde\alpha$ is the image of
$\alpha$ under the above homomorphism
$\widetilde\pi_1(M)=\pi_0(LM)\to\pi_0(\widetilde\Pp)$.  The advantage
of the Lagrangian approach in the present context is that we can use
any (regular) family of $\widetilde\om$-compatible almost complex
structures $\{\widetilde J_t\}_{0\le t\le 1/2}$ to define the Floer
homology groups, and are not restricted to those arising from periodic
families of almost complex structures on $M$ via~(\ref{eq:Jtilde}).
Hence we can apply the results of Pozniak.

Let $\widetilde J\in\widetilde\Jj$ and $H\in\Hh$ and denote by
$\R\to\Diff(M,\om):t\mapsto\psi_t$ the Hamiltonian isotopy generated
by $H$. The {\bf graph} of a loop $x\in LM$ is the set
$$
\Gamma(x) := \left\{\left(t,\psi_t^{-1}(x(t)),
      \psi_1\circ\psi_{1-t}^{-1}(x(1-t))\right)\,|\, 0\le
   t\le1/2\right\} \subset [0,1/2]\times\widetilde M.
$$
For a subset $P\subset LM$ we write $ \Gamma(P):=\bigcup_{x\in
  P}\Gamma(x), $ and for a map $\widetilde
u:\R\times[0,1/2]\to\widetilde M$ we write
$$
\Gamma(\widetilde u):=\left\{(t,\widetilde u(s,t))\,|\,
   s\in\R,\,0\le t\le1/2\right\}.
$$
A subset $P\subset\Pp(H)$ is called a {\bf $\widetilde J$-isolated
  periodic set} if there exists an open neighbourhood
$U\subset[0,1/2]\times\widetilde M$ of $\Gamma(P)$ such that the
following holds:
\begin{enumerate}
  \item[(P1)] The closure $\overline{U}$ is a compact subset of
   $[0,1/2]\times\widetilde M$.
  \item[(P2)] If $\widetilde u:\R\times[0,1/2]\to\widetilde M$ is a
   finite energy solution of~(\ref{eq:floer-L}) with
   $\Gamma(\widetilde u)\subset \overline{U}$ then there exists an
   $x\in P$ such that $u(s,t)=x(t)$ for every $(s,t)\in\R^2$.
\end{enumerate}
An open neighbourhood $U\subset[0,1/2]\times\widetilde M$ of
$\Gamma(P)$ that satisfies~$(P1)$ and~$(P2)$ is called {\bf
  $\widetilde J$-isolating}.  Note that every $\widetilde J$-isolated
periodic set is compact.

\begin{lem}\label{le:MB}
   Let $H\in\Hh(M)$. Then every Morse--Bott manifold $P\subset\Pp(H)$
   of periodic orbits is a $\widetilde J$-isolated periodic set for
   every almost complex structure $\widetilde J\in\widetilde\Jj$.
\end{lem}

\begin{proof} 
   We may assume without loss of generality that $P$ is connected.
   Let $\widetilde d_t$ denote the distance function of the metric
   $\inner{\cdot}{\cdot}_t:=\widetilde\om(\cdot,\widetilde J_t\cdot)$
   and consider the open set
   $$
   U:=\left\{(t,\widetilde x)\,\Big|\, 0\le t\le 1/2,\,\widetilde
      x\in \widetilde M,\, \sup_{y\in P}\widetilde d_t(\widetilde
      x,(y(1/2-t),y(1/2+t)))<\eps \right\} \subset
   [0,1/2]\times\widetilde M.
   $$
   Let $\eps>0$ be sufficiently small. Then since $C_0$ is an
   isolated fixed point set for $\psi_1$ it follows that every
   $x\in\Pp(H)$ with $\Gamma(x)\subset\overline{U}$ is an element of
   $P$.  Now the set $U$ satisfies~$(P2)$ because every finite energy
   solution $\widetilde u:\R\times[0,1/2]\to\widetilde M$
   of~(\ref{eq:floer-L}) with $\Gamma(\widetilde u)\subset
   \overline{U}$ is asymptotic to the set $P$ as $s\to\pm\infty$.
   Since $\Aa_{\widetilde H}=\Aa_H$ is constant along $P$ it follows
   that every such solution $\widetilde u$ has energy $E(\widetilde
   u)=0$ and hence has the form $\widetilde u(s,t)=\widetilde
   x(t)=(x(1/2-t),x(1/2+t)$ for some $x\in\Pp(H)$.
\end{proof}

\begin{lem}\label{le:isolated}
   Let $H\in\Hh$ and $J\in\widetilde\Jj$. Suppose that
   $P\subset\Pp(H)$ is a $\widetilde J$-isolated periodic set and
   $U\subset[0,1/2]\times\widetilde M$ is a $\widetilde J$-isolating
   neighbourhood of $\Gamma(P)$.  Then there exist a compact
   neighbourhood $V\subset U$ of $\Gamma(P)$ and a constant $\delta>0$
   such that the following holds.  If $\R\to\Hh(M):s\mapsto H_s$ and
   $\R \to \Jj(\widetilde M):s\mapsto\widetilde J_s$ are smooth
   homotopies such that
   $$
   \|H_s-H\|_{C^2} + \|\widetilde J_s-\widetilde J\|_{C^1} +
   \|\p_sH_s\|_{C^2} + \|\p_s\widetilde J_s\|_{C^1} < \delta
   $$
   and $\p_sH_s=0$ and $\p_s\widetilde J_s=0$ for $|s|\ge 1$ then
   every finite energy solution $\widetilde
   u:\R\times[0,1/2]\to\widetilde M$ of~(\ref{eq:floer-L}) with
   $(\widetilde H,\widetilde J)$ replaced by $(\widetilde
   H_s,\widetilde J_s)$ satisfies
   $$
   \Gamma(\widetilde u)\subset \overline{U} \qquad\IMP\qquad
   \Gamma(\tilde u)\subset V.
   $$
\end{lem}

\begin{proof}
   Suppose, by contradiction, that there exist sequences
   $$
   \R\to\Hh:s\mapsto H^\nu_s,\qquad
   \R\to\widetilde\Jj:s\mapsto\widetilde J^\nu_s,\qquad \widetilde
   u^\nu:\R\times[0,1/2]\to\widetilde M,
   $$
   and $(s^\nu,t^\nu)\in\R\times[0,1/2]$ such that the following
   holds:
\begin{enumerate}
  \item[(i)] $ \lim_{\nu\to\infty}\sup_{s\in\R}\left(
      \|H^\nu_s-H_s\|_{C^2} + \|\p_sH^\nu_s\|_{C^2} + \|\widetilde
      J^\nu_s-\widetilde J\|_{C^1} + \|\p_s\widetilde J^\nu_s\|_{C^1}
   \right) = 0.  $
  \item[(ii)] $\p_sH^\nu_s=0$ and $\p_s\widetilde J^\nu_s=0$ for
   $|s|\ge 1$.
  \item[(iii)] $\widetilde u^\nu$ is a finite energy solution
   of~(\ref{eq:floer-L}) with $(\widetilde H,\widetilde J)$ replaced
   by $(\widetilde H^\nu_s,\widetilde J^\nu_s)$.
  \item[(iv)] $\Gamma(\widetilde u_\nu)\subset \overline{U}$ and
   $\lim_{\nu\to\infty}\widetilde u^\nu(s^\nu,t^\nu)\in\p U$.
\end{enumerate}
Since there are no nonconstant $\widetilde J_t$-holomorphic spheres
in $\widetilde M$ and no nonconstant $\widetilde J_t$-holomorphic
disks with boundary in $\Delta$, the first derivatives of the
functions $\widetilde u_\nu$ are uniformly bounded.  Hence, by
Floer--Gromov compactness~\cite{F,G,MS,Sa}, there exists a
subsequence, still denoted by $\widetilde u^\nu$, such that the
shifted sequence $\widetilde u^\nu(s^\nu+s,t)$ converges in the
$C^1$-topology on compact sets to a finite energy solution $\widetilde
u:\R\times[0,1/2]\to\widetilde M$ of~(\ref{eq:floer-L}) such that
$\Gamma(\widetilde u)\subset\overline{U}$.  By taking a further
subsequence we may assume that $t^\nu\to t$ and hence $\widetilde
u(0,t)=\lim_{\nu\to\infty}\widetilde u^\nu(s^\nu,t^\nu)$ satisfies
$(t,\widetilde u(0,t))\in\p U\subset[0,1/2]\times(M\setminus U)$.
This contradicts~$(P2)$.
\end{proof}

Lemma~\ref{le:isolated} enables us to define the {\bf local Floer
  homology} $\HF^\loc(H;P)$ of a $\widetilde J$-isolated periodic set
$P\subset\Pp(H)$ as follows.  Choose a $\widetilde J$-isolating
neighbourhood $U\subset S^1\times M$ of $\Gamma(P)$, let $\delta>0$ be
as in Lemma~\ref{le:isolated}, choose a Hamiltonian function $H'$ such
that all periodic solutions $x\in\Pp(H')$ are nondegenerate and
$\|H'-H\|_{C^2}<\delta/4$, and choose a regular almost complex
structure $\widetilde J'\in\TJreg(H')$ such that $\|\widetilde
J'-\widetilde J\|_{C^1} < \delta/4$.  Then, by
Lemma~\ref{le:isolated}, all the Floer connecting orbits of
$(\widetilde H',\widetilde J')$ (i.e. solutions $\widetilde u'$
of~(\ref{eq:floer-L}) with $(\widetilde H,\widetilde J)$ replaced by
$(\widetilde H',\widetilde J')$) in $\overline{U}$ are actually
contained in~$V$. Denote the set of local periodic orbits of $H'$ near
$P$ by
$$
\Pp(H';U):=\left\{x'\in\Pp(H')\,|\,\Gamma(x')\subset U\right\}
$$
and consider the local Floer chain complex
$$
\CF^\loc(H';U) := \bigoplus_{x'\in\Pp(H';U)} \Z_2 x'.
$$
The boundary operator $\p^{H',\widetilde
  J';U}:\CF^\loc(H';U)\to\CF^\loc(H';U)$ is defined by counting the
index-$1$ solutions $\widetilde u'$ of~(\ref{eq:floer-L}), with
$(\widetilde H,\widetilde J)$ replaced by $(\widetilde H',\widetilde
J')$, such that $\Gamma(\widetilde u')\subset U$.  Since these
solutions can never converge to the boundary of $U$ it follows that
$\p^{H',J';U}$ is indeed a boundary operator and the local Floer
homology is defined by
$$
\HF^\loc(H',\widetilde J';U) :=
H_*(\CF^\loc(H';U),\p^{H',\widetilde J';U}).
$$
The same arguments as in Floer's original theory~\cite{F1,F2,F3}
now show that this local Floer homology is independent (up to natural
isomorphisms) of the isolating neighbourhood $U$, and of the
perturbations $H'$ and $\widetilde J'$ used to define it.  We write
$$
\HF^\loc(H;P) := \HF^\loc(H',\widetilde J';U)
$$
for the local Floer homology in a $\widetilde J$-isolating
neighbourhood $U$ of $\Gamma(P)$.  Strictly speaking, this is a {\it
  connected simple system} in the sense of Conley, namely a small
category whose objects are the triples $(H',\widetilde J';U)$ of local
perturbations and whose morphisms are the canonical (unique)
isomorphisms between the local Floer homologies
$\HF^\loc(H'_0,\widetilde J'_0,U_0)$ and $\HF^\loc(H'_1,\widetilde
J'_1,U_1)$.  The details of this construction were carried out by
Pozniak~\cite{P} in the context of Lagrangian intersections.

\begin{thm}{\bf (Pozniak~\cite{P})}\label{thm:pozniak}
   Let $H\in\Hh(M)$ and suppose that $P\subset\Pp(H)$ is a connected
   Morse--Bott manifold of periodic orbits.  Then $\HF^\loc(H;P) \cong
   H_*(P;\Z_2)$.
\end{thm}

\begin{proof}
   The local Floer homology of $H$ near $P$ is isomorphic to the local
   Floer homology of the pair of Lagrangian submanifolds
   $L_0:=\Delta\subset M\times M$ and $L_1:={\rm graph}(\psi_1)\subset
   M\times M$ near their clean intersection
   $\Lambda:=\{(x(0),x(0))\,|\,x\in P\}$.  Hence,
   by~\cite[Theorem~3.4.11]{P}, it is isomorphic to
   $H_*(\Lambda;\Z_2)\cong H_*(P;\Z_2)$.
\end{proof}

\begin{proof}[Proof of Theorem~\ref{T:global}]
   Fix a $1$-periodic almost complex structure $J\in\Jj$ and let
   $\widetilde J\in\widetilde\Jj$ be given by~(\ref{eq:Jtilde}).
   Then, by Lemma~\ref{le:MB}, $P$ is a $\widetilde J$-isolated
   periodic set.  Let $U$ be a $\widetilde J$-isolating neighbourhood
   of $\Gamma(P)$ and choose a sequence of regular perturbations
   $(H^\nu,J^\nu)$ that agree with $(H,J)$ in some neighbourhood of
   $\p\overline M$ and converge to $(H,J)$ in the $C^2$-norm. We claim
   that, for $\nu$ sufficiently large, all the Floer connecting orbits
   (i.e solutions of~(\ref{eq:floer})) for the pair $(H^\nu,J^\nu)$ in
   the homotopy class $\alpha$ and the action interval $[a,b]$ are
   contained in $U$.  Otherwise, there has to be a sequence $u^\nu$ of
   such connecting orbits passing through $M\setminus U$ and we can
   argue as in the proof of Lemma~\ref{le:isolated} that, in the limit
   $\nu\to\infty$, there must be a finite energy solution
   of~(\ref{eq:floer}) for the pair $(H,J)$ in the homotopy class
   $\alpha$ and the action interval $[a,b]$ that passes through
   $M\setminus U$.  However, every such connecting orbit has the form
   $u(s,t)=x(t)$ for some $x\in P$ and so is contained in $U$.  This
   contradiction proves the claim.  Hence
   $\HF^{[a,b)}(H;\alpha)\cong\HF^\loc(H;P)$, and hence the result
   follows from Theorem~\ref{thm:pozniak}.
\end{proof}


\subsection{The main example}\label{sb:example} 

In this section we consider the case where $M=U^*X$ is the open unit
cotangent bundle of a compact connected Riemannian $n$-manifold $X$
without boundary that satisfies either~$(T)$ (i.e. $X$ is a flat
torus) or~$(N)$ (i.e. $X$ has negative sectional curvature).  We shall
use the metric to identify the tangent bundle with the cotangent
bundle and denote a point in $U^*X$ by $x=(q,p)$ where $q\in X$ and
$p\in T_qX$.  Let $H:U^*X\to\R$ be a compactly supported Hamiltonian
function of the form
$$
H(q,p) = f(|p|),
$$
where $f:\R\to\R$ is a smooth function such that $f(r)=f(-r)$. The
corresponding Hamiltonian differential equation has the form
\begin{equation}\label{eq:fham}
     \dot q = \frac{f'(|p|)}{|p|}p,\qquad \Nabla{\dot{q}}{p}=0.
\end{equation}
Here $\nabla$ denotes the Levi-Civita connection.  Since $|p|$ is
constant it follows that $q$ is a geodesic for every solution $(q,p)$
of~(\ref{eq:fham}).  Moreover, since $f'(0)=0$, the zero section
$\{p=0\}$ consists of constant solutions.  There are other constant
solutions $x(t)\equiv(q,p)$ whenever $f'(|p|)=0$ but these will not be
relevant in the context of the present paper.

\begin{lem}\label{le:zero}
   The set $P_0 := \left\{x=(p,q):S^1\to T^*X\,|\, \dot
      q\equiv0,\,p\equiv0\right\}$ is a Morse--Bott manifold of
   periodic orbits for $H$ if and only if $f''(0)\ne 0$.
\end{lem}

\begin{proof} 
   Since $f$ is even there exists a smooth function $h:\R\to\R$ such
   that $f(r)=h(r^2)/2$.  Then $h'(r^2)=f'(r)/r$ and $h'(0)=f''(0)$.
   So equation~(\ref{eq:fham}) reads
   $$
   \dot q = h'(|p|^2)p,\qquad \Nabla{\dot{q}}p=0.
   $$
   By Remark~\ref{rmk:MB}, $P_0$ is a Morse--Bott manifold if and
   only if the space of periodic solutions of the linearized equation
   is equal to the tangent space of $P_0$ for every $x=(q,p)\in P_0$.
   This means that the space of periodic solutions of the linearized
   equation has the same dimension as $P_0$.  Now the linearized
   equation at a constant solution with $p\equiv0$ has the form
   $$
   \frac{d}{dt}\hat q = h'(0)\hat p = f''(0)\hat p,\qquad
   \frac{d}{dt}\hat p = 0.
   $$
   If $f''(0)\ne 0$ then the space of periodic solutions of this
   equation has dimension $n=\dim\,P_0$ and if $f''(0)=0$ it has
   dimension $2n$.
\end{proof}

Lemma~\ref{le:zero} continues to hold for any compact Riemannian
manifold $X$.  However, in the case~$(T)$ or~$(N)$, every nonconstant
closed geodesic is not contractible.  Let us now consider a nonzero
homotopy class $0\ne\alpha\in\widetilde{\pi}_1(X)$ and denote by
$\ell$ the (unique) length of the closed geodesics in the class
$\alpha$.  The space of solutions to equation~\eqref{eq:fham} that
represent the class $\alpha$ consists of $\left\{\big(q(t),
   p(t)\big)\right\}$, where
\begin{equation} \label{eq:P}
   \begin{cases}
      q(t) \textnormal{ is a geodesic in the class } \alpha,
      \textnormal{ parametrized so that } |\dot{q}| \equiv \ell. \\
      p(t) = \pm \frac{r}{\ell} \dot{q}(t), \textnormal{ where } r>0
      \textnormal{ is such that } f'(r)=\pm \ell.
   \end{cases}
\end{equation}
Given $r>0$ with $f'(r)=\pm \ell$ we denote
$$
P^\pm(r,\alpha) := \left\{x=(q,p):S^1\to U^*X\, \Big| \, p(t), q(t)
   \textnormal{ satisfy}~\eqref{eq:P}\right\}.
$$
In the case~$(T)$ the space $P^\pm(r,\alpha)$ is diffeomorphic to
$X$ and in the case~$(N)$ it is diffeomorphic to $S^1$.

\begin{lem}\label{le:TN}
   Assume $X$ satisfies~$(T)$ or~$(N)$ and let $\alpha\ne 0$, $\ell$,
   $r > 0$ with $f'(r) = \pm \ell$, and $P^\pm(r,\alpha)$ be as above.
   Then $P^\pm(r,\alpha)$ is a Morse--Bott manifold of periodic orbits
   for $H$ if and only if $f''(r)\ne 0.$ Moreover, $ \Aa_H(x) = f(r) -
   rf'(r) = f(r)\mp r\ell $ for every $x\in P^\pm(r,\alpha)$.
\end{lem}

\begin{proof} 
   As in the proof of Lemma~\ref{le:zero}, it follows from
   Remark~\ref{rmk:MB} that $P^\pm(r,\alpha)$ is a Morse--Bott
   manifold if and only if the space of periodic solutions of the
   linearized equation has the same dimension as $P^\pm(r,\alpha)$,
   namely $n$ in the case~$(T)$ and $1$ in the case~$(N)$.  We begin
   by linearizing equation~(\ref{eq:fham}).  Given a path $x:\R\to
   LU^*X:s\mapsto x_s=(q_s,p_s)$, we represent a variation of $x$ by a
   pair $\hat x=(\hat q,\hat p)$ of periodic vector fields along $q$
   via $\hat q:=\p_s q_s$ and $\hat p:=\Nabla{s}p_s$.  Since
   $\p_s|p|=|p|^{-1}\inner{p}{\Nabla{s}p}$ the linearized equation has
   the form
   \begin{equation}\label{eq:jacobi}
      \Nabla{t}\hat q = \frac{f'(r)}{r}
      \left(\hat p-\INNER{\frac{p}{r}}{\hat p}\frac{p}{r}\right)
      + f''(r)\INNER{\frac{p}{r}}{\hat p}\frac{p}{r},\qquad
      \Nabla{t}\hat p + \frac{f'(r)}{r}R(\hat q,p)p = 0.
   \end{equation}
   Here $R\in\Om^2(X,{\rm End}(TX))$ denotes the Riemann curvature
   tensor.  Note that in the case $f(r)=r^2/2$ we have $\Nabla{t}\hat
   q=\hat p$ and so~(\ref{eq:jacobi}) is equivalent to the standard
   Jacobi equation $\Nabla{t}\Nabla{t}\hat q+R(\hat q,p)q=0$.
   
   The periodic solutions of equation~(\ref{eq:jacobi}) form the
   kernel of the Hessian of the symplectic action (see
   Remark~\ref{rmk:MB}).  Taking the pointwise inner product of the
   second equation in~(\ref{eq:jacobi}) with $p$ and using
   $\Nabla{t}p=0$ we find that $\INNER{p}{\hat p}$ is constant.
   Hence, taking the $L^2$-inner product of the first equation
   in~(\ref{eq:jacobi}) with $p$ and using the fact that $|p|=r$, we
   find that every periodic solution of~(\ref{eq:jacobi}) satisfies
   $$
   f''(r)\INNER{p}{\hat p} \equiv 0
   $$
   Moreover, taking the $L^2$-inner product of the second equation
   in~(\ref{eq:jacobi}) with $\hat q$ and using the first equation we
   find that every periodic solution of~(\ref{eq:jacobi}) satisfies
   $$
   \int_0^1\left(|\hat p|^2 -\INNER{\frac{p}{r}}{\hat p}^2 -
      \INNER{R(\hat q,p)p}{\hat q}\right)\,dt = 0.
   $$
   Now suppose that $X$ has nonpositive sectional curvature and
   that $f''(r)\ne 0$.  Then $\hat p=0$ and $\Nabla{t}\hat q=0$ for
   every periodic solution of~(\ref{eq:jacobi}).  Hence the space of
   periodic solutions of~(\ref{eq:jacobi}) has dimension $n$ in the
   torus case (namely $\hat q$ is uniquely determined by $\hat q(0)\in
   T_{q(0)}X$) and has dimension one in the negative curvature case
   (namely, $\hat q$ is a scalar multiple of $p$).  In both cases it
   follows that the kernel of the Hessian of the symplectic action at
   every point $x\in P^\pm(r,\alpha)$ has the same dimension as
   $P^\pm(r,\alpha)$ and hence is equal to the tangent space of
   $P^\pm(r,\alpha)$ at $x$.  This is equivalent to the Morse--Bott
   nondegeneracy condition (Remark~\ref{rmk:MB}).  If, on the other
   hand, $f''(0)=0$ then the dimension of the space of periodic
   solutions of~(\ref{eq:jacobi}) is $n+1$ in the torus case and is
   $2$ in the negative curvature case.
\end{proof}


\subsection{Proof of Theorem~\ref{T:zero}}\label{sb:zero}

\begin{figure}[htp]
   \centerline{\psfig{figure=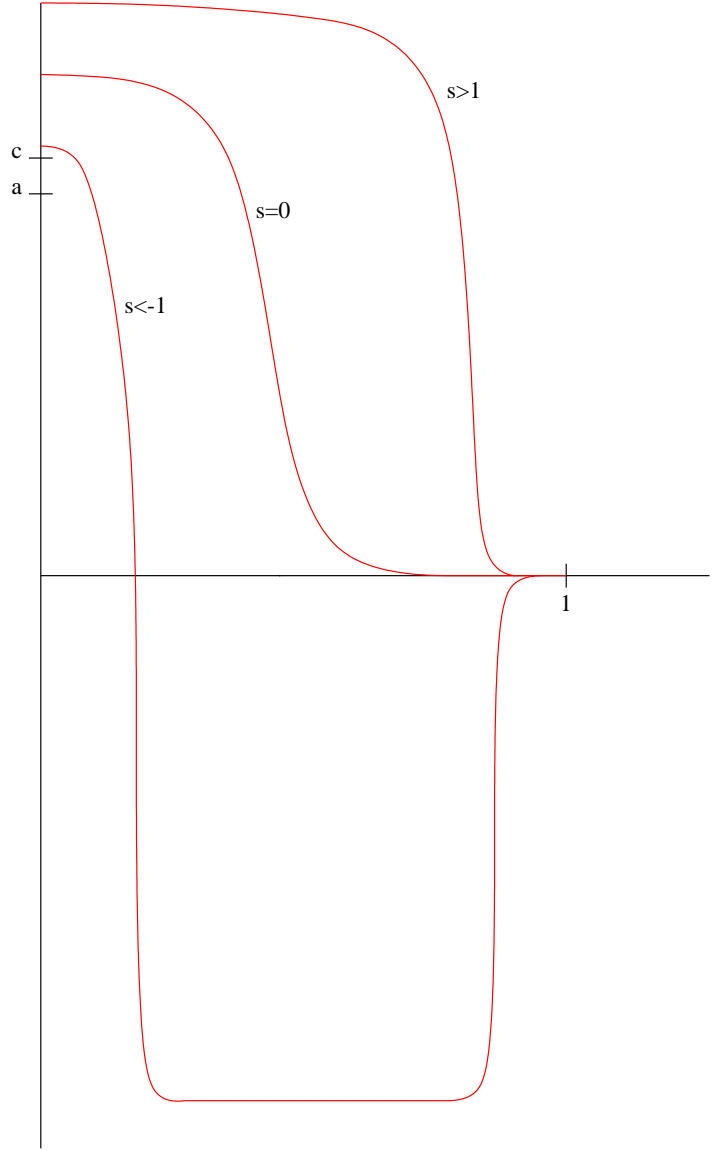}}
\caption{{A family of profile functions.}}
\label{fig:fig1}
\end{figure}

Fix a real number $c>0$ and choose a smooth family of real functions
$\{f_s(r)\}_{s \in \mathbb{R}}$, defined for $r\in \mathbb{R}$, with
the following properties (see Figure~\ref{fig:fig1}):
\begin{enumerate}
  \item[(i)] $f_s(-r)=f_s(r)$ for all $s$ and $r$.
  \item[(ii)] For every $s\in\R$
   $$
   f_s(0)> c,\qquad f_s''(0)<0,\qquad
   \lim_{s\to-\infty}f_s(0)=c,\qquad \lim_{s\to\infty}f_s(0)=\infty.
   $$
  \item[(iii)] For all $s$ and $r$ we have $\p_sf_s(r)\ge 0$.
  \item[(iv)] $f_s'(r)\le 0$ for $r\geq 0$ and, for $s\ge 1$,
   $$
   f_s(r) =
   \begin{cases}
      f_s(0)(1-r^2), &\mbox{if }0\le r\le 1-1/4s, \\
      0, &\mbox{if }r\ge 1-1/8s.
   \end{cases}
   $$
  \item[(v)] $f_s'(r)\le 0$ for $r\le1/2$, $f_s'(r)\ge 0$ for
   $r\ge1/2$, and, for $s\le -1$,
   $$
   f_s(r) =
   \begin{cases}
      f_s(0)(1-r^2), &\mbox{if }0\le r\le 1/8|s|, \\
      s, &\mbox{if }1/4|s|\le r\le 1-1/4|s|, \\
      0, &\mbox{if }r\ge 1-1/8|s|.
   \end{cases}
   $$
  \item[(vi)] For every $s$ the only critical point $r$ of $f_s$
   with $f_s(r)>0$ is $r = 0$.
\end{enumerate}
It is not hard to prove that such a family $\{f_s(r)\}_{s\in\R}$
indeed exists.

Now define $H_s:U^*X\to\R$ by
$$
H_s(q,p) := f_s(|p|).
$$
Since $X$ satisfies~$(T)$ or~$(N)$, every contractible closed
geodesic is constant.  Moreover, the constant periodic orbits
$x\in\Pp(H_s)$ have the form $x\equiv(q,p)$ where $f_s'(|p|)=0$, and
the symplectic action of such a constant solution is $
\Aa_{H_s}(x)=f_s(|p|).  $ Hence, by~(vi), the contractible periodic
solutions of $H_s$ with positive action have the form $x=(q,0)$ with
$$
\Aa_{H_s}(x)=f_s(0)>c.
$$
By Lemma~\ref{le:zero}, these solutions form a Morse--Bott manifold
of periodic orbits for $H_s$.  Hence, by Theorem~\ref{T:global}, we
have
$$\HF^{[a,\infty)}(H_s;0) \cong
\begin{cases}
   H_*(X;\Z_2).&\mbox{if }0<a < f_s(0),\\
   0,&\mbox{if }0<f_s(0)<a,
\end{cases}
$$
for every $s\in\R$. By Proposition~\ref{P:homotopy}, the monotone
homomorphism
$$
\sigma_{H_{s_1}H_{s_0}}: \HF^{[a,\infty)}(H_{s_0};\alpha)
\to\HF^{[a,\infty)}(H_{s_1};\alpha)
$$
is an isomorphism whenever $s_1\le s_0$ and
$a\notin[f_{s_1}(0),f_{s_0}(0)]$.  Now, for every $H\in\Hh(U^*X)$
there exists an $s\in\R$ such that $H\le H_s$.  Hence, by
Lemma~\ref{L:exhaust}~(ii), the homomorphism
$$
\pi_s:\SHi^{[a,\infty)}(U^*X;0)\to\HF^{[a,\infty)}(H_s;0)
$$
is an isomorphism for every $s\in\R$ such that $f_s(0)>a$.  Hence
$$
\SHi^{[a,\infty)}(U^*X;0) \cong H_*(X;\Z_2).
$$
Moreover, for every $H\in\Hh^{a,\infty}_c(U^*X, X;0)$ there exists
an $s\in\R$ such that $H_s\le H$.  Hence, by
Lemma~\ref{L:exhaust}~(i), the homomorphism
$$
\iota_s:\HF^{[a,\infty)}(H_s;0)\to\SHd^{[a,\infty);c}(U^*X, X;0)
$$
is an isomorphism for every $s\in\R$ in the case $a\le c$, and for
every $s\in\R$ with $f_s(0)<a$ in the case $a>c$.  Hence
$$\SHd^{[a,\infty);c}(U^*X, X;0) \cong
\begin{cases}
   H_*(X;\Z_2),&\mbox{if }a\le c, \\
   0,&\mbox{if }a>c.
\end{cases}
$$
In the case $a>c$ it follows that $T^{[a,\infty);c}_0=0$.  In the
case $a\le c$ it follows from Proposition~\ref{P:factor}, that the map
$T^{[a,\infty);c}_0$ can be expressed as the composition
$$
T^{[a,\infty);c}_0 = \iota_s\circ\pi_s:
\SHi^{[a,\infty)}(U^*X;0)\to\SHd^{[a,\infty);c}(U^*X, X ;0)
$$
for every $s\in\R$.  Hence, in this case, $T^{[a,\infty);c}_0$ is
an isomorphism.

It remains to prove the statement on $\widehat{C}(U^*X,X;0,a)$.
Indeed, it follows from what we have proved above that
$\A_c(U^*X,X,;0)=(0,c]$ for every $c>0$. Therefore, for every $a \in
\mathbb{R}$ we have:
$$\widehat{C}(U^*X,X;0,a)=\inf\left\{c > 0\,\bigm|\, \sup
   \A_c(U^*X,X;\alpha) > a \right\} = \max \{0,a\}.$$
The proof of
Theorem~\ref{T:zero} is complete.  \Qed


\subsection{Proof of Theorem~\ref{T:TN}}
\label{sb:TN}

\begin{figure}[htp]
   \centerline{\psfig{figure=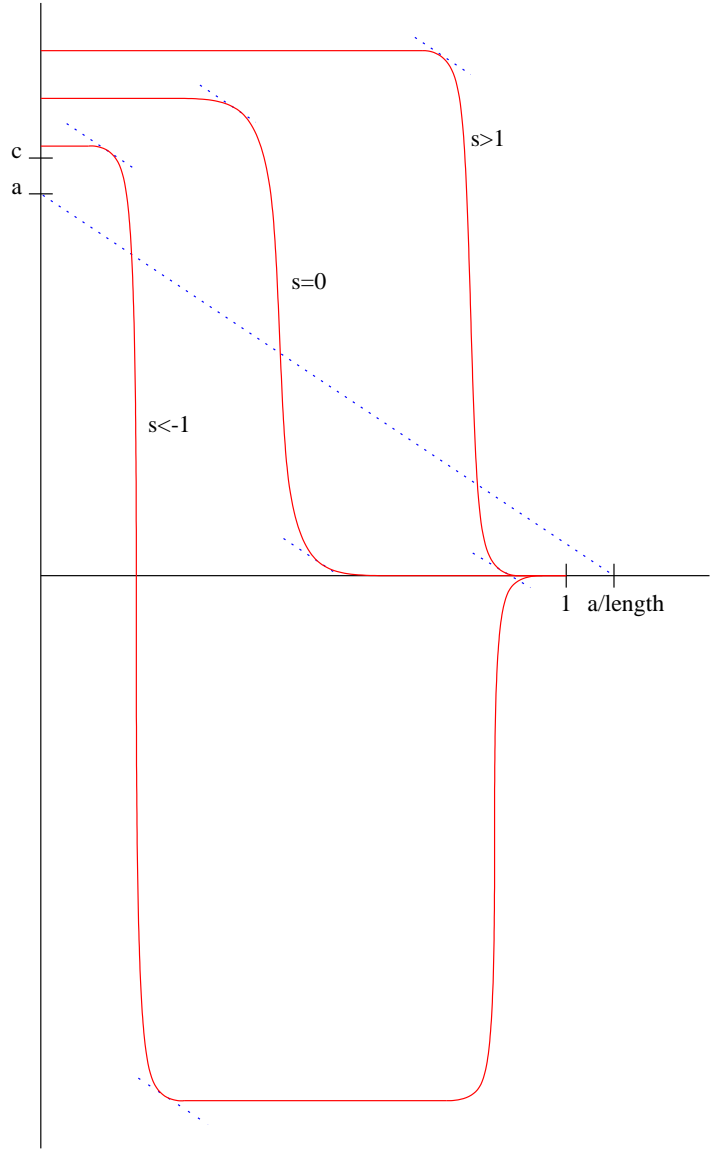}}
   \caption{{Another family of profile functions.}}
   \label{fig:fig2}
\end{figure}

Fix a nontrivial homotopy class $\alpha\in\widetilde{\pi}_1(X)$ and
let $\ell$ denote the length of the geodesics in this class. Moreover,
fix a real number $c>0$ and choose a smooth family of real functions
$\{f_s(r)\}_{s \in \mathbb{R}}$, defined for $r\in \mathbb{R}$, with
the following properties (see Figure~\ref{fig:fig2}):
\begin{enumerate}
  \item[(i)] $f_s(-r)=f_s(r)$ for all $s$ and $r$.
  \item[(ii)] For every $s\in\R$
   $$
   f_s(0)> c,\qquad \lim_{s\to-\infty}f_s(0)=c,\qquad
   \lim_{s\to\infty}f_s(0)=\infty.
   $$
  \item[(iii)] For all $s$ and $r$ we have $\p_sf_s(r)\ge 0$.
  \item[(iv)] If $s\ge 1$ then
   $$f_s(r) =
   \begin{cases}
      f_s(0), &\mbox{if }0\le r\le 1-3/8s, \\
      0, &\mbox{if }r\ge 1-1/8s,
   \end{cases}
   $$
   $f_s'(r)\le 0$ for all $r\geq 0$, and
   $$f_s''(r) =
   \begin{cases}
      <0, &\mbox{if }1-3/8s<r<1-2/8s, \\
      >0, &\mbox{if }1-2/8s<r<1-1/8s,
   \end{cases}
   $$
  \item[(v)] If $s\le -1$ then
   $$f_s(r) =
   \begin{cases}
      f_s(0), &\mbox{if }0\le r\le 1/8|s|, \\
      s, &\mbox{if }3/8|s|\le r\le 1-3/8|s|, \\
      0, &\mbox{if }r\ge 1-1/8|s|,
   \end{cases}
   $$
   $f_s'(r)\le 0$ for $r\le1/2$, $f_s'(r)\ge 0$ for $r\ge 1/2$, and
   $$f_s''(r) =
   \begin{cases}
      <0, &\mbox{if }1/8|s|<r<2/8|s|, \\
      >0, &\mbox{if }2/8|s|<r<3/8|s|.
   \end{cases}
   $$
  \item[(vi)] For every $s\in\R$ such that $f_s(0)>\ell$ there exist
   real numbers $r'_s>r_s>0$ such that
   $$
   f_s'(r_s) = f_s'(r_s') = -\ell,\qquad f_s''(r_s)<0,\qquad
   f_s''(r_s')>0,
   $$
   and $f_s'(r)\ne-\ell$ for every
   $r\in[0,\infty)\setminus\{r_s,r_s'\}$.
  \item[(vii)] For every $s \in \mathbb{R}$, the only possible points
   $r>0$ with $f_s'(r)=\ell$ must satisfy $f_s(r)<0$.
\end{enumerate}
It is not hard to show that such a family of functions
$\{f_s(r)\}_{s\in\R}$ indeed exists.
Now define $H_s:U^*X\to\R$ by
$$
H_s(q,p) := f_s(|p|).
$$
Consider first periodic orbits of $H_s$ that belong to one of the
sets $P^+(r,\alpha)$, $(r>0)$, as defined by~\eqref{eq:P}. We claim
that the corresponding action is negative. Indeed, at such a value of
$r>0$ we have $f_s'(r)=\ell$, and the action is $f_s(r)-\ell r$ which
is negative due to~(vii).

Next, denote by $P_s:=P^-(r_s,\alpha)$ and $P_s':=P^-(r_s',\alpha)$
the other two components of the set of periodic solutions in the class
$\alpha$ as defined by~(\ref{eq:P}). Then $P_s$ and $P_s'$ are both
diffeomorphic to $\T^n$ in the case~$(T)$ and to $S^1$ in the
case~$(N)$.  Moreover, by~(vi) and Lemma~\ref{le:TN}, they are
Morse--Bott manifolds of periodic orbits for $H_s$ for every $s\in\R$
and the values of the symplectic action functional on these two
critical manifolds are
$$
\Aa_{H_s}(P_s) = f_s(r_s) + r_s\ell,\qquad \Aa_{H_s}(P_s') =
f_s(r_s') + r_s'\ell.
$$
Fix a real number $a$ and denote $P:=\T^n$ in the case~$(T)$ and
$P:=S^1$ in the case~$(N)$.  We prove Theorem~\ref{T:TN} in five
steps.

\medskip
\noindent{\bf Step~1.} 
{\it If $a<\ell$ then $ \SHi^{[a,\infty)}(U^*X;\alpha) = 0.  $}

\medskip
\noindent
By~(iv), $f_s'(r) \leq 0$ for every $s\geq 1, r>0$ hence for $s\geq 1$
there are no periodic orbits of the type $P^+(r,\alpha)$.  Thus for
$s\geq 1$ the only families of periodic orbits are $P_s$ and $P_s'$.
Since both $r_s$ and $r_s'$ converge to $1$ as $s\to\infty$ it follows
that the sets $P_s$ and $P_s'$ both have action bigger than $a$ for
$s$ sufficiently large.  Hence
$$
\HF^{[a,\infty)}(H_s;\alpha)
\cong\HF^{[-\infty,\infty)}(H_s;\alpha) = 0
$$
for $s$ sufficiently large.  The last equation holds because
$\alpha\ne 0$, so $\HF^{[-\infty,\infty)}(H;\alpha)$ is independent of
$H$, and there is a Hamiltonian function with only contractible
$1$-periodic orbits.  Now Step~1 follows from
Lemma~\ref{L:exhaust}~(ii).

\medskip
\noindent{\bf Step~2.} 
{\it If $a\ge \ell$ then $ \SHi^{[a,\infty)}(U^*X;\alpha) \cong
  H_*(P;\Z_2).  $ Moreover, the homomorphism
  $$
  \pi_s:\SHi^{[a,\infty)}(U^*X;\alpha) \to
  \HF^{[a,\infty)}(H_s;\alpha)
  $$
  is an isomorphism whenever $f_s(0)>a$.  }

\medskip
\noindent
As $a > \ell > 0$ we may ignore all periodic orbits of the type
$P^+(r,\alpha)$ and consider only the families $P_s, P'_s$. The
numbers $r_s$ and $r_s'$ are the critical points of the function
$f_{s,\ell}:[0,1]\to\R$ given by
\begin{equation}\label{eq:fsl}
   f_{s,\ell}(r):=f_s(r)+r\ell.  
\end{equation}
By~(vi), the point $r_s$ is a strict local maximum and the point
$r_s'>r_s$ is a strict local minimum.  Suppose that $f_s(0)>a$.  Then
$f_{s,\ell}(0)=f_s(0)>a$ and $f_{s,\ell}(1)=\ell\le a$, hence it
follows that $f_{s,\ell}(r_s)>a$ and $f_{s,\ell}(r_s')<a$. This means
that
$$
\Aa_{H_s}(P_s)>a,\qquad \Aa_{H_s}(P_s')<a.
$$
Hence, by Theorem~\ref{T:global}, $\HF^{[a,\infty)}(H_s;\alpha)
\cong H_*(P;\Z_2)$ and, by Proposition~\ref{P:homotopy}, the
mo\-no\-tone homomorphism $ \sigma_{H_{s_1}H_{s_0}}:
\HF^{[a,\infty)}(H_{s_0};\alpha) \to\HF^{[a,\infty)}(H_{s_1};\alpha) $
is an isomorphism whenever $f_{s_i}(0)>a$ for $i=0,1$ and $s_1\le s_0$
and Step~2 follows from Lemma~\ref{L:exhaust}~(ii).

\medskip
\noindent{\bf Step~3.}
{\it If $a > c >0$ then $ \SHd^{[a,\infty);c}(U^*X,X;\alpha) = 0.  $}

\medskip
\noindent
As $a>0$ we can again ignore all orbits of type $P^+(r,\alpha)$.
Since both $r_s$ and $r_s'$ converge to $0$ as $s\to-\infty$ it
follows that the sets $P_s$ and $P_s'$ both have action less than $a$
for $-s$ sufficiently large.  Hence $ \HF^{[a,\infty)}(H_s;\alpha) = 0
$ for $-s$ sufficiently large.  Hence Step~3 follows from
Lemma~\ref{L:exhaust}~(i).

\medskip
\noindent{\bf Step~4.}
{\it If $0<a\le c$ then $ \SHd^{[a,\infty);c}(U^*X,X;\alpha) \cong
  H_*(P;\Z_2).  $ Moreover, the homomorphism
  $$
  \iota_s:\HF^{[a,\infty)}(H_s;\alpha) \to
  \SHd^{[a,\infty);c}(U^*X,X;\alpha)
  $$
  is an isomorphism for $s \ll -1$.  }

\medskip
\noindent
Since $a>0$, we may ignore as in previous steps orbits of type
$P^+(r,\alpha)$. Let $f_{s,\ell}:[0,1]\to\R$ be given
by~(\ref{eq:fsl}).  Then, by~(ii), $f_{s,\ell}(0)=f_s(0)>c\ge a$ and
hence
$$
\Aa_{H_s}(P_s) = f_{s,\ell}(r_s) > f_{s,\ell}(0) > a.
$$
If $s<\min \{-1, a-\ell/2\}$ then $f_{s,\ell}(1/2)=s+\ell/2<a$ and
hence
$$
\Aa_{H_s}(P_s') = f_{s,\ell}(r_s') < a.
$$
By Theorem~\ref{T:global}, $ \HF^{[a,\infty)}(H_s;\alpha) \cong
H_*(P_s;\Z_2) $ for $s<\min\{-1, a-\ell/2\}$.  By
Proposition~\ref{P:homotopy}, the monotone homomorphism $
\sigma_{H_{s_1}H_{s_0}}: \HF^{[a,\infty)}(H_{s_0};0)
\to\HF^{[a,\infty)}(H_{s_1};0) $ is an isomorphism for
$s_1<s_0<\min\{-1, a-\ell/2\}$. Step~4 follows now from
Lemma~\ref{L:exhaust}~(i).

\medskip
\noindent{\bf Step~5.}
{\it If $\ell\le a\le c$ then the homomorphism
  $$
  T^{[a,\infty);c}_\alpha: \SHi^{[a,\infty)}(U^*X;\alpha) \to
  \SHd^{[a,\infty);c}(U^*X,X;\alpha)
  $$
  is an isomorphism.}

\medskip
\noindent
By~(ii), $f_s(0)>c\ge a$ for every $s$.  Hence, by Step~2, $\pi_s$ is
an isomorphism for every $s\in\R$.  Moreover, by Step~4, $\iota_s$ is
an isomorphism for $s \ll -1$. By Proposition~\ref{P:factor},
$T^{[a,\infty);c}_\alpha=\iota_s\circ\pi_s$ for every $s$.  Hence
$T^{[a,\infty);c}_\alpha$ is an isomorphism.

It remains to prove the statement on $\widehat{C}(U^*X,X;\alpha,a)$.
Indeed, it follows from what we have proved above that
$\A_c(U^*X,X,;\alpha)=[\ell,c]$ for every $c>0$. Therefore, for every
$a \in \mathbb{R}$ we have:
$$\widehat{C}(U^*X,X;\alpha,a)=\inf\left\{c > 0\,\bigm|\, \sup
   \A_c(U^*X,X;\alpha) > a \right\} = \max \{\ell,a\}.$$
The proof of
Theorem~\ref{T:TN} is complete.  \Qed

\end{document}